\documentclass[11pt, a4paper]{amsart}
\usepackage[T1]{fontenc}
\usepackage[utf8]{inputenc}
\usepackage{geometry}                
\usepackage{graphicx}
\usepackage{amsmath, amsthm, amssymb, amscd}
\usepackage{bbm}
\usepackage{mathrsfs}
\usepackage{hyperref}
\usepackage[all]{xy}

\usepackage{stmaryrd}
\usepackage{float}
\usepackage{subfig}
\usepackage{tikz}

\theoremstyle{plain}
\newtheorem{thm}{Theorem}[section]
\newtheorem{lem}[thm]{Lemma}
\newtheorem{prop}[thm]{Proposition}
\newtheorem{cor}[thm]{Corollary}

\theoremstyle{definition}
\newtheorem{defn}{Definition}[section]
\newtheorem{conj}{Conjecture}[section]

\theoremstyle{remark}
\newtheorem{rem}{Remark}[section]

\numberwithin{equation}{section}

\newcommand{\gl}{\mathrm{GL}}
\newcommand{\ggl}{\mathfrak{gl}}
\newcommand{\pgl}{\mathrm{PGL}}
\newcommand{\SL}{\mathrm{SL}}

\newcommand{\bbg}{\mathbb{G}}

\newcommand{\cf}{\mathcal{F}}
\newcommand{\co}{\mathcal{O}}
\newcommand{\cp}{\mathcal{P}}

\newcommand{\cl}{\mathcal{L}}

\newcommand{\xx}{\mathscr{X}}

\newcommand{\ka}{\mathfrak{a}}
\newcommand{\kg}{\mathfrak{g}}
\newcommand{\kn}{\mathfrak{n}}
\newcommand{\kp}{\mathfrak{p}}
\newcommand{\kt}{\mathfrak{t}}

\newcommand{\km}{\mathfrak{m}}
\newcommand{\ko}{\mathfrak{o}}
\newcommand{\kl}{\mathfrak{l}}

\newcommand{\bn}{\mathbf{N}}
\newcommand{\bz}{\mathbf{Z}}
\newcommand{\br}{\mathbf{R}}

\newcommand{\bq}{\mathbf{Q}}

\newcommand{\bc}{\mathbf{C}}

\newcommand{\baa}{\mathbf{a}}

\newcommand{\rv}{\mathrm{v}}
\newcommand{\rw}{\mathrm{w}}
\newcommand{\rc}{\mathrm{c}}
\newcommand{\rd}{\mathrm{d}}
\newcommand{\rr}{\mathrm{r}}
\newcommand{\rs}{\mathrm{s}}

\newcommand{\val}{\mathrm{val}}
\newcommand{\lie}{\mathrm{Lie}}
\newcommand{\Ad}{\mathrm{Ad}}
\newcommand{\ad}{\mathrm{ad}}

\newcommand{\diag}{\mathrm{diag}}
\newcommand{\gal}{\mathrm{Gal}}

\newcommand{\sch}{\mathrm{Sch}}

\newcommand{\reg}{\mathrm{reg}}

\newcommand{\Hom}{\mathrm{Hom}}
\newcommand{\ec}{\mathrm{Ec}}

\newcommand{\rk}{\mathrm{rk}}

\newcommand{\ep}{\epsilon}
\newcommand{\vol}{\mathrm{vol}}
\newcommand{\res}{\mathrm{Res}}

\newcommand{\wtt}{\widetilde{T}}
\newcommand{\qlbar}{\overline{\mathbf{Q}}_{\ell}}
\newcommand{\fq}{\mathbf{F}_{q}}

\author{Zongbin \textsc{Chen}}

\address{Jinchunyuan west building, office 255\\
Yau Mathematical Science center, Tsinghua University\\
100084, Haidian district\\
Beijing, China} 
\email{zbchen@mail.tsinghua.edu.cn} 
\title{Truncated affine Springer fibers and Arthur's weighted orbital integrals}

\begin{document}
\maketitle


\begin{abstract}

We explain an algorithm to calculate Arthur's weighted orbital integral in terms of the number of rational points on the fundamental domain of the associated affine Springer fiber. The strategy is to count the number of rational points of the truncated affine Springer fibers in two ways: by the Arthur-Kottwitz reduction and by the Harder-Narasimhan reduction. A comparison of results obtained from these two approaches gives us recurrence relations between the number of rational points on the fundamental domains of the affine Springer fibers and Arthur's weighted orbital integrals. As an example, we calculate Arthur's weighted orbital integrals for the group $\gl_{2}$ and $\gl_{3}$.

\end{abstract}

\section{Introduction}

Let $\fq$ be the finite field with $q$ elements. Let $F=\fq(\!(\ep)\!)$ be the field of Laurent series with coefficients in $\fq$, $\co=\fq[\![\ep]\!]$ the ring of integers of $F$, $\kp=\ep \fq[\![\ep]\!]$ the maximal ideal of $\co$. We fix an algebraic closure $\overline{\mathbf{F}}_{q}$ of $\fq$, and also a compatible separable algebraic closure $\overline{F}$ of $F$. Let $\val:\overline{F}^{\times}\to \bq$ be the discrete valuation normalized by $\val(\ep)=1$.

Let $G$ be a connected split reductive algebraic group over $\fq$, assume that $\mathrm{char}(\fq)>|W|$, $W$ being the Weyl group of $G$. Let $G_{F}$ be the base change of $G$ to $F$. Let $T$ be a maximal torus of $G_{F}$. We make the assumption that the splitting field of $T$ is totally ramified over $F$.
Let $S\subset T$ be the maximal $F$-split subtorus of $T$, let $M_{0}=Z_{G_{F}}(S)$ be the centralizer of $S$ in $G_{F}$, then $M_{0}$ is a Levi subgroup of $G_{F}$ and $T$ is elliptic in $M_{0}$. Given an algebraic group, we use the Gothic letter to denote its Lie algebra.

Let $\gamma\in \kt(F)$ be a regular element, it is elliptic in $\km_{0}(F)$. Let $\cl(M_{0})$ be the set of Levi subgroups of $G$ containing $M_{0}$. For $M\in \cl(M_{0})$, consider Arthur's weighted orbital integral
\begin{equation}\label{arthur integral}
J_{M}(\gamma)=J_{M}(\gamma,\,\mathbbm{1}_{\kg(\co)})=
\int_{T(F)\backslash G(F)} \mathbbm{1}_{\kg(\co)}\big(\Ad(g)^{-1}\gamma\big)\rv_{M}(g)\frac{dg}{dt},
\end{equation}
where $\mathbbm{1}_{\kg(\co)}$ is the characteristic function of the lattice $\kg(\co)$ in $\kg(F)$, $\rv_{M}(g)$ is Arthur's weight factor, and $dg$ and $dt$ are Haar measures on $G(F)$ and $T(F)$ respectively. One of our main results states that it can be expressed in terms of the number of rational points of the \emph{fundamental domains} $F_{\gamma}^{L}$  of the affine Springer fiber $\xx_{\gamma}^{L}$, $L\in \cl(M_{0})$. The main idea is to count the number of rational points of the truncated affine Springer fibers in two different ways: by the \emph{Arthur-Kottwitz reduction} and by the \emph{Harder-Narasimhan reduction}.

Before entering into the details of our approach, we give examples of results that can be obtained in this way. The calculations for the group $G=\gl_{2}$ is easy, the results are summarized in theorem \ref{count gl2}, \ref{gl2 cal ramified}. But for the group $G=\gl_{3}$, the calculations are already quite non-trivial.
There are 3 cases to deal with: the element $\gamma$ can be split, mixed or elliptic. 
When $\gamma$ is split, we can find a set of simple roots $\{\alpha_{1},\alpha_{2}\}$ in the root system $\Phi(G,T)$ of $G$ with respect to $T$ such that 
$$
\val(\alpha_{1}(\gamma))=\val((\alpha_{1}+\alpha_{2})(\gamma))\leq \val(\alpha_{2}(\gamma)).
$$
We call $(n_{1},n_{2})=(\val(\alpha_{1}(\gamma)), \val(\alpha_{2}(\gamma)))$ the \emph{root valuation} of $\gamma$.

\begin{thm}
Let $G=\gl_{3}$, $T$ the maximal torus of diagonal matrices. Let $\gamma\in \kt(\co)$ be a regular element with root valuation $(n_{1},n_{2})\in \bn^{2}$, with $n_{1}\leq n_{2}$.
Up to an explicit volume factor, we have
\begin{eqnarray*}
J_{T}(\gamma)\doteq\sum_{i=1}^{n_{1}}i(q^{2i-1}+q^{2i-2})
+\sum_{i=n_{1}+n_{2}}^{2n_{1}+n_{2}-1}
(4n_{1}+2n_{2}-4i-3)q^{i}
+(n_{1}^{2}+2n_{1}n_{2})q^{2n_{1}+n_{2}}.
\end{eqnarray*}
For $\alpha\in \Phi(G,T)$, let $M_{\alpha}$ be the unique Levi subgroup containing $T$ with root system $\{\pm \alpha\}$, then up to an explicit volume factor, 
\begin{eqnarray*}
J_{M_{\alpha_{1}}}(\gamma)=J_{M_{\alpha_{1}+\alpha_{2}}}(\gamma)&\doteq&(n_{1}+n_{2})q^{2n_{1}+n_{2}}-q^{n_{1}+n_{2}}(1+q+\cdots+q^{n_{1}-1})\\
&&-q^{2n_{1}}(1+q+\cdots+q^{n_{2}-1}),
\end{eqnarray*}
and
$$
J_{M_{\alpha_{2}}}(\gamma)\doteq
2n_{1}q^{2n_{1}+n_{2}}-2q^{n_{1}+n_{2}}(1+q+\cdots+q^{n_{1}-1}).
$$

\end{thm}


When $\gamma$ is mixed, i.e. $T$ is isomorphic to $F^{\times}\times \mathrm{Res}_{E_{2}/F}E_{2}^{\times}$, where $E_{2}=\fq(\!(\ep^{\frac{1}{2}})\!)$ is the unique totally ramified extension of $F$ of degree 2, it can be conjugate to a matrix of the form
\begin{equation}\label{gl3 mixed}
\gamma=\begin{bmatrix}
a&&\\
&&b\\&b\ep&
\end{bmatrix}.
\end{equation}
Let $m=\val(a),\,n=\val(b)$, we have

\begin{thm}

Let $G=\gl_{3}$, let $\gamma$ be a matrix in the form $(\ref{gl3 mixed})$. When $\val(a)=m\leq n$, up to an explicit volume factor,  
\begin{eqnarray*}
J_{M_{0}}(\gamma)\doteq 2mq^{2m+n}+\sum_{j=m+n+1}^{2m+n-1}2(j-m-n)q^{j}
-\sum_{j=0}^{2m-1}\left(\bigg\lfloor \frac{j}{2} \bigg\rfloor+1\right) q^{j}.
\end{eqnarray*}
Similarly, when $\val(a)=m> n$, up to an explicit volume factor, 
\begin{eqnarray*}
J_{M_{0}}(\gamma)\doteq (2n+1)q^{3n+1}+\sum_{j=2n+1}^{3n}(2j-4n-1)q^{j}-\sum_{j=0}^{2n}\left(\bigg\lfloor \frac{j}{2} \bigg\rfloor+1\right) q^{j},
\end{eqnarray*}
where $\lfloor x\rfloor$ denotes the maximal integer less than or equal to $x$.   

\end{thm}

When $\gamma$ is anisotropic, Arthur's weighted orbital integral is just the orbital integral, and the result was essentially obtained by Goresky, Kottwitz and MacPherson \cite{gkm2}. See theorem \ref{gl3 ram count 1}, \ref{gl3 ram count 2} for the counting result.

Now let me explain our approach to the calculation of Arthur's weighted orbital integrals using the geometry of the affine Springer fibers. For simplicity, we restrict to  $J_{M_{0}}(\gamma)$. 
The affine Springer fiber $\xx_{\gamma}$ is the closed sub-scheme of the affine grassmannian $\xx=G(F)/G(\co)$ defined by the equation 
$$
\xx_{\gamma}=\big\{g\in G(F)/G(\co)\,\big|\,\Ad(g^{-1})\gamma \in \kg(\co)\big\}. 
$$ 
They can be used to geometrize Arthur's weighted orbital integrals. The group $T(F)$ acts on $\xx_{\gamma}$ by left translation. For $\mu\in X_{*}(S)$, we write $\ep^{\mu}$ for $\mu(\ep)\in S(F)$. The map $\mu\to \ep^{\mu}$ identifies $X_{*}(S)$ with a subgroup of $S(F)\subset T(F)$, which we denote by $\Lambda$. It acts freely on $\xx_{\gamma}$ and the quotient $\Lambda\backslash\xx_{\gamma}$ is a projective scheme of finite type over $\fq$ (see \cite{kl} \S3).
A simple reformulation shows that
$$
\int_{T(F)\backslash G(F)} \mathbbm{1}_{\kg(\co)}\big(\Ad(g)^{-1}\gamma\big)\rv_{M_{0}}(g)\frac{dg}{dt}=c\cdot\sum_{[g]\in\Lambda\backslash \xx_{\gamma}(\fq)} \rv_{M_{0}}(g), 
$$
where $[g]$ denotes the point $gG(\co)\in \xx$ and $c$ is a volume factor.

But this expression doesn't facilitate the calculations of Arthur's weighted orbital integral. We have to proceed in an indirect way. 
Let $\xi\in \ka_{M_{0}}^{G}$ be a generic element (for the definition of $\ka_{M_{0}}^{G}$, see the section of notations), Laumon and Chaudouard \cite{cl2} introduce a variant of the weighted orbital integral
\begin{equation}\label{cl variant}
J_{M_{0}}^{\xi}(\gamma)=J_{M_{0}}^{\xi}(\gamma,\,\mathbbm{1}_{\kg(\co)})=
\int_{T(F)\backslash G(F)} \mathbbm{1}_{\kg(\co)}\big(\Ad(g)^{-1}\gamma\big)\rw_{M_{0}}^{\xi}(g)\frac{dg}{dt},
\end{equation}
with a slightly different weight factor $\rw_{M_{0}}^{\xi}(g)$. The two weight factors are closely related to each other. When $G$ is semisimple, Laumon and Chaudouard show that 
$$
J_{M_{0}}(\gamma)=\mathrm{vol}(\ka_{M_{0}}/X_{*}(M_{0}))\cdot J_{M_{0}}^{\xi}(\gamma). 
$$

The variant $J_{M_{0}}^{\xi}(\gamma)$ has a better geometric interpretation. In fact, we can introduce a notion of $\xi$-stability on the affine Springer fiber $\xx_{\gamma}$, and show that
$$
J_{M_{0}}^{\xi}(\gamma)=\vol_{dt}\big(T(F)^{1}\big)^{-1}\cdot|\xx_{\gamma}^{\xi}(\fq)|. 
$$ 
The advantage of this variant is clear: it is a plain count rather than a weighted count. Moreover, we can use the Harder-Narasimhan reduction to get $|\xx_{\gamma}^{\xi}(\fq)|$ recursively from $|\xx_{\gamma}(\fq)|$, if only the latter is finite. Unfortunately, this is not the case, as can be seen from the fact that the free abelian group $\Lambda$ acts freely on $\xx_{\gamma}$.

Let $\Pi$ be a positive $(G,M_{0})$-orthogonal family, we can introduce a truncation $\xx_{\gamma}(\Pi)$ to overcome the finiteness issue. When $\Pi$ is sufficiently regular, we can reduce the calculation of the rational points on $\xx_{\gamma}(\Pi)$ to that of the fundamental domains $F_{\gamma}^{L}$, by the \emph{Arthur-Kottwitz reduction}. Recall that the fundamental domain $F_{\gamma}$ is introduced in \cite{chen2} to play the role of an irreducible component of $\xx_{\gamma}$ (All the irreducible components of $\xx_{\gamma}$ are isomorphic because $T(F)$ acts transitively on a dense open sub-scheme of it). The Arthur-Kottwitz reduction is a construction that decomposes $\xx_{\gamma}(\Pi)$ into locally closed sub-schemes, which are iterated affine fibrations over the fundamental domains $F_{\gamma}^{L}$, $L\in \cl(M_{0})$. The counting result is summarized in corollary \ref{AK count}. In particular, it shows that $\xx_{\gamma}(\Pi)$ depends quasi-polynomially on the truncation parameter.

On the other hand, the Harder-Narasimhan reduction doesn't behave well on $\xx_{\gamma}(\Pi)$. In fact, near the boundary, the Harder-Narasimhan strata are generally not affine fibrations over truncations of $\xx_{\gamma}^{L,\,\xi^{L}}$. To overcome this difficulty, we cut $\xx_{\gamma}(\Pi)$ into two parts: \emph{the tail} and \emph{the main body}. Roughly speaking, the tail is the union of the ``boundary irreducible components'' of $\xx_{\gamma}(\Pi)$, and the main body is its complement. The Harder-Narasimhan reduction works well on the main body, and we can use it to count the number of rational points. The result is summarized in theorem \ref{HN main body}, it can be expressed in terms of $|\xx_{\gamma}^{L,\,\xi^{L}}(\fq)|$, $L\in \cl(M_{0})$. The counting points on the tail proceeds by the Arthur-Kottwitz reduction, and can be expressed in terms of $|F_{\gamma}^{L}(\fq)|$'s. But we are not able to obtain an explicit expression, we get a recursion.

These two different approaches to counting rational points on $\xx_{\gamma}(\Pi)$ give us a recursive equation that involves the $|F_{\gamma}^{L}(\fq)|$'s and the $|\xx_{\gamma}^{L,\,\xi^{L}}(\fq)|$'s. Solving it, we can express the latter ones in terms of the former ones. The problem of calculating $J_{M_{0}}(\gamma)$ is thus reduced to counting points on $F_{\gamma}$. 

The geometry of $F_{\gamma}$ is simpler than that of $\xx_{\gamma}^{\xi}$: Goresky, Kottwitz and MacPhersion \cite{gkm1} have conjectured that the cohomology of $\xx_{\gamma}$ is pure in the sense of Deligne. As we have shown in \cite{chen2}, this is equivalent to the cohomological purity of $F_{\gamma}$. In fact, it is even expected that $F_{\gamma}$ admits a \emph{Hessenberg paving} (This notion was introduced by Goresky, Kottwitz and MacPhersion \cite{gkm2}). On the contrary, $\xx_{\gamma}^{\xi}$ is generally not cohomologically pure, as one can see in case $G=\SL_{2}$, or from the appearance of minus sign in the counting points result of theorem \ref{integral gl3}, \ref{beta}, \ref{alpha}. Although one can still look at the quotient $\xx_{\gamma}^{\xi}/A_{M_{0}}$, $A_{M_{0}}$ being the maximal $F$-split torus of the center of $M_{0}$, it is clear that the quotient no longer admits a torus action, hence it has much less structure to explore than $F_{\gamma}$.

When the torus $T$ splits, we \cite{chen3} make a conjecture on the Poincar\'e polynomial of $F_{\gamma}$, assuming the cohomological purity of $F_{\gamma}$. This gives a conjectural expression for $|F_{\gamma}(\fq)|$. We reproduce it here for the convenience of the reader. Following Chaudouard and Laumon \cite{cl1}, under the purity assumption, the cohomology of $F_{\gamma}$ can be expressed in terms of its $1$-skeleton under the $T$-action. 
Indeed, the $T$-equivariant cohomology $H_{T}^{*}(F_{\gamma},\qlbar)$\footnote{Here we actually mean the \emph{geometric} $H^{*}_{T_{\overline{\mathbf{F}}_{q}}}\big(F_{\gamma,\overline{\mathbf{F}}_{q}},\qlbar\big)$, we don't specify the base change to $\overline{\mathbf{F}}_{q}$ to simplify the notation. Similar convention for the other  cohomology groups.} will then be a free $H_{T}^{*}(\mathrm{pt}, \qlbar)$-algebra, and we have
\begin{equation}\label{eqh1}
H^{*}(F_{\gamma},\,\qlbar)=H^{*}_{T}(F_{\gamma},\,\qlbar)\otimes_{H^{*}_{T}(\mathrm{pt},\,\qlbar)}\qlbar.
\end{equation}
The torus $T$ acts on $F_{\gamma}$ with finitely many fixed points, but the $1$-dimensional $T$-orbits form a higher dimensional variety which we denote by $F_{\gamma}^{T,1}$. The bigger torus $\wtt=T\times \bbg_{m}$, $\bbg_{m}$ being the rotational torus, acts on $F_{\gamma}^{T,1}$ with finitely many $1$-dimensional $\wtt$-orbits, let $F_{\gamma}^{\wtt,1}$ be their union. Let $F_{\gamma}^{\wtt}$ be the set of $\wtt$-fixed points on $F_{\gamma}$. Let
\begin{equation}\label{eqh2}
H^{*}_{\widetilde{T}}(F_{\gamma},\,\qlbar):=H^{*}_{T}(F_{\gamma},\,\qlbar)\otimes_{H^{*}_{T}(\mathrm{pt},\,\qlbar)}H^{*}_{\widetilde{T}}(\mathrm{pt},\,\qlbar).
\end{equation}
Then the localization theorem of Goresky, Kottwitz and MacPheron \cite{gkm local} implies an exact sequence of equivariant cohomology
\begin{equation}\label{eqh3}
0\to H_{\wtt}^{*}(F_{\gamma},\,\qlbar)\to H^{*}_{\wtt}(F_{\gamma}^{\wtt},\,\qlbar)\to H^{*}_{\wtt}(F_{\gamma}^{\wtt,1},F_{\gamma}^{\wtt};\,\qlbar).
\end{equation}
Let $\Gamma$ be the graph with vertices $F_{\gamma}^{\wtt}$ and with edges $F_{\gamma}^{\wtt,1}$. Two vertices are linked by an edge if and only if they lie on the closure of the corresponding $1$-dimensional $\wtt$-orbit. We call it the \emph{moment graph} of $F_{\gamma}$ with respect to the action of $\wtt$. The above result implies that the information about the cohomology of $F_{\gamma}$ is encoded in $\Gamma$. A direct calculation of the cohomology via equations (\ref{eqh1}), (\ref{eqh2}), (\ref{eqh3}) turns out to be very hard, and we look for a combinatorial way to get around it.

Let $\mathfrak{o}$ be a total order among the vertices of the graph $\Gamma$, it will serve as the paving order. We associate to it an \textit{acyclic} oriented graph $(\Gamma, \ko)$ such that the source of each arrow is greater than its target with respect to $\ko$. For $v\in \Gamma$, denote by $n^{\ko}_{v}$ the number of arrows having source $v$.

\begin{defn}
The formal Betti number $b^{\ko}_{2i}$ associated to the order $\ko$ is defined as
$$
b_{2i}^{\ko}=\sharp\{v\in \Gamma:\,n^{\ko}_{v}=i\}.
$$
We call
$$
P^{\ko}(t)=\sum_{i}b_{2i}^{\ko}t^{2i}
$$
the \emph{formal Poincaré polynomial} associated to the order $\ko$.
\end{defn}

\begin{defn}

For $P_{1}(t),\,P_{2}(t)\in \bz[t]$, we say that $P_{1}(t)< P_{2}(t)$ if the leading coefficient of $P_{2}(t)-P_{1}(t)$ is positive. 

\end{defn}

\begin{conj}\label{conj poincare}

Let $P(t)$ be the Poincaré polynomial of $F_{\gamma}$, then
$$
P(t)=\min_{\ko}\{P^{\ko}(t)\},
$$
where $\ko$ runs through all the total orders among the vertices of $\Gamma$.
\end{conj}

The conjecture can be thought of as a kind of \emph{Morse inequality}, it has been verified in a lot of examples. For the group $G=\gl_{2}$ and $\gamma=\diag(\gamma_{1},\gamma_{2})$ with $\gamma_{1},\gamma_{2}\in F$ and $\val(\gamma_{1}-\gamma_{2})=n\in \bn$, the moment graph of $F_{\gamma}$ contains $n+1$ vertices, which are pairwise connected by an edge. It is clear that the conjecture holds in this case. In general, the moment graph of $F_{\gamma}$ is easy to describe, and we have an algorithm to find an order $\ko$ which conjecturally should attain the minimum. Although we are not able to prove both conjectures at the moment, they have been very helpful for the construction of affine pavings of $F_{\gamma}$ in concrete examples.

Under the purity assumption, conjecture \ref{conj poincare} implies counting points result for $F_{\gamma}$. Indeed, the equations (\ref{eqh1}), (\ref{eqh2}), (\ref{eqh3}) are $\gal(\overline{\mathbf{F}}_{q}/\fq)$-equivariant, and the Frobenius endomorphism acts on $H^{2i}_{\wtt}(F_{\gamma}^{\wtt},\,\qlbar)$ by $q^{i}$ (the odd degree cohomologies vanish), hence it acts on $H^{*}(F_{\gamma},\,\qlbar)$ in the same way and so
$$
|F_{\gamma}(\fq)|=\min_{\ko}\{P^{\ko}(q^{1/2})\}.
$$
Together with the recurrence relation between $|\xx_{\gamma}^{\xi}(\fq)|$ and $|F_{\gamma}(\fq)|$, it gives a conjectural complete answer to the calculation of Arthur's weighted orbital integrals in the split case.


\subsection*{Notations}

We fix a split maximal torus $A$ of $G$ over $\fq$. Without loss of generality, we suppose that $A\subset M_{0}$. Let $\Phi=\Phi(G,A)$ be the root system of $G$ with respect to $A$, let $W$ be the Weyl group of $G$ with respect to $A$. For any subgroup $H$ of $G$ which is stable under the conjugation of $A$, we note $\Phi(H,A)$ for the roots appearing in $\lie(H)$. We fix a Borel subgroup $B_{0}$ of $G$ containing $A$. Let $\Delta=\{\alpha_{1},\cdots,\alpha_{r}\}$ be the set of simple roots with respect to $B_{0}$, let $\{\varpi_{i}\}_{i=1}^{r}$ be the corresponding fundamental weights. To an element $\alpha\in \Delta$, we have a unique maximal parabolic subgroup $P_{\alpha}$ of $G$ containing $B_{0}$ such that $\Phi(N_{P_{\alpha}},A)\cap \Delta=\{\alpha\}$, where $N_{P_{\alpha}}$ is the unipotent radical of $P_{\alpha}$. This gives a bijective correspondence between the simple roots in $\Delta$ and the maximal parabolic subgroups of $G$ containing $B_{0}$. Any semi-standard maximal parabolic subgroup $P$ of $G$ is conjugate to certain $P_{\alpha}$ by an element $w\in W$, the element $w\varpi_{\alpha}$ doesn't depend on the choice of $w$, we denote it by $\varpi_{P}$.

We use the $(G,M)$ notation of Arthur. Let $\cf(A)$ be the set of parabolic subgroups of $G$ containing $A$, let $\cl(A)$ be the set of Levi subgroups of $G$ containing $A$. For every $M\in \cl(A)$, we denote by $\cp(M)$ the set of parabolic subgroups of $G$ whose Levi factor is $M$, by $\cl(M)$ the set of Levi subgroups of $G$ containing $M$, and by $\cf(M)$ the set of parabolic subgroups of $G$ containing $M$. For $P\in \cp(M)$, we denote by $P^{-}\in \cp(M)$ the opposite of $P$ with respect to $M$.

Let $X^*(M)=\Hom(M, \bbg_m)$ and $X_{*}(M)=\Hom(X^{*}(M), \bz)$.
Let $\ka_M^{*}=X^*(M)\otimes\br$ and $\ka_M=X_{*}(M)\otimes\br$. The restriction $X^{*}(M)\to X^{*}(A)$ induces an injection $\ka_{M}^{*}\hookrightarrow \ka_{A}^{*}$. Let $(\ka_{A}^{M})^{*}$ be the subspace of $\ka_{A}^{*}$ generated by $\Phi(M,A)$. We have the decomposition in direct sums
$$
\ka_{A}^{*}=(\ka_{A}^{M})^{*}\oplus \ka_{M}^{*}.
$$

The canonical pairing 
$
X_{*}(A)\times X^{*}(A)\to \bz
$ 
can be extended bilinearly to $\ka_{A}\times \ka_{A}^{*}\to \br$. For $M\in \cl(A)$, we can embed $\ka_{M}$ in $\ka_{A}$ as the orthogonal subspace to $(\ka_{A}^{M})^{*}$. Let $\ka_{A}^{M}\subset \ka_{A}$ be the subspace orthogonal to $\ka_{M}^{*}$. We have the dual decomposition
$$
\ka_{A}=\ka_{M}\oplus \ka_{A}^{M},
$$
let $\pi_{M},\,\pi^{M}$ be the projections to the two factors. More generally, for $L,M\in \cf(A),\,M\subset L$, we also have a decomposition
$$
\ka_{M}^{G}=\ka_{L}^{G}\oplus \ka_{M}^{L}.
$$
Let $\pi_{M,L},\,\pi_{M}^{L}$ be the projections to the two factors. If the context is clear, we also simplify them to $\pi_{L},\,\pi^{L}$.

We identify $X_{*}(A)$ with $A(F)/A(\co)$ by sending $\chi$ to $\chi(\ep)$. With this identification, the canonical surjection $A(F)\to A(F)/A(\co)$ can be viewed as
\begin{equation}\label{indexT}
A(F)\to X_{*}(A).
\end{equation}

We use $\Lambda_{G}$ to denote the quotient of $X_{*}(A)$ by the coroot lattice of $G$ (the subgroup of $X_{*}(A)$ generated by the coroots of $A$ in $G$). It is independent of the choice of $A$, this is the algebraic fundamental group introduced by Borovoi \cite{bo}. According to Kottwitz \cite{k1}, we have a canonical homomorphism
\begin{equation}\label{indexM}
\nu_{G}: G(F)\to \Lambda_{G},
\end{equation}
which is characterized by the following properties: it is trivial on the image of $G_{\mathrm{sc}}(F)$ in $G(F)$ ($G_{\mathrm{sc}}$ is the simply connected cover of the derived group of $G$), and its restriction to $A(F)$ coincides with the composition of (\ref{indexT}) with the projection of $X_{*}(A)$ to $\Lambda_{G}$. Since the morphism (\ref{indexM}) is trivial on $G(\co)$, it descends to a map
$$
\nu_{G}:\xx\to \Lambda_{G},
$$
whose fibers are the connected components of $\xx$. For $\mu\in \Lambda_{G}$, we denote the connected component $\nu_{G}^{-1}(\mu)$ by $\xx^{\mu}$.

Finally, we suppose that $\gamma\in \kt(\co)$ satisfies $\gamma\equiv 0 \mod \ep$ to avoid unnecessary complications.

\subsection*{Acknowledgements}We want to thank G\'erard Laumon for the discussions which have led to this work, and we want to thank an anonymous referee for the careful reading and very helpful suggestions.

\section{(Weighted) orbital integrals and the affine Springer fibers}

We recall briefly the geometrization of the (weighted) orbital integrals using the affine Springer fibers. We fix a regular element $\gamma\in \kt(\co)$ as in the introduction. Let $P_{0}=M_{0}N_{0}$ be the unique element in $\cp(M_{0})$ which contains $B_{0}$.

\subsection{Orbital integrals}

We begin by fixing the Haar measures. Let $dg$ be the Haar measure on $G(F)$ normalized by the condition $\vol_{dg}\big(G(\co)\big)=1$. For the group $T(F)$, the definition is more involved as there is no natural $\co$-structure on $T$. Let $F^{\mathrm{ur}}=\overline{\mathbf{F}}_{q}(\!(\ep)\!)$, it is the completion of the maximal unramified extension of $F$. Let $\sigma$ be the Frobenius automorphism of both $\overline{\mathbf{F}}_{q}/\fq$ and $F^{\mathrm{ur}}/F$. We fix a separable algebraic closure $\overline{F^{\mathrm{ur}}}$ of $F^{\mathrm{ur}}$, let $I_{F}=\gal(\overline{F^{\mathrm{ur}}}/F^{\mathrm{ur}})$, it is the inertia subgroup of $\Gamma=\gal(\overline{F}/F)$. According to Kottwitz \cite{k2} \S 7.6, we have an exact sequence
\begin{equation}\label{kott 1}
1\to T(F^{\mathrm{ur}})_{1}\to T(F^{\mathrm{ur}})\xrightarrow{w_{T}}X_{*}(T)_{I_{F}}\to 1,
\end{equation}
which implies another exact sequence if we take the $\langle \sigma \rangle$-invariants, 
\begin{equation}\label{kott 2}
1\to T(F)_{1}\to T(F) \xrightarrow{w_{T}} (X_{*}(T)_{I_{F}})^{\langle \sigma \rangle}\to 1,
\end{equation}
with $T(F)_{1}:=T(F)\cap T(F^{\mathrm{ur}})_{1}$. We fix the Haar measure $dt$ on $T(F)$ by setting $\vol_{dt}\big(T(F)_{1}\big)=1$. 
The group $\Lambda$ is discrete and cocompact in $T(F)$. The volume of the quotient $\Lambda\backslash T(F)$ is calculated in \cite{gkm1} \S 15.3:
$$
\mathrm{vol}_{dt}\big(\Lambda\backslash T(F)\big)=\frac{|\mathrm{coker}[X_{*}(S)_{\Gamma}\to X_{*}(T)_{\Gamma}]|}{|\ker[X_{*}(S)_{\Gamma}\to X_{*}(T)_{\Gamma}]|}.
$$

Consider the orbital integral
\begin{equation}\label{def orb}
I^{G}_{\gamma}=\int_{T(F)\backslash G(F)} \mathbbm{1}_{\kg(\co)}\big(\Ad(g^{-1})\gamma\big)\frac{dg}{dt}. 
\end{equation}
It can be interpreted as counting points on the affine Springer fiber:  

\begin{prop}[Goresky, Kottwitz, MacPherson \cite{gkm1}]\label{geom orbital 1}
$$
I_{\gamma}^{G}=\frac{|\ker[X_{*}(S)_{\Gamma}\to X_{*}(T)_{\Gamma}]|}{|\mathrm{coker}[X_{*}(S)_{\Gamma}\to X_{*}(T)_{\Gamma}]|}\cdot|\Lambda\backslash \big(\xx_{\gamma}(\fq)\big)|.
$$
\end{prop}

The $T(F)$-action on $\xx_{\gamma}$ can be exploited to further simplify the computations. 
Let $\xx_{\gamma}^{\mathrm{reg}}$ be the open sub-scheme of $\xx_{\gamma}$ consisting of the points $[g]\in \xx_{\gamma}$ such that the image of $\Ad(g^{-1})\gamma$ under the reduction $\kg(\co)\to\kg$ is regular nilpotent. 

\begin{prop}[Bezrukavnikov \cite{b}]

The group $T(F)$ acts transitively on $\xx_{\gamma}^{\reg}$.

\end{prop}

\begin{prop}[Ng\^o \cite{ngo} Prop. 3.10.1]
The open sub-scheme $\xx_{\gamma}^{\reg}$ is dense in $\xx_{\gamma}$.
\end{prop}

Consequently, all the irreducible components of $\xx_{\gamma}$ are isomorphic to each other and they are parametrized by $\pi_{0}(T(F))$. In particular, all the connected components of $\xx_{\gamma}$ are isomorphic and they can be translated to each other under the $T(F)$-action. In the calculation of the orbital integral (\ref{def orb}), we can thus restrict to the central connected component of $\xx_{\gamma}$, often this simplifies calculations.

The calculation of $I_{\gamma}^{G}$ can be reduced to that of $I_{\gamma}^{M_{0}}$, it dates back at least to Harish-Chandra that
\begin{equation}\label{geom orbital reduced}
I_{\gamma}^{G}
=q^{\frac{1}{2}\val(\det(\ad(\gamma)|\kg_{F}/\km_{0,F}))}\cdot I_{\gamma}^{M_{0}}.
\end{equation}
Geometrically, this is a reflection of the existence of an affine fibration $f_{P}:\xx_{\gamma}\to \xx_{\gamma}^{M_{0}}$ for each $P\in \cp(M_{0})$.
Recall that for $Q=LN_{Q}\in \cf(A)$, we have the retraction
$$
f_{Q}:\xx\to \xx^{L}
$$ 
which sends $[g]=gK$ to $[h]:=hL(\co)$, where $g=nhk,\,n\in N_{Q}(F),\, h\in L(F),\, k\in K$ is the Iwasawa decomposition. We want to point out that the retraction $f_{Q}$ is not a morphism between ind-$\fq$-schemes, but its restriction to the inverse image of each connected component of $\xx^{L}$:
$$
f_{Q}:f_{Q}^{-1}(\xx^{L,\nu})\to \xx^{L,\nu},\quad \nu\in \Lambda_{L},
$$
is actually a morphism over $\fq$ between ind-$\fq$-schemes. Moreover, these retractions satisfy obvious transitivity property.

 
Restricted to the affine Springer fibers, the retraction $f_{Q}$ sends $\xx_{\gamma}$ to $\xx_{\gamma}^{L}$. To see this, for $[g]\in \xx_{\gamma}$, let $g=nhk$ be the Iwasawa decomposition as above. We can write $g=hn'k$ with $n'=h^{-1}nh\in N_{Q}(F)$. Now that $\Ad(h^{-1})\gamma\in \kl(F)$, we have 
$$
\Ad({n'}^{-1})\Ad(h^{-1})\gamma=\Ad(h^{-1})\gamma+n''
$$
for some $n''\in \kn_{Q}(F)$. This implies that
$$
\Ad(h)^{-1}\gamma \in [\kg(\co)+\kn_{Q,F}]\cap \kl(F)=\kl(\co),
$$
which means that $f_{Q}([g])=[h]\in \xx_{\gamma}^{L}$.

\begin{prop}[Kazhdan-Lusztig \cite{kl} \S5, Prop.1]\label{KL retraction}

For any $\nu\in \Lambda_{L}$, the retraction 
$$
f_{Q}:\xx_{\gamma}\cap f_{Q}^{-1}(\xx_{\gamma}^{L,\nu})\to \xx_{\gamma}^{L,\nu}
$$ 
is an iterated affine fibration over $\fq$ of relative dimension
$
\val(\det(\ad(\gamma)|\kn_{Q}(F))).
$

\end{prop}

The reader can also consult \cite{chen2}, Prop. 3.2 for a proof.

\subsection{Arthur's weighted orbital integral}

\subsubsection{The weight factor $\rv_{M}$}

Let $M\in \cl(M_{0})$, roughly speaking, the weight factor $\rv_{M}(g)$ is the volume of a polytope in $\ka_{M}$ generated by the point $[g]\in \xx$. Let $H_{M}: M(F)\to \ka_{M}$ be the unique map\footnote{Our definition differs from the conventional one by a minus sign. But as we will see, it simplifies  computations.} satisfying
$$
\chi(H_{M}(m))=\val(\chi(m)),\quad \forall\, \chi\in X^{*}(M),\,m\in M(F).
$$
Notice that it is a group homomorphism.
Moreover, it is invariant under the right $K$-action, so it induces a map from $\xx^{M}$ to $\ka_{M}$, still denoted by $H_{M}$. For $P=MN\in \cf(A)$, let $H_{P}:\xx\to \ka_{M}$ be the composition 
$$
H_{P}:\xx\xrightarrow{f_{P}}\xx^{M}\xrightarrow{H_{M}}\ka_{M},
$$
As shown in \cite{cl1}, Lemma 6.1, the map $H_{M}$ is constant on each connected component of $\xx^{M}$, so it has a factorization $
H_{M}:\xx^{M}\xrightarrow{\nu_{M}}\Lambda_{M}\to \ka_{M}$.
A simple calculation of the restriction of the map to $\xx^{A}\subset \xx^{M}$ shows that the map $\Lambda_{M}\to \ka_{M}$ is just the one induced from the natural inclusion $X_{*}(A)\hookrightarrow \ka_{A}= X_{*}(A)\otimes \br$. Hence $H_{P}$ is also the composition
$$
H_{P}:\xx\xrightarrow{f_{P}}\xx^{M}\xrightarrow{\nu_{M}}\Lambda_{M}\to \ka_{M}.
$$

The map $H_{P}$ has the following remarkable property.
There is a notion of adjacency among the parabolic subgroups in $\cp(M)$: Two parabolic subgroups $P_{1}=MN_{1},\,P_{2}=MN_{2}\in \cp(M)$ are said to be \emph{adjacent} if both of them are contained in a parabolic subgroup $Q=LN_{Q}$ such that $L\supset M$ and $\rk(L)=\rk(M)+1$. Given such an adjacent pair, we define an element $\beta_{P_{1},P_{2}}\in \Lambda_{M}$ in the following way: Consider the collection of elements in $\Lambda_{M}$ obtained from coroots of $A$ in $\kn_{1}\cap \kn_{2}^{-}$, we define $\beta_{P_{1},P_{2}}$ to be the minimal element in this collection, i.e. all the other elements are positive integral multiples of it. Note that $\beta_{P_{2},P_{1}}=-\beta_{P_{1},P_{2}}$, and if $M=A$, then $\beta_{P_{1},P_{2}}$ is the unique coroot which is positive for $P_{1}$ and negative for $P_{2}$. We denote also by $\beta_{P_{1},P_{2}}$ for its image in $\ka_{M}$ if no confusion is caused.

\begin{prop}[Arthur \cite{a} Lemma 3.6]\label{arthur orthogonal}

Let $P_{1},\,P_{2}\in \cp(M)$ be two adjacent parabolic subgroups. For any $x\in \xx$, we have
$$
H_{P_{1}}(x)-H_{P_{2}}(x)=n(x,P_{1},P_{2})\cdot \beta_{P_{1},P_{2}},
$$
with $n(x, P_{1}, P_{2})\in \bz_{\geq 0}$.

\end{prop}

The reader can consult \cite{chen2} Prop. 2.1 for a proof.
For any point $x\in \xx$, we write $\ec_{M}(x)$ for the convex hull in $\ka_{M}$ of the $H_{P}(x),\,P\in \cp(M)$. For any $Q\in \cf(M)$, we denote by $\ec_{M}^{Q}(x)$ the face of $\ec_{M}(x)$ whose vertices are $H_{P}(x),\,P\in \cp(M),\,P\subset Q$. When $M=A$, we omit the subscript $A$ to simplify the notation.

To define the volume, we need to choose a Lebesgue measure on $\ka_{M}^{G}$. We fix a $W$-invariant inner product $\langle\, \cdot\,,\cdot\,\rangle$ on the vector space $\ka_{A}^{G}$. Notice that $\ka_{A}^{M}$ and $\ka_{M}$ are orthogonal to each other with respect to the inner product for any $M\in \cl(A)$. We fix a Lebesgue measure on $\ka_{M}^{G}$ normalised by the condition that the lattice generated by the orthonormal bases in $\ka_{M}^{G}$ has covolume $1$.

The weight factor $\rv_{M}(g)$ is the volume of the projection $\pi_{M}^{G}\big(\ec_{M}(g)\big)\subset\ka_{M}^{G}$. We have to pass to $\ka_{M}^{G}$ because the polytope $\ec_{M}(g)$ will lie in a hyperplane of $\ka_{M}$ if $G$ has non-trivial connected center.
The weight factor $\rv_{M}(g)$ has the following invariance properties: It is invariant under the right action of $K$, i.e.
$$
\rv_{M}(gk)=\rv_{M}(g),\quad \forall\,k\in K.
$$
This is evident from the definition of $\rv_{M}(g)$. It is not so evident but also true that 
$$
\rv_{M}(mg)=\rv_{M}(g),\quad \forall\, m\in M(F).
$$
Indeed, for any $P\in \cp(M)$, we have
$
f_{P}(mg)=mf_{P}(g). 
$
As $H_{M}$ is a group homomorphism, this implies
$$
H_{P}(mg)=H_{M}(m)+H_{P}(g),
$$ 
so 
$
\ec_{M}(mg)$ is just the translation of $\ec_{M}(g)$ by $H_{M}(m).
$
In particular, they have the same volume.

Similar to proposition \ref{geom orbital 1}, we can interpret Arthur's weighted orbital integral as 
\begin{align*}
\int_{T(F)\backslash G(F)} \mathbbm{1}_{\kg(\co)}\big(\Ad(g)^{-1}\gamma\big)\rv_{M}(g)\frac{dg}{dt}=\frac{|\ker[X_{*}(S)_{\Gamma}\to X_{*}(T)_{\Gamma}]|}{|\mathrm{coker}[X_{*}(S)_{\Gamma}\to X_{*}(T)_{\Gamma}]|}
\sum_{[g]\in \Lambda\backslash \big(\xx_{\gamma}(\fq)\big)} \rv_{M}(g), 
\end{align*}
i.e. it is a weighted count of the rational points on the affine Springer fiber. Notice also that $J_{G}(\gamma)=I_{\gamma}^{G}$ as $\rv_{G}(g)=1$ for all $g\in G(F)$. 



\subsubsection{A variant}

In their work on the weighted fundamental lemma \cite{cl2}, Laumon and Chaudouard introduce a variant of the weighted orbital integral.

Assume that $G$ is semisimple, let $\xi\in \ka_{M}$ be a generic element. For $g\in G(F)$, they introduce the weight factor
$$
\rw_{M}^{\xi}(g)=|\{\lambda\in X_{*}(M)\mid \lambda+\xi\in \ec_{M}(g)\}|.
$$ 
It is the number of integral points in the polytope $\ec_{M}(g)-\xi$. Similar to $\rv_{M}(g)$, the weight factor $\rw_{M}^{\xi}(g)$ is invariant under the right $K$-action and the left $M(F)$-action. In particular, it descends to a function on $\xx$. Consider the following weighted orbital integral
$$
J_{M}^{\xi}(\gamma)=
\int_{T(F)\backslash G(F)} \mathbbm{1}_{\kg(\co)}\big(\Ad(g)^{-1}\gamma\big)\rw_{M}^{\xi}(g)\frac{dg}{dt}.
$$

\begin{rem}\label{general weight factor}

For general reductive algebraic group $G$, $\xi\in \ka_{M}^{G}$, as $G(F)=M(F)\cdot G_{\mathrm{der}}(F)$, we can define the weight factor $\rw_{M}^{\xi}$ uniquely by requiring it to be invariant under the left $M(F)$-action and the right $K$-action, and that as a function on $\xx$ its restriction to $\xx^{G_{\mathrm{der}}}$ coincides with the above definition for $G_{\mathrm{der}}$. In other words, for generic $\xi\in \ka_{M}^{G}$, we define
$$
\rw_{M}^{\xi}(g)=|\{\lambda\in X_{*}(M_{G_{\mathrm{der}}})\mid \lambda+\xi\in \pi_{M}^{G}(\ec_{M}(g))\}|,
$$ 
where $M_{G_{\mathrm{der}}}=M\cap G_{\mathrm{der}}$.
 Notice that the weight factor $\rv_{M}$ satisfies these conditions as well, and this justifies our definition in the general case. 
\end{rem}

The variant $J_{M}^{\xi}(\gamma)$ has a better geometric interpretation.


\begin{lem}\label{surj lattice}

Let $T(F)^{1}=T(F)\cap \ker (H_{M_{0}})$, then it is of finite volume and we have an exact sequence
$$
1\to T(F)^{1}\to T(F)\xrightarrow{H_{M_{0}}} X_{*}(M_{0})\to 1.
$$

\end{lem}

\begin{proof}

The first assertion is due to the fact that $T$ is anisotropic modulo the center of $M_{0}$. For the second assertion, only the surjectivity is non-trivial. Recall that we have the exact sequence
$$
1\to T(F)_{1}\to T(F)\xrightarrow{w_{T}} \big(X_{*}(T)_{I_{F}}\big)^{\langle \sigma\rangle}\to 1,
$$
and that the map $w_{T}$ is defined via a map $v_{T}: T(F^{\mathrm{ur}})\to \Hom(X^{*}(T)^{I_{F}},\bz)$, similar to the definition of $H_{M_{0}}$. Hence the morphism $H_{M_{0}}$ factors through $w_{T}$, and the surjectivity results from those of $w_{T}$ and the homomorphism $$\big(X_{*}(T)_{I_{F}}\big)^{\langle \sigma\rangle}=X_{*}(T)_{I_{F}}\to \Hom(X^{*}(M_{0}),\bz)=X_{*}(M_{0}).$$

\end{proof}

As a consequence, let $T(F)_{M}^{1}=T(F)\cap \,\ker (H_{M})$, let $\Lambda^{H_{M}}=\Lambda\cap \,\ker (H_{M})$, the quotient $\Lambda^{H_{M}}\backslash T(F)_{M}^{1}$ is of finite volume and we have an exact sequence
$$
1\to T(F)_{M}^{1}\to T(F)\xrightarrow{H_{M}} X_{*}(M)\to 1.
$$

\begin{prop}\label{geom weighted orbital}

We have the equality
$$
J_{M}^{\xi}(\gamma)=\mathrm{vol}_{dt}\big(\Lambda^{H_{M}}\backslash T(F)_{M}^{1}\big)^{-1}\cdot
 \Big|\Lambda^{H_{M}}\backslash \Big\{[g]\in \xx_{\gamma}(\fq)\big |\, \xi\in \ec_{M}(g)\Big\}\Big|.
$$
In particular, 
$$
J_{M_{0}}^{\xi}(\gamma)=\mathrm{vol}_{dt}\big(T(F)^{1}\big)^{-1}\cdot
 \Big|\Big\{[g]\in \xx_{\gamma}(\fq)\big |\, \xi\in \ec_{M_{0}}(g)\Big\}\Big|.
$$
\end{prop}

\begin{proof}

Let $\mathbbm{1}_{M,\,g}$ be the characteristic function of $\ec_{M}(g)$. As 
$$\ec_{M}(tg)=\ec_{M}(g)+H_{M}(t),\quad \forall\,t\in T(F),\,g\in G(F),$$ 
we have 
\begin{eqnarray*}
\sum_{t\in\, T(F)_{M}^{1}\backslash T(F)} \mathbbm{1}_{M,tg}(\xi)&=& \big|\big\{\lambda\in X_{*}(M)\mid \xi\in \ec_{M}(g)+\lambda \big\}\big|  \\
&=&\rw_{M}^{\xi}(g).
\end{eqnarray*}

Now we can rewrite
\begin{eqnarray*}
J_{M}^{\xi}(\gamma)&=&
\int_{T(F)\backslash G(F)} \mathbbm{1}_{\kg(\co)}\big(\Ad(g)^{-1}\gamma\big)\rw_{M}^{\xi}(g)\frac{dg}{dt} \\
&=&
\int_{T(F)\backslash G(F)} \mathbbm{1}_{\kg(\co)}\big(\Ad(g)^{-1}\gamma\big) \sum_{t\in\, T(F)_{M}^{1}\backslash T(F)} \mathbbm{1}_{M,tg}(\xi)
\frac{dg}{dt}     \\
&=&
\int_{T(F)_{M}^{1}\backslash  G(F)} \mathbbm{1}_{\kg(\co)}\big(\Ad(g)^{-1}\gamma\big)  \mathbbm{1}_{M,g}(\xi)\,\frac{dg}{dt}         \\
&=&
\mathrm{vol}_{dt}\big(\Lambda^{H_{M}}\backslash T(F)_{M}^{1}\big)^{-1}
\int_{\Lambda^{H_{M}}\backslash  G(F)} \mathbbm{1}_{\kg(\co)}\big(\Ad(g)^{-1}\gamma\big)  \mathbbm{1}_{M,g}(\xi)\,dg          \\
&=&
\mathrm{vol}_{dt}\big(\Lambda^{H_{M}}\backslash T(F)_{M}^{1}\big)^{-1}\cdot
 \Big|\Lambda^{H_{M}}\backslash \Big\{[g]\in \xx_{\gamma}(\fq)\big |\, \xi\in \ec_{M}(g)\Big\}\Big|.
\end{eqnarray*}

\end{proof}

In particular, $J_{M}^{\xi}(\gamma)$ is a plain count of a subset of $\xx_{\gamma}(\fq)$. In \S \ref{HN whole}, we will see that the condition $\xi\in \ec_{M}(g)$ behaves as a stability condition (We believe that it is in fact a stability condition in the sense of Mumford). In particular, there is a Harder-Narasimhan type decomposition of $\xx_{\gamma}$ associated with it.

\begin{rem}

It is time to explain why we have imposed the assumption that $T$ is totally ramified over $F$. Without it, the Frobenius $\sigma\in \gal(\overline{\mathbf{F}}_{q}/\fq)$ acts non-trivially on $X_{*}(T)_{I_{F}}$, and the morphism $T(F)\to X_{*}(M_{0})$ in lemma \ref{surj lattice} might fail to be surjective. (Indeed, it does fail for $T$ an unramified maximal torus in $\gl_{n}$.) As a consequence, the interpretation of $J_{M}^{\xi}(\gamma)$ as  in proposition \ref{geom weighted orbital} no longer holds.  
\end{rem}

For completeness, we compute the volume factors in proposition \ref{geom weighted orbital}. We have the exact sequence
$$
1\to \Lambda^{H_{M}}\to T(F)_{M}^{1}\to \Lambda\backslash T(F) \xrightarrow{H_{M}} X_{*}(M)/H_{M}(\Lambda)\to 1.
$$
Because $\Lambda$ is of finite index in $X_{*}(M_{0})$ and the morphism $X_{*}(M_{0})\to X_{*}(M)$ is surjective, the quotient $X_{*}(M)/H_{M}(\Lambda)$ is finite, and so
\begin{eqnarray}
\mathrm{vol}_{dt}\big(\Lambda^{H_{M}}\backslash T(F)_{M}^{1}\big)
&=&
\vol_{dt}\big( \Lambda\backslash T(F) \big)\cdot |X_{*}(M)/H_{M}(\Lambda)|^{-1}
\nonumber
\\
&=&
\frac{|\mathrm{coker}[X_{*}(S)_{\Gamma}\to X_{*}(T)_{\Gamma}]|}{|\ker[X_{*}(S)_{\Gamma}\to X_{*}(T)_{\Gamma}]|\cdot |X_{*}(M)/H_{M}(X_{*}(S))|}.  \label{volume factor}
\end{eqnarray}

\subsubsection{Comparison of weighted orbital integrals}

The weight factors $\rv_{M}$ and $\rw_{M}^{\xi}$ are closely related, we can compare the associated weighted orbital integrals.

\begin{thm}[Chaudouard-Laumon \cite{cl2}]\label{cl comparision}

We have the equality
$$
J_{M}(\gamma)= \vol(\ka_{M}/X_{*}(M)) \cdot J_{M}^{\xi}(\gamma).
$$

\end{thm}

\begin{rem}

For general reductive algebraic group $G$, with the definition of $\rw_{M}^{\xi}$ as explained in remark \ref{general weight factor}, the comparison theorem becomes
$$
J_{M}(\gamma)= \vol\big(\ka_{M_{G_{\mathrm{der}}}}^{G_{\mathrm{der}}}/X_{*}(M_{G_{\mathrm{der}}})\big) \cdot J_{M}^{\xi}(\gamma),
$$
as can be seen from the proof below.
\end{rem}

Laumon and Chaudouard work over the ring of ad\`eles, but their proof carries over to the local setting. We reproduce their proof here, but to simplify the exposition, we assume moreover that $G$ is simply connected. The key is to rewrite the convex polytope $\ec_{M}(g)$ as alternating differences of translations of cones. We need some notations. For $P=MN_{P}\in \cf(A)$, take a Borel subgroup $B\in \cp(A)$ contained in $P$. Let $\Delta_{B}$ be the simple roots of $\Phi(G,A)$ with respect to $B$, let $\Delta_{B,P}=\Delta_{B}\cap \Phi(N_{P},A)$ and $\Delta_{B, P}^{\vee}$ the associated coroots. The restriction $X^{*}(A)\to X^{*}(A_{M})$ induces a bijection from $\Delta_{B,P}$ to a subset of $X^{*}(A_{M})$ denoted $\Delta_{P}$. Similarly, the projection $\ka_{A}\to \ka_{M}$ induces a bijection from $\Delta_{B,P}^{\vee}$ to a subset $\Delta_{P}^{\vee}$. Obviously, the definition of $\Delta_{P}$ and $\Delta_{P}^{\vee}$ is independent of the choice of $B$. Moreover, they form basis of $\ka_{M}^{*}$ and $\ka_{M}$ respectively. Let $(\varpi_{\alpha})_{\alpha\in \Delta_{P}}$ be the basis of $\ka_{M}^{*}$ dual to $\Delta_{P}^{\vee}$.

For a generic element $\lambda\in \ka_{M}^{*}$, let 
$$
\Delta_{P}^{\lambda}=\{\alpha\in \Delta_{P}\mid \langle\lambda, \alpha\rangle<0\},
$$ 
and let $\varphi_{P}^{\lambda}$ be the characteristic function of the cone
$$
\{a\in \ka_{M}\mid \varpi_{\alpha}(a)>0,\,\forall \,\alpha\in \Delta_{P}^{\lambda};\;\varpi_{\alpha}(a)\leq 0,\,\forall \,\alpha\in \Delta_{P}\backslash \Delta_{P}^{\lambda}\}.
$$
According to Arthur \cite{a2}, the characteristic function of the convex polytope $\ec_{M}(g)$ is equal to the function
$$
a\in \ka_{M}\longmapsto \sum_{P\in \cp(M)} (-1)^{|\Delta_{P}^{\lambda}|}\varphi_{P}^{\lambda}(-H_{P}(g)+a).
$$
The proof is best illustrated by Figure 11.1 at \cite{a3}, page 63. It relies on the combinatorial identity 
$$
\sum_{F\subset S} (-1)^{|F|}=\begin{cases} 1,&\text{ if } S=\emptyset,\\
0,&\text{ otherwise},
\end{cases}
$$ 
for any finite set $S$. Now we can rewrite
\begin{eqnarray}
\rw_{M}^{\xi}(g)
&=&
\sum_{\chi\in X_{*}(M)}\sum_{P\in \cp(M)}(-1)^{|\Delta_{P}^{\lambda}|} \varphi^{\lambda}_{P}(-H_{P}(g)+\chi+\xi)\label{weight 1}
\\
\rv_{M}(g)
&=& \int_{\ka_{M}} \sum_{P\in \cp(M)} (-1)^{|\Delta_{P}^{\lambda}|}\varphi_{P}^{\lambda}(-H_{P}(g)+a)\,da. \label{weight 2}
\end{eqnarray}
We introduce an extra exponential factor to treat the infinite sum in (\ref{weight 1}), let 
$$
S_{P}(\lambda)=\sum_{\chi\in X_{*}(M)} \varphi^{\lambda}_{P}(-H_{P}(g)+\chi+\xi)e^{\langle\lambda,\chi\rangle}.
$$
The series converges absolutely for generic $\lambda$, hence
$$
\rw_{M}^{\xi}(g)=\lim_{\lambda\to 0}\sum_{P\in \cp(M)}(-1)^{|\Delta_{P}^{\lambda}|} S_{P}(\lambda),
$$
where the limit is taken for generic $\lambda\in \ka_{M}^{*}$.

We can calculate $S_{P}(\lambda)$ explicitly. Let $\xi=[\xi]_{P}+\{\xi\}_{P}$ with $[\xi]_{P}\in X_{*}(M)$ and $\{\xi\}_{P}=\sum_{\alpha\in \Delta_{P}}r_{\alpha}\alpha^{\vee}$ for some $0< r_{\alpha}<1$. After a simple change of variables, we get
\begin{eqnarray*}
S_{P}(\lambda)&=& e^{\langle \lambda,\, H_{P}(g)-[\xi]_{P}  \rangle} \sum_{\chi\in X_{*}(M)} \varphi_{P}^{\lambda}(\chi+\{\xi\}_{P})e^{\langle\lambda,\chi\rangle}\\
&=&
e^{\langle \lambda,\, H_{P}(g)-[\xi]_{P}  \rangle}
\sum_{(m_{\alpha})_{\alpha\in \Delta_{P}}}
e^{\langle\lambda,\,\sum_{\alpha}m_{\alpha}\alpha^{\vee}\rangle},
\end{eqnarray*}
where $(m_{\alpha})_{\alpha\in \Delta_{P}}$ runs over the integers satisfying $m_{\alpha}\geq 0$ for $\alpha\in \Delta_{P}^{\lambda}$ and $m_{\alpha}\leq -1$ for $\alpha\in \Delta_{P}\backslash\Delta_{P}^{\lambda}$. The geometric series can be calculated to be
$$
S_{P}(\lambda)=(-1)^{|\Delta_{P}^{\lambda}|}e^{\langle \lambda,\, H_{P}(g)-[\xi]_{P}  \rangle}\prod_{\alpha\in \Delta_{P}}\frac{1}{e^{\langle\lambda,\,\alpha^{\vee}\rangle}-1}.
$$
Let $\rc_{P}(\lambda)=\prod_{\alpha\in \Delta_{P}} (e^{\langle\lambda,\,\alpha^{\vee}\rangle}-1)$. Taking everything together, we get
\begin{equation}\label{w alt}
\rw_{M}^{\xi}(g)=
\lim_{\lambda\to 0}
\sum_{P\in \cp(M)}
\rc_{P}(\lambda)^{-1}e^{\langle \lambda,\, H_{P}(g)-[\xi]_{P}  \rangle}.\end{equation}

Similarly, we can rewrite (\ref{weight 2}) as
\begin{eqnarray*}
\rv_{M}(g)
&=& 
\lim_{\lambda\to 0}
\int_{\ka_{M}} \sum_{P\in \cp(M)} (-1)^{|\Delta_{P}^{\lambda}|}\varphi_{P}^{\lambda}(-H_{P}(g)+a) e^{\langle \lambda,\,a  \rangle}   \,da
\\
&=&
\lim_{\lambda\to 0}
\sum_{P\in \cp(M)} (-1)^{|\Delta_{P}^{\lambda}|} \int_{\ka_{M}} \varphi_{P}^{\lambda}(-H_{P}(g)+a) e^{\langle \lambda,\,a  \rangle}   \,da
\\
&=&
\lim_{\lambda\to 0}
\sum_{P\in \cp(M)}  
e^{\langle \lambda,\,H_{P}(g)  \rangle}\cdot
\vol(\ka_{M}/X_{*}(M))\prod_{\alpha\in \Delta_{P}} \langle \lambda,\,\alpha^{\vee}\rangle^{-1}
\end{eqnarray*}
Let $\rd_{P}(\lambda)=\vol(\ka_{M}/X_{*}(M))^{-1}\prod_{\alpha\in \Delta_{P}} \langle \lambda,\,\alpha^{\vee}\rangle$, we get
\begin{equation}\label{v alt}
\rv_{M}(g)=\lim_{\lambda\to 0}
\sum_{P\in \cp(M)}  \rd_{P}(\lambda)^{-1}\cdot
e^{\langle \lambda,\,H_{P}(g)  \rangle}. 
\end{equation}

To deal with limits of the form (\ref{w alt}) and (\ref{v alt}) systematically, we need Arthur's notion of \emph{$(G,M)$-family} \cite{a4}. It is a family of smooth functions $(\rr_{P}(\lambda))_{P\in \cp(M)}$ on $\ka_{M}^{*}$ which satisfy for any adjacent parabolic subgroups $(P,P')$ the property that $\rr_{P}(\lambda)=\rr_{P'}(\lambda)$ for any $\lambda$ on the hyperplane defined by the unique coroot in $\Delta_{P}^{\vee}\cap (-\Delta_{P'}^{\vee})$. For any such family, we define
$$
\rr_{M}(\lambda)=\sum_{P\in \cp(M)} \rd_{P}(\lambda)^{-1} \, \rr_{P}(\lambda),
$$
for generic $\lambda\in \ka_{M}^{*}$. Arthur has shown in \cite{a4} that the function extends smoothly over all $\ka_{M}^{*}$. Let
$$
\rr_{M}=\lim_{\lambda\to 0} \rr_{M}(\lambda).
$$
It generalizes the equation (\ref{v alt}). Indeed, the functions
$$
\rv_{P}(\lambda, g)=e^{\langle \lambda,\,H_{P}(g)  \rangle},\quad P\in \cp(M),
$$
form a $(G,M)$-family, and the resulting $\rv_{M}(g)$ is exactly Arthur's weight factor. From this point of view, we call $\rr_{M}$ the \emph{volume} of the $(G,M)$-family $(\rr_{P}(\lambda))_{P\in \cp(M)}$.

Notice that the summands in (\ref{w alt}) and (\ref{v alt}) differ by a factor
$$
\rw_{P}(\lambda,\xi)=\frac{\rd_{P}(\lambda)}{\rc_{P}(\lambda)}e^{-\langle \lambda,\,[\xi]_{P} \rangle},
$$
and that they form a $(G,M)$-family. Let $\rw_{P}(\lambda,g,\xi)=\rv_{P}(\lambda,g)\rw_{P}(\lambda,\xi),\,P\in \cp(M)$, they form a $(G,M)$-family and equation (\ref{w alt}) can be rewritten as
\begin{equation}\label{w alt 2}
\rw_{M}^{\xi}(g)=\rw_{M}(g,\xi).
\end{equation} 
In other words, we have expressed the lattice points counting weight factor $\rw_{M}^{\xi}(g)$ as the volume of the product of two $(G,M)$-families.




We need a result of Arthur on the volume of the product of two $(G,M)$-families. Let $\{\rr_{P}(\lambda)\}_{P\in \cp(M)}$, $\{\rs_{P}(\lambda)\}_{P\in \cp(M)}$ be two $(G,M)$-families. For $Q=LN_{Q}\in \cf(M)$, let
$$
\rr_{R}^{Q}(\lambda)=\rr_{RN_{Q}}(\lambda), \quad \forall \, R\in \cp^{L}(M).
$$ 
It is easy to see that $\rr_{R}^{Q}(\lambda),\,R\in \cp^{L}(M)$ form a $(L,M)$-family. The function $\rr_{M}^{Q}(\lambda)$ and the volume $\rr_{M}^{Q}$ are defined in a similar way. From the $(G,M)$-family $\{\rs_{P}(\lambda)\}_{P\in \cp(M)}$, Arthur has defined a smooth function $\rs_{Q}'(\lambda)$ on $\ka_{Q}^{*}$. The definition is quite involved and we refer the reader to \cite{a4} \S6. Let $\rs_{Q}'=\rs_{Q}'(0)$.

\begin{lem}[Arthur \cite{a4} Lem. 6.3 and Cor. 6.4]\label{vol prod}

Let $\{\rr_{P}(\lambda)\}_{P\in \cp(M)}$, $\{\rs_{P}(\lambda)\}_{P\in \cp(M)}$ be two $(G,M)$-families, let $\rr\cdot \rs$ be the product of the two $(G,M)$-families, then for any $\lambda\in \ka_{M}^{*}$, we have 
$$
(\rr\cdot \rs)_{M}(\lambda)=\sum_{Q\in \cf(M)}\rr_{M}^{Q}(\lambda)\rs_{Q}'(\lambda).
$$
In particular, 
$$
\rs_{M}(\lambda)=\sum_{P\in \cp(M)}\rs_{P}'(\lambda).
$$

\end{lem}


In our situation, this implies
\begin{equation}\label{w mul 1}
\rw_{M}^{\xi}(g)=\big(\rv(g)\cdot \rw(\xi)\big)(0)=\sum_{Q\in \cf(M)}\rv_{M}^{Q}(g)\rw_{Q}'(\xi),
\end{equation}
and
\begin{equation}\label{w mul 2}
\rw_{M}(\xi)=\sum_{P\in \cp(M)}\rw_{P}'(\xi).
\end{equation}

Similar results hold for Levi subgroups $L$ containing $M$: 
\begin{eqnarray}
\rw_{L}^{\xi_{L}}(g)&=&\sum_{R\in \cf(L)}\rv_{L}^{R}(g)\rw_{R}'(\xi_{L}), \label{w mul 3}\\
\rw_{L}(\xi_{L})&=&\sum_{Q\in \cp(L)}\rw_{Q}'(\xi_{L}), \label{w mul 4}
\end{eqnarray}
with $\rw_{R}'(\xi_{L})$ deduced from the $(G,L)$-family
$$
\rw_{Q}(\lambda,\xi_{L})=\frac{\rd_{Q}(\lambda)}{\rc_{Q}(\lambda)}e^{-\langle \lambda,\,[\xi_{L}]_{Q} \rangle}, \quad \forall\, Q\in \cp(L).
$$
Setting $g=e\in G$ in equation (\ref{w mul 3}), and notice that 
$$
\rv_{L}^{R}(e)=\begin{cases}1, & \text{ if } R\in \cp(L),\\
0,& \text{ otherwise}, 
\end{cases}
$$
we get 
\begin{equation}\label{w mul 5}
\rw_{L}^{\xi_{L}}(e)=\sum_{Q\in \cp(L)}
\rw_{Q}'(\xi_{L})=\rw_{L}(\xi_{L}),
\end{equation}
where the second equality is just the equation (\ref{w mul 4}).

\begin{lem}\label{sum of angle}

$$
\sum_{Q\in \cp(L)}
\rw_{Q}'(\xi)
=\frac{\vol(\ka_{L}/X_{*}(L))}{\vol(\ka_{M}/X_{*}(M))}\cdot \rw_{L}^{\xi_{L}}(e)=
\begin{cases}
\vol(\ka_{M}/X_{*}(M))^{-1},& \text{if } L=G,\\
0,&\text{otherwise}.
\end{cases}
$$

\end{lem}

\begin{proof}

Recall that given a $(G,M)$-family $\{\rs_{P}(\lambda)\}_{P\in \cp(M)}$, we can define a $(G,L)$-family by setting 
$$
\rs_{Q}(\lambda)=\rs_{P}(\lambda), \quad \forall\, \lambda\in \ka_{L}^{*}\subset \ka_{M}^{*},
$$
for any $P\in \cp(M),\,P\subset Q$. Moreover, the function $\rs_{Q}'(\lambda)$ deduced from the $(G,M)$-family $\{\rs_{P}(\lambda)\}_{P\in \cp(M)}$ is the same as that from the $(G,L)$-family $\{\rs_{Q}(\lambda)\}_{Q\in \cp(L)}$ by formula (6.3) in \cite{a4}.
In this way, we get the $(G,L)$-family
$\{\rw_{Q}(\lambda,\xi)\}_{Q\in \cp(L)}$ and the equality
$$
\sum_{Q\in \cp(L)}\rw_{Q}'(\xi)=\rw_{L}(\xi)=\lim_{\lambda\to 0} \sum_{Q\in \cp(L)} \rd_{Q}(\lambda)^{-1}\cdot \frac{\rd_{P}(\lambda)}{\rc_{P}(\lambda)}\cdot e^{-\langle \lambda,\,[\xi]_{P}  \rangle},
$$
by the second assertion of lemma \ref{vol prod}, where for each $Q\in \cp(L)$ we take $P\in \cp(M),\,P\subset Q$ and the limit is taken for $\lambda\in \ka_{L}^{*}$ generic. Now that
$$
\frac{\rd_{P}(\lambda)}{\rc_{P}(\lambda)}
=\frac{\vol(\ka_{L}/X_{*}(L))}{\vol(\ka_{M}/X_{*}(M))}\cdot
\frac{\rd_{Q}(\lambda)}{\rc_{Q}(\lambda)}
$$
and $\langle \lambda,\,[\xi]_{P}  \rangle=\langle \lambda,\,[\xi_{L}]_{Q}  \rangle$ for any $\lambda\in \ka_{L}^{*}$, we get
$$
\sum_{Q\in \cp(L)}\rw_{Q}'(\xi)=\frac{\vol(\ka_{L}/X_{*}(L))}{\vol(\ka_{M}/X_{*}(M))}\cdot \rw_{L}(\xi_{L})
=\frac{\vol(\ka_{L}/X_{*}(L))}{\vol(\ka_{M}/X_{*}(M))}\cdot \rw_{L}^{\xi_{L}}(e),
$$
where the last equality follows from equation (\ref{w mul 5}).

\end{proof}

By (\ref{w mul 1}), we can rewrite $J_{M}^{\xi}(\gamma)$ as
\begin{eqnarray*}
J_{M}^{\xi}(\gamma)
&=&
\int_{T(F)\backslash G(F)} \mathbbm{1}_{\kg(\co)}\big(\Ad(g)^{-1}\gamma\big)\rw_{M}^{\xi}(g)\frac{dg}{dt}\\
&=&
\int_{T(F)\backslash G(F)} \mathbbm{1}_{\kg(\co)}\big(\Ad(g)^{-1}\gamma\big)\Bigg[\sum_{Q\in \cf(M)}\rv_{M}^{Q}(g)\rw_{Q}'(\xi)
\Bigg]\frac{dg}{dt}.
\end{eqnarray*}
As $\rv_{M}^{Q}$ is also left $T(F)$-invariant, we can define
$$
J_{M}^{Q}(\gamma)
=
\int_{T(F)\backslash G(F)} \mathbbm{1}_{\kg(\co)}\big(\Ad(g)^{-1}\gamma\big)\rv_{M}^{Q}(g)\frac{dg}{dt}.
$$
Let $Q=LN_{Q}$ be the standard Levi decomposition. Let $dl$ be the Haar measure on $L(F)$ normalized by $\vol_{dl}(L(\co))=1$, let $dn $ be the Haar measure on $N_{Q}(F)$ normalized by $\vol_{dn}(N_{Q}(\co))=1$. Using Iwasawa decomposition, we can rewrite $J_{M}^{Q}(\gamma)$ as
\begin{eqnarray}
J_{M}^{Q}(\gamma)
&=&
\int_{T(F)\backslash L(F)} \int_{N_{Q}(F)}\int_{K} 
\mathbbm{1}_{\kg(\co)}\big(\Ad(nlk)^{-1}\gamma\big)\rv_{M}^{Q}(nlk)\, dk\cdot dn \cdot \frac{dl}{dt}    \nonumber
\\
&=&
\int_{T(F)\backslash L(F)} \int_{N_{Q}(F)}\int_{K}
\mathbbm{1}_{\kg(\co)}\big(\Ad(nl)^{-1}\gamma\big)\rv_{M}^{Q}(l)\, dk\cdot dn \cdot \frac{dl}{dt}     \nonumber
\\
&=&
\int_{T(F)\backslash L(F)} 
\Big[\int_{N_{Q}(F)}\int_{K}
\mathbbm{1}_{\kg(\co)}\big(\Ad(nl)^{-1}\gamma\big) \,dk\cdot dn\Big]  \rv_{M}^{L}(l)\,\frac{dl}{dt},   \label{JQ}
\end{eqnarray}
where in the second and third lines we have used the equalities $\rv_{M}^{Q}(nlk)=\rv_{M}^{Q}(l)$ and $\rv_{M}^{Q}(l)=\rv_{M}^{L}(l)$ respectively, they follow directly from definitions. Notice that 
\begin{eqnarray*}
\int_{N_{Q}(F)}\int_{K}
\mathbbm{1}_{\kg(\co)}\big(\Ad(nl)^{-1}\gamma\big) \,dk\cdot dn
&=&\big|\big\{[nl]\in N_{Q}(F)lK/K\mid \Ad(nl)^{-1}\gamma\in \kg(\co)\big\}\big|
\\
&=&\big|\big(f_{Q}^{-1}([l])\cap \xx_{\gamma}\big)(\fq)\big|
\\
&=& q^{\val(\det(\ad\gamma\mid\kn_{Q,F}))}\cdot \mathbbm{1}_{\mathfrak{l}(\co)}\big(\Ad(l)^{-1}\gamma\big),
\end{eqnarray*}
where the last equality follows from Prop. \ref{KL retraction}. Continuing the calculation (\ref{JQ}), we get
\begin{eqnarray*}
J_{M}^{Q}(\gamma)
&=&
q^{\val(\det(\ad\gamma\mid\kn_{Q,F}))}\cdot\int_{T(F)\backslash L(F)} 
\mathbbm{1}_{\mathfrak{l}(\co)}\big(\Ad(l)^{-1}\gamma\big) 
\rv_{M}^{L}(l)\,\frac{dl}{dt}
\\
&=&
q^{\val(\det(\ad\gamma\mid\kn_{Q,F}))}\cdot J_{M}^{L}(\gamma).
\end{eqnarray*}
Combining all the above calculations, we get
\begin{eqnarray*}
J_{M}^{\xi}(\gamma)
&=&
\sum_{L\in \cl(M)}\sum_{Q\in \cp(L)} J_{M}^{Q}(\gamma) \cdot \rw_{Q}'(\xi)=\sum_{L\in \cl(M)}\sum_{Q\in \cp(L)}q^{\val(\det(\ad\gamma\mid\kn_{Q,F}))} J_{M}^{L}(\gamma) \cdot \rw_{Q}'(\xi)\\
&=&
\sum_{L\in \cl(M)}q^{\frac{1}{2}\val(\det(\ad\gamma\mid\kg_{F}/\mathfrak{l}_{F}))} J_{M}^{L}(\gamma) \sum_{Q\in \cp(L)}\rw_{Q}'(\xi)\\
&=& {\vol(\ka_{M}/X_{*}(M))}^{-1} \cdot J_{M}(\gamma),
\end{eqnarray*}
where the last equality follows from lemma \ref{sum of angle}. This finishes the proof of theorem \ref{cl comparision}.

\section{Counting points by Arthur-Kottwitz reduction}\label{Section AK count}

From now on, we will assume that $G_{\mathrm{der}}$ is simply connected. 
The general case can be reduced to this one by focusing on each connected  component. This extra assumption gives some technical convenience, for example, $M_{0,\mathrm{der}}$ will be simply connected, $\Lambda_{M_{0}}$ will be torsion free and we get an inclusion $\Lambda_{M_{0}}\hookrightarrow \ka_{M_{0}}$. Moreover, we have $X_{*}(M_{0})=\Lambda_{M_{0}}$, according to \cite{cl2}, lemma 11.6.1.

Fix $M\in \cl(M_{0})$, let $\Pi$ be a sufficiently regular positive $(G,M)$-orthogonal family. We count the number of points on $\Lambda^{H_{M}}\backslash\xx_{\gamma}^{\nu_{0}}(\Pi)$, $\nu_{0}\in \Lambda_{G}$. Generalizing our work \cite{chen2}, we show that it can be reduced to counting points on the \emph{intermediate fundamental domains} $F_{\gamma}^{L,M},\,L\in \cl(M)$, and the counting result depends \emph{quasi-polynomially} on the truncation parameter. Moreover, counting points on $\Lambda^{H_{M}}\backslash F_{\gamma}^{L,M}$ can be further reduced to that of the fundamental domains $F_{\gamma}^{M'}$ for some $M'\in \cl(M_{0})$ ``transversal'' to $M$.

\subsection{Truncations on the affine grassmannian}


Recall the following definition of Arthur \cite{a}, which is a formalization of the orthogonal  properties in proposition \ref{arthur orthogonal}.

\begin{defn}
A family $\Pi=(\lambda_{P})_{P\in \cp(M)}$ of elements in $\ka_{M}^{G}$ is called a \emph{positive} $(G,M)$-\emph{orthogonal family} if it satisfies
$$
\lambda_{P_{1}}-\lambda_{P_{2}}=n_{P_{1},P_{2}}\cdot \pi_{M}^{G}(\beta_{P_{1},P_{2}}), \quad \text{ with } n_{P_{1},P_{2}}\in \br_{\geq 0},
$$
for any two adjacent parabolic subgroups $P_{1},\, P_{2}\in \cp(M)$. 
\end{defn}

Given such a positive $(G,M)$-orthogonal family, we will denote again by $\Pi$ the convex hull of the $\lambda_{P}$'s. For $Q=LN_{Q}\in \cf(M)$, parallel to $\ec_{M}^{Q}(x)$, we denote by $\Pi^{Q}$ the face of $\Pi$ whose vertices are $\lambda_{P},\,P\in \cp(M),\,P\subset Q$. With the projection $\pi_{M}^{L}$, it can be seen as a positive $(L,M)$-orthogonal family. This sets up a bijection between the set $\cf(M)$ and the set of the faces of $\Pi$. 
Moreover, we denote by $\lambda_{Q}$ or $\lambda_{Q}(\Pi)$ the element $\pi_{M,L}(\lambda_{P'})$ for any $P'\in \cp(M),\,P'\subset Q$. One can show that $\big(\lambda_{Q}(\Pi)\big)_{Q\in \cp(L)}$ forms a positive $(G,L)$-orthogonal family.
Later on, we also use the notation $(\lambda_{\bar{w}})_{\bar{w}\in W/W_{M}}$ for $(\lambda_{\bar{w}\cdot P})_{\bar{w}\in W/W_{M}}$, and we use the notation $\lambda_{w}(\Pi)$ or $\lambda_{w\cdot P}(\Pi)$ to indicate the vertex of $\Pi$ indexed by $w\cdot P$.

Following Chaudouard and Laumon \cite{cl1}, we define the \emph{truncated affine grassmannian} $\xx(\Pi)$ to be
$$
\xx(\Pi)=\big\{x\in \xx\,\big|\, \pi_{M}^{G}\big(\ec_{M}(x)\big)\subset \Pi\big\}.
$$
We want to point out that its connected components are also parametrized by $\Lambda_{G}$, but they are not isomorphic in general. However, there is periodicity in the connected components: Let $G^{\ad}$ be the adjoint group of $G$, let $c_{G}:\Lambda_{G}\to \Lambda_{G^{\ad}}$ be the projection induced by the natural projection $T\to T/Z_{G}$. For $\nu,\nu'\in \Lambda_{G}$, we have
$$
\xx^{\nu}(\Pi)=\xx^{\nu'}(\Pi), \quad \text{if } c_{G}(\nu)=c_{G}(\nu'),
$$ 
because they can be translated to each other by elements in $Z_{G}(F)$.

For regular element $\gamma\in \kt(\co)$, we can truncate the affine Springer fiber $\xx_{\gamma}$ similarly by defining 
$$
\xx_{\gamma}(\Pi)=\xx_{\gamma}\cap \xx(\Pi),
$$
and the same observation on the connected components of $\xx(\Pi)$ holds also for $\xx_{\gamma}(\Pi)$.

\subsection{The intermediate fundamental domain}\label{intermediate}



We generalize our construction of the fundamental domain $F_{\gamma}$ in \cite{chen2}.\footnote{In \cite{chen2}, we have confused $\Lambda,\,\Lambda_{M_{0}}$ and $\pi_{0}(T(F))$. 
With our current notations, there are morphisms $\Lambda\to \Lambda_{M_{0}}$ and $\Lambda\to \pi_{0}(T(F))$. Generally, they are not isomorphic. In particular, $F_{\gamma}$ is \emph{not} a fundamental domain for the $\Lambda$-action, i.e. $\xx_{\gamma}\neq \bigcup_{\lambda\in \Lambda}\lambda\cdot F_{\gamma}$. 
Moreover, the group $\pi_{0}(T(F))$ may have complicated torsion subgroup, this implies that $F_{\gamma}$ may have complicated irreducible components as well, contrary to our expectation there. Actually, there should be a bijection between $\pi_{0}(F_{\gamma})$ and $\pi_{0}(F_{\gamma}^{M_{0}})$, and both are isomorphic to $\pi_{0}(T(F))_{\mathrm{tor}}$.
Nevertheless, other results of \cite{chen2} hold if we assume that $G_{\mathrm{der}}$ is simply connected, and the general case can be reduced to that one. This extra assumption is to make sure that for any Levi subgroup $M\in \cl(M_{0})$ we have $\Lambda_{M}$ being torsion free and we get an inclusion $\Lambda_{M}\hookrightarrow \ka_{M}$, they hold as $M_{\mathrm{der}}$ is simply connected. Moreover, we have $X_{*}(M)=\Lambda_{M}$, according to \cite{cl2} lemma 11.6.1.
}

Let $P_{1}=MN_{1},P_{2}=MN_{2}\in \cp(M)$ be two adjacent parabolic subgroups. 
Let $m_{\alpha}$ be the unique positive integer such that the image of $\alpha^{\vee}$ in $\Lambda_{M}$ is equal to $m_{\alpha}\cdot \beta_{P_{1},P_{2}}$. Let
$$
n(\gamma, P_{1},P_{2})=\sum_{\alpha\in \Phi(N_{1}, T_{\overline{F}})\cap \Phi(N_{2}^{-}, T_{\overline{F}})} \val(\alpha(\gamma))\cdot m_{\alpha}.
$$
It can be verified that $n(\gamma, P_{1},P_{2})$ is an integer.

\begin{prop}[Goresky-Kottwitz-MacPherson]\label{gkm bound}

Let $x\in \xx_{\gamma}$. 

\begin{enumerate}

\item For any two adjacent parabolic subgroups $P_{1}, P_{2}\in \cp(M)$, we have
$$
n(x,P_{1},P_{2})\leq n(\gamma,P_{1},P_{2}).
$$

\item The point $x$ is regular in $\xx_{\gamma}$ if and only if the following two conditions hold:
\begin{enumerate}
\item
the point $f_{P}(x)$ is regular in $\xx^{M}_{\gamma}$ for all $P\in \cp(M)$;

\item
for any two adjacent parabolic subgroups $P_{1}, P_{2}$ in $\cp(M)$, one has
$$
n(x,P_{1},P_{2})=n(\gamma,P_{1},P_{2}).
$$

\end{enumerate}
\end{enumerate}
\end{prop}

Notice that although Goresky, Kottwitz and MacPherson work over the field $F=\bc(\!(\ep)\!)$, their proof works for any field $F=k(\!(\ep)\!)$ with $\mathrm{char}(k)>|W|$. 
Their result motivates our definition:

\begin{defn}
Take a regular point $x_{0}\in \xx_{\gamma}^{\reg}$. Let 
$$
F_{\gamma}^{G,M}=\{x\in \xx_{\gamma}\mid \ec_{M}(x)\subset \ec_{M}(x_{0}), \,\nu_{G}(x)=\nu_{G}(x_{0})\}.
$$
We call it an \emph{intermediate fundamental domain} of $\xx_{\gamma}$ with respect to $M$.
\end{defn}

We should have used the notation $F_{\gamma,x_{0}}^{G,M}$ to indicate the dependence on $x_{0}$, but they are isomorphic to each other for any choice of the regular point $x_{0}$. Indeed, for any two regular points $x_{1},x_{2}$, we can find $t\in T(F)$ such that $x_{1}=t\cdot x_{2}$. Now that $\ec_{M}(tx)=\ec_{M}(x)+H_{M}(t), \,\forall\,x\in \xx$, the intermediate fundamental domain given by $x_{1}$ is just the translation by $t$ of that given by $x_{2}$. Notice that for $M=M_{0}$, we recover the fundamental domain $F_{\gamma}$. For simplicity, we assume that $\nu_{G}(x_{0})=0$.


Unlike the fundamental domain, the intermediate $F_{\gamma}^{G,M}$ is no longer of finite type for $M\supsetneq M_{0}$. Nonetheless, we have

\begin{prop}\label{quot finite 1}

The free discrete abelian group $\Lambda^{H_{M}}$ acts freely on $F_{\gamma}^{G,M}$, and the quotient $\Lambda^{H_{M}}\backslash F_{\gamma}^{G,M}$ is of finite type.

\end{prop}

\begin{proof}

Recall that $\Lambda^{H_{M}}=\Lambda\cap \ker(H_{M})$ by definition, hence it preserves  $F_{\gamma}^{G,M}$ because left translation by $m\in \ker(H_{M})$ doesn't change the polytope $\ec_{M}(x)$ due to the property
$$
H_{P}(mx)=H_{M}(m)+H_{P}(x),\quad \forall\, m\in M(F),\,x\in \xx,\,P\in \cp(M).
$$

For the finiteness issue. Let $\Lambda_{M_{0}}^{H_{M}}\subset \Lambda_{M_{0}}$ be the kernel of the natural projection $\Lambda_{M_{0}}\to \Lambda_{M}$. By definition, we have
$$
\pi_{M_{0},M}^{-1}(\ec_{M}(x_{0}))=\bigcup_{\nu\in \Lambda_{M_{0}}^{H_{M}}}(\nu+\ec_{M_{0}}(x_{0})),
$$
which implies that
$$
F_{\gamma}^{G,M}=\bigcup_{\nu\in \Lambda_{M_{0}}^{H_{M}}}\xx_{\gamma}^{0}\big(\nu+\ec_{M_{0}}(x_{0})\big).
$$
Now that $\Lambda\cong X_{*}(S)$ and $X_{*}(S)\hookrightarrow \Lambda_{M_{0}}$ is of finite index, the quotient $\Lambda^{H_{M}}\backslash \Lambda_{M_{0}}^{H_{M}}$ is of finite cardinal. Hence the quotient $\Lambda^{H_{M}}\backslash F_{\gamma}^{G,M}$ is dominated by union of finitely many translations of $F_{\gamma}$ under the natural projection $F_{\gamma}^{G,M}\to \Lambda^{H_{M}}\backslash F_{\gamma}^{G,M}$. As $F_{\gamma}$ is of finite type, so is the quotient $\Lambda^{H_{M}}\backslash F_{\gamma}^{G,M}$.


\end{proof}

Similar proof applies to:

\begin{prop}\label{quot finite 2}

Let $\Pi$ be a regular positive $(G,M)$-orthogonal family. For any $\nu\in \Lambda_{G}$, the free discret abelian group $\Lambda^{H_{M}}$ acts freely on $\xx_{\gamma}^{\nu}(\Pi)$, and the quotient $\Lambda^{H_{M}}\backslash \xx_{\gamma}^{\nu}(\Pi)$ is of finite type. In particular, 
$$
|\big(\Lambda^{H_{M}}\backslash \xx_{\gamma}^{\nu}(\Pi)\big)(\fq)|<\infty.
$$
\end{prop}

\begin{rem}

In the definition of (weighted) orbital integral, we are concerned more about analogues of  $\Lambda^{H_{M}}\backslash \big(\xx_{\gamma}^{\nu}(\Pi)(\fq)\big)$, but notice that there is bijection between 
$$
\big(\Lambda^{H_{M}}\backslash \xx_{\gamma}^{\nu}(\Pi)\big)(\fq)\quad \text{and} \quad\Lambda^{H_{M}}\backslash \big(\xx_{\gamma}^{\nu}(\Pi)(\fq)\big),
$$
because $\Lambda^{H_{M}}$ acts freely on $\xx_{\gamma}^{\nu}(\Pi)$ and the Galois group $\gal(\overline{\mathbf{F}}_{q}/\fq)$ acts trivially on $\Lambda^{H_{M}}$. We will decompose the scheme $\Lambda^{H_{M}}\backslash \xx_{\gamma}^{\nu}(\Pi)$ in different ways, the bijection above implies that we can deduce equality of rational points over $\fq$ from the decomposition of schemes. 
\end{rem}

In the following, we simplify the notation $\ec_{M}(x_{0})$ to $\Sigma_{\gamma}^{G,M}$. For $\nu\in \Lambda_{G}$, let 
$$
F_{\gamma}^{G,M, \nu}:=\xx_{\gamma}^{\nu}(\Sigma_{\gamma}^{G,M}).
$$
As we have explained before, it depends only on the class $c_{G}(\nu)\in \Lambda_{G^{\ad}}$. For $M=M_{0}$, we simplify $\Sigma_{\gamma}^{G,M_{0}}$ to $\Sigma_{\gamma}$ and $F_{\gamma}^{G,M_{0}, \nu}$ to $F_{\gamma}^{\nu}$.

\subsection{The Arthur-Kottwitz reduction}\label{AK red}

Recall that we can reduce the geometry of $\xx_{\gamma}$ to that of its fundamental domain by the Arthur-Kottwitz reduction \cite{chen2}. The construction can be generalized to our current setting.


Let $Q_{0}$ be the unique parabolic subgroup in $\cp(M)$ which contains $P_{0}$. Let $\varsigma\in \ka_{M}^{G}$ be such that $\alpha(\varsigma)$ is positive but almost equal to $0$  for any $\alpha\in \Delta_{Q_{0}}$. Let $D_{M}=(\lambda_{P})_{P\in \cp(M)}$ be the $(G,M)$-orthogonal family given by
\begin{equation}\label{slight expansion}
\lambda_{P}=H_{P}(x_{0})+w\cdot \varsigma,
\end{equation}
where $w\in W$ is any element satisfying $P=w\cdot Q_{0}$. For $Q=LN_{Q}\in \cf(M)$, define $R^{G,M}_{Q}$ to be the subset of $ \ka_{M}^{G}$ satisfying conditions
\begin{eqnarray*}
\pi_{M}^{L}(a)&\subset& D_{M}^{Q};\\
\alpha(\pi_{M,L}(a))&\geq& \alpha(\pi_{M,L}(\lambda_{Q})),\quad\forall\,\alpha\in \Delta_{Q}.
\end{eqnarray*}
This gives us a partition which dates back at least to Arthur \cite{a2}:
\begin{equation}\label{partAK}
\ka_{M}^{G}=\bigcup_{Q\in \cf(M)} R^{G,M}_{Q}.
\end{equation} 
It induces a disjoint partition of $\Lambda_{M}$ via the map $\Lambda_{M}\to \ka_{M}^{G}$, as we have perturbed $(H_{P}(x_{0}))_{P\in \cp(M)}$ with $\varsigma$. The Fig. \ref{partitiongl3} gives an illustration of the partition for the group $\gl_{3}$ and $M_{0}=T=A$.

\begin{figure}[h]
\begin{center}
\begin{tikzpicture}[scale=0.6]

\draw (-2.5,3.17)--(0.5,3.17);
\draw (-4,0.57)--(-1.5,-3.75);
\draw (1.5,-3.75)--(3,-1.15);
\draw (-4,0.57)--(-2.5,3.17);
\draw (-1.5,-3.75)--(1.5,-3.75);
\draw (0.5,3.17)--(3,-1.15);

\draw [red] (-2.5,3.17)--(-2.5, 5.17);
\draw [red] (0.5,3.17)--(0.5, 5.17);
\draw [red] (-1.5,-3.75)--(-1.5, -5.75);
\draw [red] (1.5,-3.75)--(1.5, -5.75);
\draw [red] (-4.23,4.17)--(-2.5,3.17);
\draw [red] (2.23,4.17)--(0.5,3.17);
\draw [red] (-4,0.57)--(-5.75, 1.57);
\draw [red] (-4,0.57)--(-5.75, -0.43);
\draw [red] (-1.5,-3.75)--(-3.23,-4.75);
\draw [red] (1.5,-3.75)--(3.23,-4.75);
\draw [red] (4.73,0.15)--(3,-1.15);
\draw [red] (4.73,-2.15)--(3,-1.15);

\node [blue] at (-2.5, -5.5) {$R_{B^{-}}$};
\node [blue] at (-4.5, -2.6) {$R_{P^{-}}$};
\node [blue] at (3,2) {$R_{P}$};
\node [blue] at (1.4,4.6) {$R_{B}$};
\node [blue] at (0,0) {$D_{0}$};



\end{tikzpicture}
\caption{Partition of $\ka_{A}^{G}$ for $\gl_{3}$.}
\label{partitiongl3}
\end{center}
\end{figure}

Similar to the key lemma 3.1 of \cite{chen2}, we have the following result due to proposition \ref{gkm bound}:

\begin{lem}\label{AK basic}
For any $x\in \xx_{\gamma}$, there exists a unique $Q\in \cf(M)$ such that 
$$
\pi_{M}^{G}\big(\ec_{M}^{Q}(x)\big)\subset R^{G,M}_{Q}.
$$ 
\end{lem}

The referee has suggested an equivalent form of the lemma, which is much easier to understand and to prove: Let $\ka_{M}^{G}=\bigcup_{Q\in \cf(M)} R'_{Q}$ be the partition as above attached to the \emph{positive} $(G,M)$-orthogonal family $\big(H_{P}(x_{0})-H_{P}(x)+w\cdot \varsigma\big)_{P\in \cp(M)}$, then the statement is equivalent to the existence of a unique $Q\in \cf(M)$ such that $0\in R'_{Q}$. Here the positiveness of the $(G,M)$-orthogonal family is due to proposition \ref{gkm bound}.
Let 
$$
S_{Q}^{G,M}:=\{x\in \xx_{\gamma}\mid \pi_{M}^{G}\big(\ec_{M}^{Q}(x)\big)\subset R_{Q}^{G,M}\}.
$$
We get thus a disjoint partition 
\begin{equation}\label{AK0}
\xx_{\gamma}=\xx_{\gamma}(D_{M})\cup\bigcup_{\substack{Q\in \cf(M)\\Q\neq G}}S^{G,M}_{Q}.
\end{equation} 
For each parabolic subgroup $Q=LN_{Q}\in \cf(M)$, consider the restriction of the retraction $f_{Q}: \xx\to \xx^{L}$ to $S^{G,M}_{Q}$, its image is $S^{G,M}_{Q}\cap \xx_{\gamma}^{L}$. Recall that the connected components of $\xx_{\gamma}^{L}$ are fibers of the map $\nu_{L}:\xx_{\gamma}^{L}\to \Lambda_{L}$. For $\nu\in \Lambda_{L}$, let $\xx_{\gamma}^{L,\nu}$ be its fiber at $\nu$. Let 
$$
S_{Q}^{G,M,\nu}=S_{Q}^{G,M}\cap f_{Q}^{-1}(\xx_{\gamma}^{L,\nu}),
$$ 
we have
$$
S_{Q}^{G,M,\nu}\cap \xx_{\gamma}^{L,\nu}=\xx_{\gamma}^{L,\nu}(D_{M}^{Q}).
$$

\begin{prop}\label{klretraction}

The strata $S_{Q}^{G,M,\nu}$ are locally closed sub-schemes of $\xx_{\gamma}$, and the retraction $f_{Q}: S_{Q}^{G,M,\nu}\to \xx_{\gamma}^{L,\nu}(D_{M}^{Q})$ is an iterated affine fibration over $\fq$ of dimension
$$
\val(\det(\ad(\gamma\mid \kn_{Q,F}))).
$$
\end{prop}

Indeed, by the bound on $\ec_{M}(x)$ given by proposition \ref{gkm bound}, we get
$$
S_{Q}^{G,M,\nu}=\xx_{\gamma}\cap f_{Q}^{-1}(\xx_{\gamma}^{L,\nu}(D_{M}^{Q})).
$$
It is an iterated affine fibration over $\xx_{\gamma}^{L,\nu}(D_{M}^{Q})$ by proposition \ref{KL retraction}.

The decomposition (\ref{AK0}) can thus be refined to 
\begin{equation}\label{AK}
\xx_{\gamma}=\xx_{\gamma}(D_{M})\cup \bigcup_{\substack{Q=LN_{Q}\in \cf(M)\\ Q\neq G}}\bigcup_{\nu\in \Lambda_{L}\cap \pi_{L}(R^{G,M}_{Q})} S_{Q}^{G,M,\nu}.
\end{equation}
where we have loosely used $\Lambda_{L}\cap \pi_{L}(R_{Q}^{G,M})$ to mean elements in $\Lambda_{L}$ whose projection to $\ka_{L}^{G}$ lies in $\pi_{L}(R_{Q}^{G,M})$. Similar notations will be used later on. 
The decomposition (\ref{AK}) will also be called \emph{the Arthur-Kottwitz reduction}. Notice that the stratum $S_{Q}^{G,M,\nu}$ is an iterated affine fibration over $\xx_{\gamma}^{L,\nu}(D_{M}^{Q})=F_{\gamma}^{L,M,\nu}$, and the later is related to $F_{\gamma}^{L,M}$ again by Arthur-Kottwitz reduction, similar to that explained in lemma 3.4 of \cite{chen2}.

As in \cite{chen2}, the existence of Arthur-Kottwitz reduction implies:

\begin{cor}

For any $\gamma\in \kt(\co)$, suppose that $F_{\gamma}^{L,M}$ is cohomologically pure for any proper Levi subgroup $L\in \cl(M)$. Then $\xx_{\gamma}$ is cohomologically pure if and only if $F_{\gamma}^{G,M}$ is.
\end{cor}

We can restrict the Arthur-Kottwitz reduction to the truncated affine Springer fibers.
A positive $(G,M)$-orthogonal family $\Pi=(\mu_{P})_{P\in \cp(M)}$ is said to be \emph{regular} with respect to $D_{M}$ if $\mu_{P}\in R^{G,M}_{P},\,\forall\,P\in \cp(M)$. In this case, each $S_{Q}^{G,M,\nu}$ is either contained in $\xx_{\gamma}^{\nu}(\Pi)$ or disjoint from it. So we have
\begin{equation}\label{AK truncated}
\xx_{\gamma}(\Pi)=\xx_{\gamma}(D_{M})\cup \bigcup_{\substack{Q=LN_{Q}\in \cf(M)\\ Q\neq G}}\bigcup_{\substack{\nu\in \Lambda_{L}\cap \pi_{L}(R^{G,M}_{Q})
\\ \cap \pi_{L}(\Pi)}} S_{Q}^{G,M,\nu}.
\end{equation}

The reduction can be further restricted to each connected component of $\xx_{\gamma}(\Pi)$:  
\begin{equation}\label{AK connected}
\xx_{\gamma}^{\nu_{0}}(\Pi)=\xx_{\gamma}^{\nu_{0}}(D_{M})\cup \bigcup_{\substack{Q=LN_{Q}\in \cf(M)\\ Q\neq G}}\bigcup_{\substack{\nu\in \Lambda^{\nu_{0}}_{L}\cap \pi_{L}(R^{G,M}_{Q})
\\ \cap \pi_{L}(\Pi)}} S_{Q}^{G,M, \nu}.
\end{equation}

As we have explained, left translation by elements in $\ker(H_{M})$ doesn't change the polytope $\ec_{M}(x)$, hence the group $\Lambda^{H_{M}}$ acts on each item of equation (\ref{AK connected}). Now that we have finiteness results--proposition \ref{quot finite 1} and \ref{quot finite 2}, 
combined with proposition \ref{klretraction} and the periodicity of $F_{\gamma}^{L,M,\nu}$ in $\nu\in \Lambda_{L}$, the equation (\ref{AK connected}) implies counting points equality:

\begin{cor}\label{AK count}

We have the equality
\begin{eqnarray*}
|\big(\Lambda^{H_{M}}\backslash\xx_{\gamma}^{\nu_{0}}(\Pi)\big)(\fq)|
&=&
|\big(\Lambda^{H_{M}}\backslash F_{\gamma}^{G,M,\nu_{0}}\big)(\fq)|+\sum_{\substack{Q=LN_{Q}\in \cf(M)\\Q\neq G}}\sum_{\mu\in \Lambda_{L^{\ad}}} 
q^{\val(\det(\ad\gamma\mid\kn_{Q,F}))}\cdot \\
&&|\big(\Lambda^{H_{M}}\backslash F_{\gamma}^{L,M,\mu}\big)(\fq)|\cdot 
|\Lambda^{\nu_{0}}_{L}\cap \pi_{L}(R_{Q}^{G,M})\cap \pi_{L}(\Pi)\cap c_{L}^{-1}(\mu)|.\nonumber
\end{eqnarray*}

\end{cor}

Notice that the term $|\Lambda^{\nu_{0}}_{L}\cap \pi_{L}(R_{Q}^{G,M})\cap \pi_{L}(\Pi)\cap c_{L}^{-1}(\mu)|$ counts the number of lattice points in a polytope. Well-known  techniques from toric geometry tells us that the counting result depends \emph{quasi-polynomially} on the size of the polytope.

\begin{rem}\label{general AK}

As the above constructions rely ultimately on the bound of $\ec_{M}(x)$ given by proposition \ref{gkm bound}, they continue to work if we replace $\Sigma_{\gamma}^{G,M}$ from the beginning by any integral positive $(G,M)$-orthogonal family $\Sigma$ which satisfies
\begin{equation}\label{gkm inequality}
\lambda_{P_{1}}(\Sigma)-\lambda_{P_{2}}(\Sigma)=n_{P_{1},P_{2}}\cdot\beta_{P_{1},P_{2}},\quad \text{ with } n_{P_{1},P_{2}}\geq n(\gamma,P_{1},P_{2}),
\end{equation} 
for any two adjacent parabolic subgroups $P_{1},P_{2}\in \cp(M)$. 
The resulting decomposition will also be called the Arthur-Kottwitz reduction.

\end{rem}

\subsection{Counting points on the intermediate fundamental domains}

Although the intermediate fundamental domains $F_{\gamma}^{G,M}$ looks like something new, it turns out that counting points of $\Lambda^{H_{M}}\backslash F_{\gamma}^{G,M}$ can be reduced to that of the fundamental domains.

As explained in the proof of proposition \ref{quot finite 1}, we have 
$$
\pi_{M_{0},M}^{-1}(\Sigma_{\gamma}^{G,M})=\bigcup_{\nu\in \Lambda_{M_{0}}^{H_{M}}}(\nu+\Sigma_{\gamma}),
$$
and
$$
F_{\gamma}^{G,M}=\bigcup_{\nu\in \Lambda_{M_{0}}^{H_{M}}}\xx_{\gamma}^{0}\big(\nu+\Sigma_{\gamma}\big).
$$
Let $P_{0}^{M}=P_{0}\cap M$, let $\beta_{1}^{\vee},\cdots,\beta_{r'}^{\vee}\in \Lambda_{M_{0}}^{H_{M}}$ be a basis of $\ka_{M_{0}}^{M}$ which is positive with respect to $P_{0}^{M}$. For $\mu_{1},\mu_{2}\in \Lambda_{M_{0}}^{H_{M}}$, we say that $\mu_{1}\preccurlyeq_{P_{0}^{M}} \mu_{2}$ if $\mu_{2}-\mu_{1}$ is a linear combination of $\beta_{i}^{\vee}$'s with positive coefficients. This defines a partial order $\preccurlyeq_{P_{0}^{M}}$ on $\Lambda_{M_{0}}^{H_{M}}$.
For $\mu\in \Lambda_{M_{0}}^{H_{M}}$, let 
$$
\Lambda_{M_{0},\preccurlyeq \mu}^{H_{M}}:=\{\nu\in \Lambda_{M_{0}}^{H_{M}}\mid \nu\preccurlyeq_{P_{0}^{M}} \mu\} 
\quad \text{and} \quad
\Pi_{\gamma,\preccurlyeq \mu}^{G,M}:=\bigcup_{\nu\in \Lambda_{M_{0},\preccurlyeq \mu}^{H_{M}}}\big(\nu+\Sigma_{\gamma}\big).
$$
Then $\Pi_{\gamma,\preccurlyeq \mu}^{G,M}$ is a semi-infinite polytope in $\ka_{M_{0}}^{G}$ and $\Lambda_{M_{0},\preccurlyeq \mu}^{H_{M}}$ is the integral points in it.
Similar definitions for $\Lambda_{M_{0},\prec \mu}^{H_{M}}$
and $\Pi_{\gamma,\prec \mu}^{G,M}$. But notice that $\Pi_{\gamma,\prec \mu}^{G,M}$ is not a semi-infinite polytope, it is the union of finitely many semi-infinite polytopes of the form $\Pi_{\gamma,\preccurlyeq \mu'}^{G,M}$, $\mu'\in \Lambda_{M_{0}}^{H_{M}}$.
Let
$$
F_{\gamma,\preccurlyeq \mu}^{G,M}:=\xx_{\gamma}^{0}(\Pi_{\gamma,\preccurlyeq \mu}^{G,M})
=\bigcup_{\nu\in \Lambda_{M_{0},\preccurlyeq \mu}^{H_{M}}}\xx_{\gamma}^{0}\big(\nu+\Sigma_{\gamma}\big),
$$
and similarly
$$
F_{\gamma,\prec \mu}^{G,M}:=\bigcup_{\nu\in \Lambda_{M_{0},\prec \mu}^{H_{M}}}\xx_{\gamma}^{0}\big(\nu+\Sigma_{\gamma}\big).
$$ 
It is the union of finitely many $F_{\gamma,\preccurlyeq \mu'}^{G,M}$, $\mu'\in \Lambda_{M_{0}}^{H_{M}}$. Let 
$$
F_{\gamma,\,\mu}^{G,M}:=F_{\gamma,\preccurlyeq \mu}^{G,M}\backslash F_{\gamma,\prec \mu}^{G,M}.
$$
Being difference of closed sub-schemes, $F_{\gamma,\,\mu}^{G,M}$ is locally closed in $\xx_{\gamma}^{0}$. As $F_{\gamma,\preccurlyeq \mu}^{G,M}$ is semi-infinite unions of translations of the fundamental domains, they are all isomorphic; Similarly for $F_{\gamma,\prec \mu}^{G,M}$. Hence $F_{\gamma,\,\mu}^{G,M}$ are all isomorphic. Moreover,
$$
F_{\gamma,\preccurlyeq \mu}^{G,M}=\bigsqcup_{\nu\in \Lambda_{M_{0},\preccurlyeq \mu}^{H_{M}}} 
F_{\gamma,\,\nu}^{G,M}
$$
by induction, and 
$$
F_{\gamma}^{G,M}=\lim_{\mu\to \infty}
F_{\gamma,\preccurlyeq \mu}^{G,M}=\bigsqcup_{\nu\in \Lambda_{M_{0}}^{H_{M}}} 
F_{\gamma,\,\nu}^{G,M}.
$$
From all these we conclude:

\begin{prop}\label{reduce 1}

$F_{\gamma,\,\mu}^{G,M}$ is isomorphic to each other for all $\mu\in \Lambda_{M_{0}}^{H_{M}}$. In particular,
$$
|\big(\Lambda^{H_{M}}\backslash F_{\gamma}^{G,M}\big)(\fq)|=|\Lambda^{H_{M}}\backslash \Lambda_{M_{0}}^{H_{M}}|\cdot |F_{\gamma,\,\mu}^{G,M}(\fq)|.
$$

\end{prop}

Counting points on $F_{\gamma,\,\mu}^{G,M}$ can be reduced to the fundamental domains via a process similar to the Arthur-Kottwitz reduction. To begin with, notice that $\Pi_{\gamma,\preccurlyeq \mu}^{G,M}$ is bounded only in directions that are positive with respect to $P_{0}^{M}$. Indeed, its vertices are indexed by $P\in \cp(M_{0})$ satisfying $P\cap M=P_{0}^{M}$ and its faces indexed by $Q\in \cf(M_{0})$ such that  $Q\cap M\supset P_{0}^{M}$. Then, we define a semi-infinite polytope $\Pi_{\gamma,\preccurlyeq \mu^{-}}^{G,M}$ which is a translation of $\Pi_{\gamma,\preccurlyeq \mu}^{G,M}$: Let $\Delta^{M}$ be the set of simple roots in $\Phi(M,A)$ with respect to $B_{0}\cap M$, let $\Delta^{M}_{P_{0}^{M}}=\Delta^{M}\cap \Phi(N_{P_{0}^{M}},A)$, with $N_{P_{0}^{M}}$ being the unipotent radical of $P_{0}^{M}$. For $\alpha\in \Delta^{M}_{P_{0}^{M}}$, let $\{\omega_{\alpha}^{\vee}\}$ be the corresponding fundamental coweights. Let 
\begin{equation}\label{schrink mu}
\mu^{-}=\mu-\pi_{2}\Big(\textstyle\sum_{\alpha\in \Delta^{M}_{P_{0}^{M}}}\omega_{\alpha}^{\vee}\Big),
\end{equation} 
where $\pi_{2}$ is the projection to the second factor in the orthogonal decomposition $\ka_{A}^{G}=\ka_{A}^{M_{0}}\oplus \ka_{M_{0}}^{M}\oplus \ka_{M}^{G}$. Then $\Pi_{\gamma,\preccurlyeq \mu^{-}}^{G,M}$ is a translation of $\Pi_{\gamma,\preccurlyeq \mu}^{G,M}$ by the same vector. Now let $\varsigma\in \ka_{M_{0}}^{G}$ be a generic element such that $\alpha(\varsigma)$ is positive but almost equal to $0$ for any $\alpha\in \Delta_{P_{0}}$. 
We perturb the semi-infinite polytope $\Pi_{\gamma,\preccurlyeq \mu^{-}}^{G,M}$ to a similar one $\Pi'$, with vertices
$$
\lambda_{P}(\Pi')=\lambda_{P}\big(\Pi_{\gamma,\preccurlyeq \mu^{-}}^{G,M}\big)+w\cdot \varsigma,
\quad
\forall\,P\in \cp(M_{0}),\,P\cap M=P_{0}^{M},
$$
where $w\in W$ is any element satisfying $P=w\cdot P_{0}$. Both $\Pi_{\gamma,\preccurlyeq \mu}^{G,M}$ and $\Pi'$ can be seen as limits of positive $(G,M_{0})$-orthogonal families containing $\Sigma_{\gamma}$, hence we can apply an analogue of the Arthur-Kottwitz reduction to get a decomposition of the complement $\xx_{\gamma}^{0}(\Pi_{\gamma,\preccurlyeq \mu}^{G,M})\backslash \xx_{\gamma}^{0}(\Pi')$.
For $Q=LN_{Q}\in \cf(M_0)$ satisfying $Q\cap M\supset P_{0}^{M}$, define $R_{\Pi',\,Q}$ to be the subset of $ \ka_{M_0}^{G}$ satisfying conditions
\begin{eqnarray*}
\pi_{M_{0}}^{L}(a)&\subset& {\Pi'}^{Q};\\
\alpha(\pi_{M_{0},L}(a))&\geq& \alpha(\pi_{M_{0},L}(\lambda_{Q}(\Pi'))),\quad\forall\,\alpha\in \Delta_{Q}.
\end{eqnarray*}
This gives us a partition 
\begin{equation}\label{part AK pi}
\ka_{M_0}^{G}=\Pi'\;\cup \bigcup_{\substack{Q\in \cf(M_0),\,Q\neq G\\ Q\cap M\supset P_{0}^{M}}} R_{\Pi',\,Q}.
\end{equation} 
It induces a disjoint partition of $\Lambda_{M_{0}}$. For $G=\gl_{3}$, $\gamma$ split and $M=M_{\alpha_{12}}$, we get Fig. \ref{AK induction picture-semi-infinite}.

\begin{figure}[h]
\begin{center}
\begin{tikzpicture}[every node/.style={scale=0.8}]

\draw (4.75,-1.25)--(-2.25,-1.25);
\draw[dashed] (-2.25,-1.25)--(-4.25,-1.25);
\draw (-2.75,1.3)--(4.25,1.3);
\draw[dashed] (-4.75,1.3)--(-2.75,1.3);
\draw (4.25,1.3)--(5.25,-0.4);
\draw (5.25,-0.4)--(4.75,-1.25);

\draw (3.25,1.3)--(4.25,-0.4);
\draw (4.25,-0.4)--(3.75,-1.25);

\draw [red] (3.3,1.35)--(3.3,3.05);
\draw [red] (3.3,1.35)--(5.1,2.2);
\draw [red] (3.3,1.35)--(-2.85,1.35);

\draw [red] (4.33,-0.4)--(5.83,0.45);
\draw [red] (4.33,-0.4)--(5.83,-1.25);

\draw [red] (3.8,-1.3)--(5.3,-2.15);
\draw [red] (3.8,-1.3)--(3.8,-3);

\draw [red] (3.8,-1.3)--(4.33,-0.4);
\draw [red] (3.3,1.35)--(4.33,-0.4);
\draw [red] (3.8,-1.3)--(-2.35,-1.3);

\node [red] at (1, 0) {$\Pi'$};
\node [red] at (4.1, 2.55) {$R_{\Pi',B_{0}}$};

\end{tikzpicture}
\caption{Arthur-Kottwitz reduction for $\Pi_{\gamma,\preccurlyeq \mu}^{G,M}$.}
\label{AK induction picture-semi-infinite}
\end{center}
\end{figure}

Running the same construction as in \S \ref{AK red}. For $Q=LN_{Q}\in \cf(M_0),\,Q\neq G$ and $Q\cap M\supset P_{0}^{M}$, let 
$$
S_{\Pi',Q}:=\{x\in \xx_{\gamma}\mid \pi_{M_{0}}^{G}\big(\ec_{M_{0}}^{Q}(x)\big)\subset R_{\Pi', Q}\},
$$
and let 
$$
S_{\Pi',Q}^{\nu}=S_{\Pi',Q}\cap f_{Q}^{-1}(\xx_{\gamma}^{L,\nu}),\quad \forall\,\nu\in \Lambda_{L}.
$$ 
We get a disjoint partition
\begin{equation}\label{AK pi}
\xx_{\gamma}^{0}(\Pi_{\gamma,\preccurlyeq \mu}^{G,M})
=
\xx_{\gamma}^{0}(\Pi')\cup \bigcup_{\substack{Q\in \cf(M_0),\,Q\neq G\\ Q\cap M\supset P_{0}^{M}}}\bigcup_{\substack{\nu\in \Lambda^{0}_{L}\cap \pi_{L}(R_{\Pi',Q})
\\ \cap \pi_{L}(\Pi_{\gamma,\preccurlyeq \mu}^{G,M})}} 
S_{\Pi', Q}^{\nu}.
\end{equation}
The strata $S_{\Pi', Q}^{\nu}$ are locally closed sub-schemes of $\xx_{\gamma}$, and the retraction $f_{Q}: S_{\Pi',Q}^{\nu}\to \xx_{\gamma}^{L,\nu}({\Pi'}^{Q})$ is an iterated affine fibration over $\fq$ of dimension
$
\val(\det(\ad(\gamma\mid \kn_{Q,F}))).
$

\begin{prop}\label{reduce 2}
We have the equality
$$
F_{\gamma,\,\mu}^{G,M}=\bigcup_{\substack{Q\in \cf(M_0),\,Q\neq G\\ Q\cap M=P_{0}^{M}}}\bigcup_{\substack{\nu\in \Lambda^{0}_{L}\cap \pi_{L}(R_{\Pi',Q})
\\ \cap \pi_{L}(\Pi_{\gamma,\preccurlyeq \mu}^{G,M})}} 
S_{\Pi', Q}^{\nu}.
$$
Moreover, the index set $\Lambda^{0}_{L}\cap \pi_{L}(R_{\Pi',Q}) \cap \pi_{L}(\Pi_{\gamma,\preccurlyeq \mu}^{G,M})$ consists of at most one element, it is non-empty if and only if $Q$ is not contained in any $Q'\in \cf(M)$.

\end{prop}

\begin{proof}

For the first assertion. By construction, the points $x\in F_{\gamma,\,\mu}^{G,M}$ are characterized by the property  
$$
\ec_{M_{0}}^{Q}(x)\subset (\mu+\Sigma_{\gamma})^{Q}
$$ 
for some $Q\in \cf(M_0),\,Q\neq G,\, Q\cap M=P_{0}^{M}$. Since this is also the property characterizing points on the right hand side of the equality, we get the equality as claimed. The second assertion follows from the observation that $\Pi'$ is a slight expansion of $\Pi_{\gamma,\preccurlyeq \mu^{-}}^{G,M}$, hence the regions $R_{\Pi',Q}$, $Q$ contained in some maximal parabolic subgroup in $\cf(M)$, contain no elements in $\Lambda_{L}^{0}\cap \pi_{L}(\Pi_{\gamma,\preccurlyeq \mu}^{G,M})$.  

\end{proof}

If the index set $\Lambda^{0}_{L}\cap \pi_{L}(R_{\Pi',Q}) \cap \pi_{L}(\Pi_{\gamma,\preccurlyeq \mu}^{G,M})$ is non-empty, we denote by $\nu_{Q}$ the unique element in it. Let $\mu_{Q}\in \Lambda_{L}^{0}$ be the unique element such that $\alpha(\mu_{Q})=1$ for all $\alpha\in \Delta_{Q}$, then we have
$$
\xx_{\gamma}^{L,\nu_{Q}}({\Pi'}^{Q})\cong F_{\gamma}^{L,\mu_{Q}}.
$$
Combining with the fact that 
$$
f_{Q}: S_{\Pi',Q}^{\nu_{Q}}\to \xx_{\gamma}^{L,\nu_{Q}}({\Pi'}^{Q})
$$ 
is an iterated affine fibration, we get

\begin{cor}\label{reduce 3}

For $M\in \cl(M_{0}),\,M\neq M_{0}$, we have the equality
$$
|F_{\gamma,\,\mu}^{G,M}(\fq)|=
\sum_{\substack{Q=LN_{Q}\in \cf(M_0),\,Q\neq G\\ \text{satisfying }(*)}}
q^{\frac{1}{2}\val(\det(\ad\gamma|\kg_{F}/\mathfrak{l}_{F}))}\cdot |F_{\gamma}^{L,\mu_{Q}}(\fq)|,
$$
where $(*)$ refers to the conditon that $Q\cap M=P_{0}^{M}$ and $Q\nsubseteq Q',\,\forall\, Q'\in \cf(M)$. 
\end{cor}

Notice that the equation doesn't involve the fundamental domain $F_{\gamma}$. Together with corollary \ref{AK count} and proposition \ref{reduce 1}, we get an expression of $|\big(\Lambda^{H_{M}}\backslash \xx_{\gamma}^{\nu_{0}}(\Pi)\big)(\fq)|$ in terms of $|F_{\gamma}^{L,\mu_{Q}}(\fq)|$, $L\in \cl(M_{0}),\, L\neq G$. Recall that counting points on $F_{\gamma}^{L,\mu_{Q}}$ can be reduced to that of $F_{\gamma}^{L'}$, $L'\in \cl(M_{0}), L'\subset L$, by the Arthur-Kottwitz reduction, we get an expression of $|\big(\Lambda^{H_{M}}\backslash \xx_{\gamma}^{\nu_{0}}(\Pi)\big)(\fq)|$ in terms of $|F_{\gamma}^{L'}(\fq)|, L'\in \cl(M_{0})$.

\section{Counting points by Harder-Narasimhan reduction}\label{Section HN count}

The number of points $|\big(\Lambda^{H_{M}}\backslash\xx_{\gamma}^{\nu_{0}}(\Pi)\big)(\fq)|,\,\nu_{0}\in \Lambda_{G}$, can also be counted by the Harder-Narasimhan reduction. The comparison with results from last section gives us a recursive relation between Arthur's weighted orbital integrals and the number of rational points on the fundamental domains.

\subsection{Harder-Narasimhan reduction on the affine Springer fibers}\label{HN whole}

We have introduced a notion of $\xi$-stability on the affine Grassmannian and constructed the associated Harder-Narasimhan reduction in \cite{chen1}. 
In this section, we generalize it to a broader set-up. 
The following lemma is an analogue of \cite{cl2}, proposition 5.6.1. Let $S$ be an affine $\fq$-scheme and $x\in \xx(S)$. For every point $s\in S$, let $x_{s}\in \xx(k(s))$ be the base change of $x$ to the residue field $k(s)$ of S at $s$. Let $\mathcal{C}_{x}$ be the map on $S$ which sends every point $s\in S$ to the convex polytope $\ec(x_{s})$.

\begin{lem}\label{open nature of ec}

Suppose that $S$ is noetherian. The map $\mathcal{C}_{x}$ from $S$ to the set of convex polytopes in $\ka_{A}^{G}$ ordered by inclusion is lower semi-continuous. In other words, for any convex polytope $\Xi$, the set 
$$
\{s\in S\mid\mathcal{C}_{x}(s)\supset \Xi\}
$$
is open. 
\end{lem}

\begin{proof}

To begin with, we show that $\mathcal{C}_{x}$ is constructible and it takes only finitely many values. Passing to the irreducible components of $S$, we can suppose that $S$ is irreducible. Let $\eta$ be the generic point of $S$. Let $g_{\eta}\in G\big(k(\eta)(\!(\ep)\!)\big)$ be a representative of $x_{\eta}$. For $B=AN\in \cp(A)$, we have the Iwasawa decomposition 
$$
g_{\eta}=n_{\eta}a_{\eta}k_{\eta},
$$
where $n_{\eta}\in N\big(k(\eta)(\!(\ep)\!)\big)$, $a_{\eta}\in A\big(k(\eta)(\!(\ep)\!)\big)$ and $k_{\eta}\in G\big(k(\eta)[\![\ep]\!]\big)$. Because $\eta$ is the generic point and the map $\nu_{A}:\xx^{A}\to X_{*}(A)$ is essentially the valuation map, there exists an open sub-scheme $U$ of $S$ such that $H_{B}(x_{s})=\nu_{A}(a_{\eta})$ for any $x\in U$. As $\ec(x_{s})$ is the convex hull of $H_{B}(x_{s}), \,B\in \cp(A)$, the map $\mathcal{C}_{x}$ takes the constant value $\ec(x_{\eta})$ on the intersection of all such open sub-schemes $U$. This proves the constructibility of $\mathcal{C}_{x}$. By the noetherian induction, the map $\mathcal{C}_{x}$ takes only finitely many values.

To finish the proof, we only need to show that the map $\mathcal{C}_{x}$ decreases under specialisation. In other words, let $S$ be the spectrum of a discret valuation ring, let $s$ be its special point and $\eta$ its generic point, then
$$
\ec(x_{s})\subset \ec(x_{\eta}).
$$
This is equivalent to the assertion that 
\begin{equation}\label{specialisation inequality}
f_{B}(x_{s})\prec_{B} f_{B}(x_{\eta}), \quad \forall\, B\in \cp(A),
\end{equation}
where $\prec_{B}$ is the order on $X_{*}(A)$ such that
$\mu_{1}\prec_{B}\mu_{2}$ if and only if $ \mu_{2}-\mu_{1}$ is a positive linear combination of positive coroots with respect to $B$.

Let $\mu=f_{B}(x_{\eta})\in X_{*}(A)$. By definition, we have 
$$
x_{\eta}\in U_{B}(\!(\ep)\!)\ep^{\mu}G[\![\ep]\!]/G[\![\ep]\!],
$$
where $U_{B}$ is the unipotent radical of $B$. So
$$
x_{s}\in \overline{x_{\eta}}\subset \overline{U_{B}(\!(\ep)\!)\ep^{\mu}G[\![\ep]\!]/G[\![\ep]\!]}=\bigcup_{\substack{\lambda\in X_{*}(A)\\ \lambda\prec_{B}\mu}} 
U_{B}(\!(\ep)\!)\ep^{\lambda}G[\![\ep]\!]/G[\![\ep]\!],
$$
which implies the relation (\ref{specialisation inequality}).

\end{proof}

\begin{defn}

Let $\xi\in \ka_{M}^{G}$ be a generic element.
A point $x\in \xx$ is said to be $\xi$-\emph{stable} if the polytope $\pi^{G}\big(\ec_{M}(x)\big)$ contains $\xi$. 

\end{defn}

As $\ec_{M}(x)=\pi_{M}(\ec(x))$, the subset
$$
\xx^{\xi}=\big\{x\in \xx\mid \xi\in \pi^{G}\big(\ec_{M}(x)\big)\big\}
$$
is an open sub-ind-$\fq$-scheme of $\xx$ by lemma \ref{open nature of ec}. This been shown, all the other constructions of \cite{chen1} generalize.

\begin{rem}

When $M=A$, we recover the $\xi$-stability of \cite{chen1}. In that work, we prove that the notion of $\xi$-stability co\"incides with the notion of stability for a twisted action of $A$ on $\xx$. We believe that this holds also in the current setting with the torus $A_{M}$ playing the role of $A$. If this holds, we can conclude that the quotient $\xx^{\xi}/A_{M}$ exists as an ind-$\fq$-scheme.

\end{rem}

Harder-Narasimhan reduction works as well in this setting.
For $Q=LN_{Q}\in \cf(M)$, let $\Phi_Q(G,L)$ be the image of $\Phi(N_{Q},\,A)$ in $(\ka_L^G)^{*}$. For any point $a\in \ka_L^G$, we define a cone in $\ka_L^G$,
$$
D_Q(a)=\left\{y\in\ka_L^G\,|\, \alpha(y-a)\geq 0,\,\forall \alpha\in \Phi_Q(G,L)\right\}.
$$

\begin{defn}

For any geometric point $x\in \xx$, we define a semi-cylinder $C_Q(x)$ in $\ka_{M}^G$ by
$$
C_Q(x)=\pi_{M}^{L,-1}\big(\ec_{M}^{L}(f_{Q}(x))\big)\cap \pi_{M,L}^{-1}\big(D_Q(H_{Q}(x))\big).
$$
\end{defn}

By definition, we get a partition 
$$
\ka_{M}^{G}=\pi^{G}(\ec_{M}(x))\cup\bigcup_{\substack{Q\in \cf(M)\\Q\neq G}}C_Q(x),
$$ 
for which the interior of any two parts doesn't intersect. The picture is similar to Fig. \ref{partitiongl3}.
Hence for any $x\notin \xx^{\xi}$, there exists a unique parabolic subgroup  $Q\in \cf(M)$ such that $\xi\in C_Q(x)$ as $\xi$ is generic. In this case, $f_Q(x)\in \xx^L$ is $\xi^{L}$-stable, where $\xi^{L}=\pi_{M}^{L}(\xi)\in \ka_{M}^L$. Let
$$
X_Q=\{x\in \xx| \,\xi\in C_{Q}(x)\},
$$
we have the decomposition of the affine grassmannian
\begin{equation}\label{decomp1}
\xx=\xx^{\xi}\sqcup \bigsqcup_{\substack{Q\in \cf(M)\\Q\neq G}}X_Q.
\end{equation}

For $Q\in \cp(L)$, let $Q^{-}$ be the parabolic subgroup opposite to $Q$ with respect to $L$. Let $\Lambda^{\xi}_{L,Q}=(\pi_{L}^{G})^{-1}\big(D_{Q^{-}}(\xi_{L})\big)\cap \Lambda_{L}$, we have the disjoint partition
$$
\Lambda_L=\bigsqcup_{Q\in \cp(L)}\Lambda^{\xi}_{L,Q}.
$$
For $\lambda\in \Lambda_{L}$, let
$
\xx^{L,\lambda,\xi^L}=\xx^{L,\xi^L}\cap \xx^{L,\lambda}.
$
The stratum $X_{Q}$ can be further decomposed into $N_{Q}(\!(\ep)\!)$-orbits
$$
X_{Q}=\bigsqcup_{\lambda\in \Lambda^{\xi}_{L,Q}} N_{Q}(\!(\ep)\!)\xx^{L,\,\lambda,\,\xi^{L}}.
$$
Each orbit is locally closed in $\xx$, they are infinite dimensional homogeneous affine fibrations on $\xx^{L,\lambda,\xi^{L}}$ under the retraction $f_{Q}$. The above discussions can be summarized as:

\begin{thm}\label{red1}

The affine grassmannian can be decomposed as
$$
\xx=\xx^{\xi}\sqcup \bigsqcup_{\substack{Q=LN_{Q}\in \cf(M)\\ Q\neq G}}\bigsqcup_{\lambda\in \Lambda^{\xi}_{L,Q}}
N_{Q}(\!(\ep)\!)\xx^{L,\,\lambda,\,\xi^{L}}.
$$
Each stratum $N_{Q}(\!(\ep)\!)\xx^{L,\,\lambda,\,\xi^{L}}$ is an infinite dimensional homogeneous affine fibration over $\xx^{L,\,\lambda,\,\xi^{L}}$.

\end{thm}

Now that $\gamma\in \km(F)$, we can restrict the above constructions to $\xx_{\gamma}$. Let $\xx_{\gamma}^{\xi}=\xx_{\gamma}\cap \xx^{\xi}$, it is an open sub-scheme of $\xx_{\gamma}$. As $T(F)\xrightarrow{H_{M}} X_{*}(M)$ is surjective, the connected components of $\xx_{\gamma}^{\xi}$ can be translated to each other by elements of $T(F)$. Moreover, for different choices of generic element $\xi,\,\xi'\in \ka_{M}^{G}$, the corresponding $\xx_{\gamma}^{\xi},\,\xx_{\gamma}^{\xi'}$ can be translated to each other by elements of $T(F)$. Hence $\xx_{\gamma}^{\xi}$ doesn't depend on the choice of $\xi$.

The Harder-Narasimhan reduction restricts to
\begin{equation}\label{HN}
\xx_{\gamma}=\xx_{\gamma}^{\xi}\sqcup \bigsqcup_{\substack{Q=LN_{Q}\in \cf(M)\\ Q\neq G}}\bigsqcup_{\lambda\in \Lambda^{\xi}_{L,Q}}\big(\xx_{\gamma}\cap N_{Q}(\!(\ep)\!)\xx_{\gamma}^{L,\,\lambda,\,\xi^{L}}\big).
\end{equation}
By proposition \ref{KL retraction}, the retraction 
$$
f_{Q}:\xx_{\gamma}\cap N_{Q}(\!(\ep)\!)\xx_{\gamma}^{L,\,\lambda,\,\xi^{L}}\to \xx_{\gamma}^{L,\,\lambda,\,\xi^{L}}
$$ 
is an iterated affine fibration over $\fq$ of relative dimension
$
\val(\det(\ad(\gamma)|\kn_{Q}(F))).
$

Coming back to the weighted orbital integrals. With the definition for general reductive algebraic groups as explained in remark \ref{general weight factor}, proposition \ref{geom weighted orbital} can be reformulated as

\begin{prop}\label{cl geom} 

Let $\xi\in \ka_{M}^{G}$ be a generic element, then
$$
J_{M}^{\xi}(\gamma)=\mathrm{vol}_{dt}\big(\Lambda^{H_{M}}\backslash T(F)_{M}^{1}\big)^{-1}\cdot
 \big|\Lambda^{H_{M}}\backslash \big((\xx^{G_{\mathrm{der}}}\cap\xx_{\gamma}^{\xi})(\fq)\big)\big|.
$$
In particular, let $\xi_{0}\in \ka_{M_{0}}^{G}$ be a generic element, then
\begin{equation*}
J_{M_{0}}^{\xi_{0}}(\gamma)=\mathrm{vol}_{dt}\big(T(F)^{1}\big)^{-1}\cdot
 \big|(\xx^{G_{\mathrm{der}}}\cap\xx_{\gamma}^{\xi})(\fq)\big|.
\end{equation*}

\end{prop}

\begin{proof}

When $G$ is semisimple, the proposition is a reformulation of proposition \ref{geom weighted orbital}. 
The complexity arises when $G$ has non-trivial connected center.

As $T$ is totally ramified over $F$, with the exact sequence (\ref{kott 2}), we see that the morphism $T(F)\xrightarrow{\nu_{G}} \Lambda_{G}$ is surjective, hence $G(F)=T(F)G_{\mathrm{der}}(F)$, and so
\begin{eqnarray*}
J_{M}^{\xi}(\gamma)&=&
\int_{T(F)\backslash G(F)} \mathbbm{1}_{\kg(\co)}\big(\Ad(g)^{-1}\gamma\big)\rw_{M}^{\xi}(g)\frac{dg}{dt} \\
&=&
\int_{T_{G_{\mathrm{der}}}(F)\backslash G_{\mathrm{der}}(F)} \mathbbm{1}_{\kg(\co)}\big(\Ad(g)^{-1}\gamma\big)\rw_{M_{G_{\mathrm{der}}}}^{\xi}(g)\frac{dg}{dt},
\end{eqnarray*}
with $T_{G_{\mathrm{der}}}=T\cap G_{\mathrm{der}}$ and $M_{G_{\mathrm{der}}}=M\cap G_{\mathrm{der}}$.
Following calculations in proposition \ref{geom weighted orbital}, we get result similar to what we claim, with $\Lambda^{H_{M}}$ replaced by $\Lambda\cap G_{\mathrm{der}}(F)\cap \ker(H_{M_{G_{\mathrm{der}}}})$ and $T(F)_{M}^{1}$ replaced by $T_{G_{\mathrm{der}}}(F)_{M_{G_{\mathrm{der}}}}^{1}$. Notice that $\ker(H_{M})=M_{\mathrm{der}}(F)\cdot M(\co)$ by lemma 6.1 of \cite{cl1}, we have 
$$
\Lambda\cap G_{\mathrm{der}}(F)\cap \ker(H_{M_{G_{\mathrm{der}}}})=\Lambda\cap \ker(H_{M})=\Lambda^{H_{M}},
$$
and 
$$
T_{G_{\mathrm{der}}}(F)_{M_{G_{\mathrm{der}}}}^{1}=T(F)\cap G_{\mathrm{der}}(F)\cap \ker(H_{M_{G_{\mathrm{der}}}})=T(F)\cap \ker(H_{M})=T(F)_{M}^{1},
$$
and the proposition is proved.
\end{proof}

The volume factors have been calculated in equation (\ref{volume factor}).

\subsection{Harder-Narasimhan reduction for the truncated affine Springer fibers}\label{section HN truncated}

In contrast to the Arthur-Kottwitz reduction, the Harder-Narasimhan reduction doesn't work well on the truncated affine Springer fiber $\xx_{\gamma}(\Pi)$. We need to cut it into two parts, the \emph{tail} and the \emph{main body}. The Harder-Narasimhan reduction works well on the main body, and we can handle the tail with the Arthur-Kottwitz reduction.

For $Q\in \cf(M),\, Q\neq G$, we define the positive $(G,M)$-orthogonal family 
$
E_{Q}(\Pi)
$,
which as a polytope is the union of the translations $\pi^{G}(\Sigma_{\gamma}^{G,M}+\lambda)$, $\lambda\in \Lambda_{M}$, such that
$
\pi^{G}(\Sigma_{\gamma}^{G,M}+\lambda)^{Q}\subset \Pi^{Q}.
$
Let
$$
{}^{t}\xx_{\gamma}(\Pi)=\bigcup_{\substack{Q\in \cf(M)\\ Q\neq G}}\xx_{\gamma}(E_{Q}(\Pi)), \quad {}^{m}\xx_{\gamma}(\Pi)=\xx_{\gamma}(\Pi)\backslash ^{t}\xx_{\gamma}(\Pi).
$$
We call them the \emph{tail} and the \emph{main body} of $\xx_{\gamma}(\Pi)$, they are closed and open sub-schemes of $\xx_{\gamma}(\Pi)$ respectively. Fig. \ref{tail picture} gives an example of $E_{Q}(\Pi)$ for the group $G=\gl_{3}$ when $M=A$.

\begin{figure}[h]
\begin{center}
\begin{tikzpicture}[scale=0.5,every node/.style={scale=0.7}]

\draw (0.75,-1.3)--(-0.25,-1.3);
\draw (0.75,-1.3)--(1.25,-0.4);
\draw (1.25,-0.4)--(0.25,1.3);
\draw (0.25,1.3)--(-0.75,1.3);
\draw (-0.75,1.3)--(-1.25,0.4);
\draw (-1.25,0.4)--(-0.25,-1.3);

\draw (4.75,-8.1)--(9.25,-0.4);
\draw (9.25,-0.4)--(4.25,8.1);
\draw (4.25,8.1)--(-4.75,8.1);
\draw (-4.75,8.1)--(-9.25,0.4);
\draw (-9.25,0.4)--(-4.25,-8.1);
\draw (-4.25,-8.1)--(4.75,-8.1);

\draw [blue] (5.25,-7.25)--(4.25,-5.55);
\draw [blue] (4.25,-5.55)--(-4.75,-5.55);
\draw [blue] (-4.75,-5.55)--(-5.25,-6.4);

\draw [red] (3.75,-8.1)--(2.75,-6.4);
\draw [red] (2.75,-6.4)--(7.25,1.25);
\draw [red] (7.25,1.25)--(8.25,1.25);

\draw [blue] (8.75,-1.25)--(7.75,-1.25);
\draw [blue] (7.75,-1.25)--(2.75,7.25);
\draw [blue] (2.75,7.25)--(3.25,8.1);

\draw [red] (5.25,6.4)--(4.75,5.55);
\draw [red] (4.75,5.55)--(-4.25,5.55);
\draw [red] (-4.25,5.55)--(-5.25,7.25);

\draw [blue] (-3.75,8.1)--(-2.75,6.4);
\draw [blue] (-2.75,6.4)--(-7.25,-1.25);
\draw [blue] (-7.25,-1.25)--(-8.25,-1.25);

\draw [red] (-8.75,1.25)--(-7.75,1.25);
\draw [red] (-7.75,1.25)--(-2.75,-7.25);
\draw [red] (-2.75,-7.25)--(-3.25,-8.1);

\node [blue] at (4, 6.8) {$E_{B_{0}}(\Pi)$};
\node [blue] at (0, 6.8) {$E_{P_{1}}(\Pi)$};
\node [blue] at (6,3.6) {$E_{P_{2}}(\Pi)$};

\node [blue] at (0,0) {$\Sigma_{\gamma}$};

\end{tikzpicture}
\caption{$E_{P}(\Pi)$ for $\gl_{3}$ when $M_{0}=A$.}
\label{tail picture}
\end{center}
\end{figure}

Before proceeding, we make precise the condition of $\Pi$ being sufficiently regular. We would like it to satisfy the following conditions:
\begin{enumerate}
\item

$\Pi$ is $\Sigma_{\gamma}$-regular.

\item

For all $P,Q\in \cf(M)$, $E_{P}(\Pi)\cap E_{Q}(\Pi)=E_{P\cap Q}(\Pi)$.

\item

The complement $\Pi\backslash \bigcup_{\substack{Q\in \cf(M)\\ Q\neq G}}E_{Q}(\Pi)$ is a polytope associated to a positive $(G,M)$-orthogonal family, let $\Pi_{0}$ be a slight shrinking of it (The definition is similar to equation (\ref{slight expansion}), with the  sign ``$+$'' replaced by ``$-$''). 
We require that $\Pi_{0}$ is sufficiently large: For all $Q=LN_{Q}\in \cf(M)$, the face $\Pi_{0}^{Q}$ contains the translations of $\Sigma_{\gamma}^{Q}$ in $\ka_{M}^{L}$ which have $\xi^{L}$ as one of its vertices.

\end{enumerate}

\begin{rem}\label{Pi suff reg}

As $\Pi_{0}$ is convex, the condition (3) implies that for any $\nu\in \Lambda_{L,Q}^{\xi}\cap \pi_{L}(\Pi_{0})$ the intersection $\Pi_{0}\cap \pi_{M_{0},L}^{-1}(\nu)$ contains translations of $\Sigma_{\gamma}^{Q}$ in the hyperplane $\pi_{M_{0},L}^{-1}(\nu)$ which have $\xi^{L}$ as one of its vertices.
By definition of $\xi$-stability, this implies
$$
\xx_{\gamma}^{L,\nu,\xi^{L}}\subset \xx_{\gamma}^{L,\nu}\big(\Pi_{0}\cap \pi_{M_{0},L}^{-1}(\nu)\big),\quad \text{ for all } \nu\in \Lambda_{L,Q}^{\xi}\cap \pi_{L}(\Pi_{0}).
$$
Actually, this is the reason to impose condition (3).
\end{rem}


\subsubsection{The main body}

By definition, a Harder-Narasimhan stratum $N_{Q}(\!(\ep)\!)\xx_{\gamma}^{L,\,\nu,\,\xi^{L}}\cap\xx_{\gamma}$, $\nu\in \Lambda_{L,Q}^{\xi}$, intersects non-trivially with $\xx_{\gamma}(\Pi)$ if and only if $\nu\in \Lambda_{L,Q}^{\xi}\cap \pi_{L}(\Pi)$.
So, after restriction, the Harder-Narasimhan reduction becomes
\begin{equation*}
\xx_{\gamma}(\Pi)=\xx_{\gamma}^{\xi}\sqcup \bigsqcup_{\substack{Q=LN_{Q}\in \cf(M)\\ Q\neq G}}\bigsqcup_{\lambda\in \Lambda^{\xi}_{L,Q}\cap \pi_{L}(\Pi)}\big(\xx_{\gamma}(\Pi)\cap N_{Q}(\!(\ep)\!)\xx_{\gamma}^{L,\lambda,\xi^{L}}\big).
\end{equation*}
The problem is that the retraction
$$
f_{Q}: \xx_{\gamma}(\Pi)\cap N_{Q}(\!(\ep)\!)\xx_{\gamma}^{L,\lambda,\xi^{L}}
\to \xx_{\gamma}^{L,\lambda,\xi^{L}}
$$
is not necessarily an iterated affine fibration. 
This problem disappears on the main body ${}^{m}\xx_{\gamma}(\Pi)$. We begin by analyzing the polytope $\ec(x),\,x\in \xx_{\gamma}(\Pi)$.

\begin{lem}

For $x\in \xx_{\gamma}(\Pi)$, suppose that 
$$
\pi^{G}(\ec_{M}(x))\subset \bigcup_{\substack{Q\in \cf(M)\\ Q\neq G}} E_{Q}(\Pi),
$$ 
then $\pi^{G}(\ec_{M}(x))\subset E_{Q}(\Pi)$ for some $Q\in \cf(M),\,Q\neq G$.

\end{lem}

\begin{proof}

By proposition \ref{gkm bound}, it is enough to prove the lemma for $x\in \xx_{\gamma}^{\reg}$. In this case, the polytope $\pi^{G}(\ec_{M}(x))$ is a translation of $\Sigma_{\gamma}^{G,M}$.
As $\bigcup_{\substack{Q\in \cf(M)\\ Q\neq G}} E_{Q}(\Pi)$ is the union of translation of $\Sigma_{\gamma}^{G,M}$ along the facets of $\Pi$, there must be a maximal parabolic subgroup $Q\in \cf(M)_{\max}$ such that $\pi^{G}(\ec_{M}^{Q}(x))\subset \Pi^{Q}$. By definition, this means that $\pi^{G}(\ec_{M}(x))\subset E_{Q}(\Pi)$.

\end{proof}

\begin{lem}

Let $Q=LN\in \cf(M),\, \nu \in \Lambda_{L,Q}^{\xi}$. Suppose that 
$$
{}^{m}\xx_{\gamma}(\Pi)\cap N(\!(\ep)\!)\xx_{\gamma}^{L,\nu}\neq \emptyset,
$$ 
then $\nu\in \Lambda_{L,Q}^{\xi}\cap \pi_{L}(\Pi_{0})$.

\end{lem}

\begin{proof}

We only need to show that $\nu\in \pi_{L}(\Pi_{0})$. 
Let $x\in {}^{m}\xx_{\gamma}(\Pi)\cap N(\!(\ep)\!)\xx_{\gamma}^{L,\nu}$. As $x\in {}^{m}\xx_{\gamma}(\Pi)$, we have 
$$
\pi^{G}(\ec_{M}(x))\nsubseteq E_{Q'}(\Pi),\quad \forall\,Q'\in \cf(M),\,Q'\neq G. 
$$
By the previous lemma, this is equivalent to 
\begin{equation}\label{ec not included}
\pi^{G}(\ec_{M}(x))\nsubseteq \bigcup_{\substack{Q'\in \cf(M)\\Q'\neq G}}E_{Q'}(\Pi).
\end{equation}

Suppose that $\nu\notin \pi_{L}(\Pi_{0})$, then $\nu\in \pi_{L}(E_{Q_{0}}(\Pi))$ for some $Q_{0}\in \cf(M),\,Q_{0}\supset L$. As $\nu\in \Lambda_{L,Q}^{\xi}$, the parabolic subgroup $Q_{0}$ need to satisfy $Q_{0}\supset Q^{-}$. Now that $x\in N(\!(\ep)\!)\xx_{\gamma}^{L,\nu}$ and that $\ec_{L}(x)$ is a positive $(G,L)$-orthogonal family, we have
$$
\alpha(H_{Q'}(x)-\nu)\geq 0, \quad \forall\,\alpha\in \Delta_{Q},\,Q'\in \cp(L).
$$
As $Q_{0}\supset Q^{-}$, this implies that 
$$
H_{Q'}(x)\subset \pi_{L}(E_{Q_{0}}(\Pi)),\quad \forall\,Q'\in \cp(L).
$$
Hence $\pi_{L}(\ec_{M}(x))\subset \pi_{L}(E_{Q_{0}}(\Pi))$, so 
$$
\ec_{M}(x)\subset \pi_{L}^{-1}(\pi_{L}(E_{Q_{0}}(\Pi)))\subset\bigcup_{\substack{Q'\subset \cf(M),\,Q'\neq G\\Q'\cap Q_{0}\neq \emptyset}} E_{Q'}(\Pi).
$$
This is in contradiction with the relation (\ref{ec not included}), hence $\nu$ must lie in $\pi_{L}(\Pi_{0})$.

\end{proof}

Restricting the Harder-Narasimhan reduction (\ref{HN}) to the main part ${}^{m}\xx_{\gamma}(\Pi)$, we have
$$
{}^{m}\xx_{\gamma}(\Pi)=\xx_{\gamma}^{\xi}\sqcup \bigsqcup_{\substack{Q=LN_{Q}\in \cf(M)\\ Q\neq G}}\bigsqcup_{\lambda\in \Lambda^{\xi}_{L,Q}\cap \pi_{L}(\Pi_{0})}
\big({}^{m}\xx_{\gamma}(\Pi)\cap N_{Q}(\!(\ep)\!)\xx_{\gamma}^{L,\lambda,\xi^{L}}\big).
$$
The retraction $f_{Q}$ behaves much better on the stratum ${}^{m}\xx_{\gamma}(\Pi)\cap N_{Q}(\!(\ep)\!)\xx_{\gamma}^{L,\lambda,\xi^{L}}$:

\begin{prop}

Let $Q=LN_{Q}\in \cf(M),\, \nu \in \Lambda_{L,Q}^{\xi}\cap \pi_{L}(\Pi_{0})$. We have
$$
{}^{m}\xx_{\gamma}(\Pi)\cap N_{Q}(\!(\ep)\!)\xx_{\gamma}^{L,\nu,\xi^{L}}=\xx_{\gamma}\cap N_{Q}(\!(\ep)\!)\xx_{\gamma}^{L,\nu,\xi^{L}}.
$$
Hence the retraction 
$$
f_{Q}: {}^{m}\xx_{\gamma}(\Pi)\cap N_{Q}(\!(\ep)\!)\xx_{\gamma}^{L,\nu,\xi^{L}}\to \xx_{\gamma}^{L,\nu,\xi^{L}}
$$
is an iterated affine fibration over $\fq$.

\end{prop}

\begin{proof}

Notice that the second assertion is the corollary of the first one, as follows from proposition \ref{KL retraction}. It is thus enough to show the first one. In particular, it is enough to show
$$
\xx_{\gamma}\cap N_{Q}(\!(\ep)\!)\xx_{\gamma}^{L,\nu,\xi^{L}}\subset {}^{m}\xx_{\gamma}(\Pi)\cap N_{Q}(\!(\ep)\!)\xx_{\gamma}^{L,\nu,\xi^{L}},
$$
as the inclusion in the other direction is obvious.

Let $x\in \xx_{\gamma}\cap N_{Q}(\!(\ep)\!)\xx_{\gamma}^{L,\nu,\xi^{L}}$, 
we claim that $\ec_{M}(x)\subset \Pi$. 
According to remark \ref{Pi suff reg}, the condition (3) of $\Pi$ being sufficiently regular implies
\begin{equation}\label{ecm inclusion}
\ec_{M}^{L}(f_{Q}(x))\subset \Pi_{0}\cap \pi_{L}^{-1}(\nu), 
\end{equation}
because $f_{Q}(x)\in \xx_{\gamma}^{L,\nu,\xi^{L}}$.
This implies that $\ec_{M}(x)\subset \Pi$ by proposition \ref{gkm bound} because of the inclusion
$$
\bigcup_{\substack{\lambda\in \Lambda_{M}\\
\text{satisfying }(*)}}(\lambda+\Sigma_{\gamma}^{G,M})\subset \Pi,
$$
where the condition $(*)$ refers to 
$$
(\lambda+\Sigma_{\gamma}^{G,M})\cap \Pi_{0}\neq \emptyset.
$$
The inclusion (\ref{ecm inclusion}) also implies that 
$$
\ec_{M}(x)\nsubseteq \bigcup_{\substack{Q\in \cf(M)\\Q\neq G}}E_{Q}(\Pi).
$$
So $x\in {}^{m}\xx_{\gamma}(\Pi)$, and the proof is concluded.

\end{proof}

We summarize the above discussions in a proposition.

\begin{prop}\label{HN for main body}

The main body has a decomposition
$$
{}^{m}\xx_{\gamma}(\Pi)=\xx_{\gamma}^{\xi}\sqcup \bigsqcup_{\substack{Q=LN_{Q}\in \cf(M)\\ Q\neq G}}\bigsqcup_{\lambda\in \Lambda^{\xi}_{L,Q}\cap \pi_{L}(\Pi_{0})}
\big({}^{m}\xx_{\gamma}(\Pi)\cap N_{Q}(\!(\ep)\!)\xx_{\gamma}^{L,\lambda,\xi^{L}}\big),
$$
and the retraction $f_{Q}$ on each stratum 
$$
f_{Q}: {}^{m}\xx_{\gamma}(\Pi)\cap N_{Q}(\!(\ep)\!)\xx_{\gamma}^{L,\nu,\xi^{L}}\to \xx_{\gamma}^{L,\nu,\xi^{L}}
$$
is an iterated affine fibration over $\fq$ of dimension
$
\val(\det(\ad(\gamma)\mid \kn_{Q,F})).
$

\end{prop}

Of course, we can restrict the decomposition to each connected component ${}^{m}\xx_{\gamma}^{\nu_{0}}(\Pi)$, $\nu_{0}\in \Lambda_{G}$. 
Let $\Lambda_{L,Q}^{\nu_{0},\xi}=\Lambda_{L,Q}^{\xi}\cap \Lambda_{L}^{\nu_{0}}$,
the decomposition implies
\begin{eqnarray*}
&&|\big(\Lambda^{H_{M}}\backslash {}^{m}\xx_{\gamma}^{\nu_{0}}(\Pi)\big)(\fq)|\\
&=&
|\big(\Lambda^{H_{M}}\backslash\xx_{\gamma}^{\nu_{0},\xi}\big)(\fq)|+\sum_{\substack{Q=LN_{Q}\in \cf(M)\\ Q\neq G}}\sum_{\lambda\in \Lambda^{\nu_{0},\xi}_{L,Q}\cap \pi_{L}(\Pi_{0})} 
q^{\frac{1}{2}\val(\det(\ad(\gamma)\mid \kg_{F}/\mathfrak{l}_{F}))}\cdot |\big(\Lambda^{H_{M}}\backslash\xx^{L,\lambda,\xi^{L}}_{\gamma}\big)(\fq)|\nonumber \\
&=&
|\big(\Lambda^{H_{M}}\backslash\xx_{\gamma}^{0,\xi}\big)(\fq)|
+\sum_{\substack{Q=LN_{Q}\in \cf(M)\\ Q\neq G}}
q^{\frac{1}{2}\val(\det(\ad \gamma|\kg_{F}/\mathfrak{l}_{F}))}\cdot |\big(\Lambda^{H_{M}}\backslash\xx^{L,0,\xi^{L}}_{\gamma}\big)(\fq)|\cdot |\Lambda^{\nu_{0},\xi}_{L,Q}\cap \pi_{L}(\Pi_{0})|\nonumber\\
&=&
|\big(\Lambda^{H_{M}}\backslash\xx_{\gamma}^{0,\xi}\big)(\fq)|
+\sum_{\substack{L\in \cl(M)\\ L\neq G}}
q^{\frac{1}{2}\val(\det(\ad \gamma|\kg_{F}/\mathfrak{l}_{F}))}\cdot |\big(\Lambda^{H_{M}}\backslash\xx^{L,0,\xi^{L}}_{\gamma}\big)(\fq)|\cdot\sum_{Q\in \cp(L)} |\Lambda^{\nu_{0},\xi}_{L,Q}\cap \pi_{L}(\Pi_{0})|\nonumber\\
&=& 
|\big(\Lambda^{H_{M}}\backslash\xx_{\gamma}^{0,\xi}\big)(\fq)|
+\sum_{\substack{L\in \cl(M)\\ L\neq G}}
q^{\frac{1}{2}\val(\det(\ad \gamma|\kg_{F}/\mathfrak{l}_{F}))}\cdot |\big(\Lambda^{H_{M}}\backslash\xx^{L,0,\xi^{L}}_{\gamma}\big)(\fq)|\cdot|\Lambda^{\nu_{0}}_{L}\cap \pi_{L}(\Pi_{0})|. \nonumber
\end{eqnarray*}
Here for the second equality we have used the fact that all the connected components of $\xx_{\gamma}^{L,\,\xi_{L}}$ are isomorphic. Moreover, the last term in the equation counts the number of lattice points in a polytope, it can be calculated effectively with methods from toric geometry. In summary,

\begin{thm}\label{HN main body}

For any $\nu_{0}\in \Lambda_{G}$, the number of rational points on the main body is
\begin{eqnarray*}
|\big(\Lambda^{H_{M}}\backslash {}^{m}\xx_{\gamma}^{\nu_{0}}(\Pi)\big)(\fq)|=|\big(\Lambda^{H_{M}}\backslash\xx_{\gamma}^{0,\xi}\big)(\fq)|
+\sum_{\substack{L\in \cl(M)\\ L\neq G}}
q^{\frac{1}{2}\val(\det(\ad \gamma|\kg_{F}/\mathfrak{l}_{F}))}\cdot\\
 |\big(\Lambda^{H_{M}}\backslash\xx^{L,0,\xi^{L}}_{\gamma}\big)(\fq)|\cdot|\Lambda^{\nu_{0}}_{L}\cap \pi_{L}(\Pi_{0})|.
\end{eqnarray*}

\end{thm}

\subsubsection{The tail}\label{tail 1}

As the polytope $\Pi$ satisfies
$$
E_{P}(\Pi)\cap E_{Q}(\Pi)=E_{P\cap Q}(\Pi),\quad \forall\, P,Q\in \cf(M),
$$
by the inclusion-exclusion principle, we have
\begin{equation}\label{tail alternate}
|\big(\Lambda^{H_{M}}\backslash {}^{t}\xx_{\gamma}^{\nu_{0}}(\Pi)\big)(\fq)|=\sum_{\substack{Q\in \cf(M)\\ Q\neq G}} (-1)^{\rk(G)-\rk(Q)-1}|\big(\Lambda^{H_{M}}\backslash \xx_{\gamma}^{\nu_{0}}(E_{Q}(\Pi))\big)(\fq)|,
\end{equation}
where the notation $\rk$ means the semisimple rank.  
Although the polytope $E_{Q}(\Pi)$ is not $\Sigma_{\gamma}^{G,M}$-regular, we can use the general Arthur-Kottwitz reduction as explained in remark \ref{general AK} repeatedly to decompose $\xx_{\gamma}(E_{Q}(\Pi))$ into locally closed sub-schemes which are iterated affine fibrations over $F_{\gamma}^{L, M}$, $L\in \cl(M)$. 
This gives a formula for $|\big(\Lambda^{H_{M}}\backslash \xx_{\gamma}^{\nu_{0}}(E_{Q}(\Pi))\big)(\fq)|$ in terms of the $|\big(\Lambda^{H_{M}}\backslash F_{\gamma}^{L,M}\big)(\fq)|$'s, which can be further reduced to counting points on the fundamental domains by proposition \ref{reduce 1} and corollary \ref{reduce 3}. This process applies to a large family of truncated affine Springer fibers.

We introduce a family of operators on the set of all positive $(G,M)$-orthogonal families.
Recall that $Q_{0}=MU_{0}$ is the unique parabolic subgroup in $\cp(M)$ which contains $P_{0}$. 
For $L\in \cl(M)$, let $Q_{0}^{L}=Q_{0}\cap L$. For a positive $(L,M)$-orthogonal family, we say that two faces of it are \emph{conjugate} if their associated parabolic subgroups are conjugate to each other by the Weyl group $W_{L}/W_{M}$. In particular, the edges of the polytope are parametrized by minimal elements in $\cf^{L}(M)\backslash\{M\}$.
An edge is said to be \emph{of type} $\alpha\in \Delta_{A_{M}}^{Q_{0}^{L}}:=\Phi(L,A)\cap\Phi(U_{0},A)\cap \Delta$ if it is conjugate to the edge having vertices $\lambda_{Q^{L}_{0}},\lambda_{s_{\alpha}Q^{L}_{0}}$, where $s_{\alpha}$ is the simple reflection associated to $\alpha$.
Let $A_{M,\alpha}^{G,L}$ be the operator on the set of positive $(G,M)$-orthogonal families defined as follows: As a polytope, it increases by one the length of all the edges whose image in $\ka_{M}^{L}$ under the projection $\pi_{M}^{L}$ are of type $\alpha$, and keep the length of all the others invariant. To check that it actually sends a positive $(G,M)$-orthogonal family to another one, it suffices to verify for the faces of dimension $2$, but this is clear. Then we set the vertex
$$
\lambda_{Q_{0}}(A_{M,\alpha}^{G,L}(\Pi))=\lambda_{Q_{0}}(\Pi)+\frac{1}{2}\pi_{M}(\varpi_{\alpha}^{\vee}),
$$ 
to make $A_{M,\alpha}^{G,L}(\Pi)$ symmetric with respect to $\Pi$. Here $\varpi_{\alpha}^{\vee}$ is the fundamental coweight corresponding to $\alpha$.
By definition, we see that the operators $A_{M,\alpha}^{G,L}$ commute with each other. When $G=L$, we simplify the notation $A_{M,\alpha}^{G,G}$ to $A_{M,\alpha}^{G}$.


Given a tuple of non-negative integers $\underline{n}=(n_{\alpha}), \alpha\in \Delta_{A_{M}}^{Q_{0}^{L}}$, let
\begin{equation}\label{polytope}
\Sigma_{\gamma}^{\underline{n}}
=\prod_{\alpha\in \Delta_{A_{M}}^{Q_{0}^{L}}}(A_{M,\alpha}^{G,L})^{n_{\alpha}}(\Sigma_{\gamma}).
\end{equation}
It is easy to see that the polytopes $E_{Q}(\Pi)$ can be made by iterating this process. 
For $\alpha\in \Delta_{A_{M}}^{Q_{0}^{L}}$, let $1_{\alpha}$ be the tuple taking value $1$ at $\alpha$ and $0$ otherwhere. By remark \ref{general AK}, the Arthur-Kottwitz reduction works for the complement $\xx_{\gamma}(\Sigma_{\gamma}^{\underline{n}+1_{\alpha}})  \backslash   \xx_{\gamma}(\Sigma_{\gamma}^{\underline{n}})$. The process is completely the same as explained in \S \ref{AK red}, we don't repeat it here. The resulting strata are iterated affine fibrations over truncated affine Springer fibers of the form $\xx_{\gamma}^{L,\nu'}(\Sigma_{\gamma}^{L,(\underline{n}')})$, $L\in \cl(M),\nu'\in \Lambda_{L}$. Iterating this process,  $\xx_{\gamma}(\Sigma_{\gamma}^{\underline{n}})$ can be decomposed as disjoint union of locally closed sub-schemes, which are iterated affine fibrations over $F_{\gamma}^{L,M,\mu}$, $L\in \cl(M),\mu\in \Lambda_{L}$. In particular, counting points on $\Lambda^{H_{M}}\backslash \xx_{\gamma}(\Sigma_{\gamma}^{\underline{n}})$ can be reduced to counting points on $\Lambda^{H_{M}}\backslash F_{\gamma}^{L,M,\mu}$, which can be further reduced to counting points on the fundamental domains $F_{\gamma}^{L'}, L'\in \cl(M_{0}),$ as we have explained in \S\ref{intermediate}. This process applies to counting points on 
$
\Lambda^{H_{M}}\backslash \xx_{\gamma}^{\nu_{0}}(E_{Q}(\Pi)).
$
By equation (\ref{tail alternate}), it gives an expression of $|\big(\Lambda^{H_{M}}\backslash {}^{t}\xx_{\gamma}^{\nu_{0}}(\Pi)\big)(\fq)|$ in terms of $|F_{\gamma}^{L}(\fq)|, L\in \cl(M_{0})$.

\subsection{Application to Arthur's weighted orbital integral}

By theorem \ref{cl comparision} and proposition \ref{cl geom}, Arthur's weighted orbital integral $J_{M}(\gamma)$ calculates essentially $|\big(\Lambda^{H_{M}}\backslash \xx_{\gamma}^{0,\,\xi}\big)(\fq)|$, as $\xx^{G_{\mathrm{der}}}\cap \xx_{\gamma}^{\xi}$ is the union of $|\Lambda_{G_{\mathrm{der}}}|$-copies of $\xx_{\gamma}^{0,\,\xi}$. The two approaches in \S \ref{Section AK count} and \S \ref{Section HN count} to calculate $|\big(\Lambda^{H_{M}}\backslash \xx_{\gamma}^{0}(\Pi)\big)(\fq)|$ give us a recurrence relation involving $|\big(\Lambda^{H_{M}}\backslash F_{\gamma}^{L, M,\,\mu}\big)(\fq)|$ and $|\big(\Lambda^{H_{M}}\backslash \xx_{\gamma}^{L,\,0,\,\xi^{L}}\big)(\fq)|$, for $L\in \cl(M),\,\mu\in \Lambda_{L^{\mathrm{ad}}}$. If we are able to solve this recurrence relation, we will get an expression for $|\big(\Lambda^{H_{M}}\backslash \xx_{\gamma}^{0,\,\xi}\big)(\fq)|$ in terms of $|\big(\Lambda^{H_{M}}\backslash F_{\gamma}^{L, M,\,\mu}\big)(\fq)|$'s, which can be further reduced to counting points on fundamental domains as explained in \S\ref{intermediate}.

\section{Calculations for the group $\gl_{2}$}

Let $G=\gl_{2}$, let $\gamma\in \ggl_{2}(F)$ be a regular semi-simple integral element. Assume that $\mathrm{char}(k)>2$ and the splitting field of $\gamma$ is totally ramified over $F$.
The torus $T$ is isomorphic either to $F^{\times}\times F^{\times}$ or to $\res_{E/F}E^{\times}$, where $E$ is a separable totally ramified field extension over $F$ of degree $2$. We call elements $\gamma$ in these cases \emph{split} and \emph{anisotropic} respectively.

\subsection{Split elements}

We can take $T$ to be the maximal torus of $G$ of the diagonal matrices and $\gamma\in \kt(\co)$ a regular element. Let
$$
n=\val(\alpha_{12}(\gamma)),
$$
which we call the root valuation of $\gamma$. The dimension of the affine Springer fiber $\xx_{\gamma}$ is known to be
$$
\dim(\xx_{\gamma})=n.
$$

In the remaining of the section, we assume that $n\geq 1$, as the case $n=0$ reduces to the group $\gl_{1}$. 
Recall that we have calculated $F_{\gamma}$ in \cite{chen2}. Let $X_{*}(T)\cong \bz^{2}$ be the usual identification, let $(n,0)\in \bz^{2}$, let 
$$
\sch(n,0)=\overline{K \begin{pmatrix}\ep^{n}&\\ &1\end{pmatrix}K/K}.
$$
We have $F_{\gamma}\cong \sch(n,0)$, and its number of rational points is
$$
|F_{\gamma}(\fq)|=\sum_{i=0}^{n}q^{i},
$$
by the Bruhat-Tits decomposition of $\sch(n,0)$. As $\Lambda_{\pgl_{2}}=\bz/2$, $F_{\gamma}$ has only one variant $F_{\gamma}^{1}$, we can calculate its number of rational points to be
$$
F_{\gamma}^{1}(\fq)=\sum_{i=0}^{n-1}q^{i}.
$$

Let $a\in \bn$, let $\Pi$ be the positive $(G,T)$-orthogonal family defined by 
$$
\lambda_{w}(\Pi)=\lambda_{w}(\Sigma_{\gamma})+ w(a\alpha_{12}^{\vee}),\quad \forall\,w\in W.
$$
Assume that $a\gg 0$, then $\Pi$ is sufficiently regular in the sense of \S \ref{section HN truncated}. We can calculate easily
$$
Q_{\gamma}^{0}(a):=|\xx_{\gamma}^{0}(\Pi)(\fq)|=\sum_{i=0}^{n}q^{i}+2q^{n}a,
$$ 
by the Arthur-Kottwitz reduction. We see that $Q_{\gamma}^{0}(a)$ is polynomial in $a$. By theorem \ref{HN main body}, we have
$$
|{}^{m}\xx_{\gamma}^{0}(\Pi)(\fq)|=|\xx_{\gamma}^{0,\,\xi}(\fq)|+[2a-(n+1)]q^{n}.
$$
The tail is the disjoint union of two fundamental domains, so its number of rational points is
$$
|{}^{t}\xx_{\gamma}^{0}(\Pi)(\fq)|=2\sum_{i=0}^{n}q^{i}.
$$
Because
$$
|\xx_{\gamma}^{0}(\Pi)(\fq)|=|{}^{m}\xx_{\gamma}^{0}(\Pi)(\fq)|+|{}^{t}\xx_{\gamma}^{0}(\Pi)(\fq)|,
$$
we get the equation
$$
\sum_{i=0}^{n}q^{i}+2aq^{n}=|\xx_{\gamma}^{0,\,\xi}(\fq)|+[2a-(n+1)]q^{n}+2\sum_{i=0}^{n}q^{i}.
$$
Solving it, we get
\begin{equation}\label{gl2 xi}
|\xx_{\gamma}^{0,\,\xi}(\fq)|=nq^{n}-\sum_{i=0}^{n-1}q^{i}.
\end{equation}
Now that $T(F)^{1}=T(\co)=T(F)_{1}$ has volume $1$, by proposition \ref{cl geom}, we have
$$
J_{T}^{\xi}(\gamma)=|\xx_{\gamma}^{0,\,\xi}(\fq)|=nq^{n}-\sum_{i=0}^{n-1}q^{i}.
$$

On the other hand, we can use equation (\ref{geom orbital reduced}) to calculate easily the orbital integral 
$$
I^{G}_{\gamma}=q^{n}.
$$
Combined with theorem \ref{cl comparision}, the above calculations can be  summarized as:

\begin{thm}\label{count gl2}

Let $\gamma\in \ggl_{2}(F)$ be a regular semisimple integral element of root valuation $n$. It has orbital integral $
I^{G}_{\gamma}=q^{n}.
$
The number of rational points on $\xx_{\gamma}^{0}(\Pi)$ is 
$$
|\xx_{\gamma}^{0}(\Pi)(\fq)|=\sum_{i=0}^{n}q^{i}+2aq^{n},
$$
and Arthur's weighted orbital integral $J_{T}(\gamma)$ equals
$$
J_{T}(\gamma)=\vol\big(\ka_{T_{\SL_{2}}}^{\SL_{2}}/X_{*}(T_{\SL_{2}})\big)\cdot \Big[nq^{n}-\sum_{i=0}^{n-1}q^{i}\Big].
$$

\end{thm}

\subsection{Anisotropic elements}

In this case $E=\fq(\!(\ep^{\frac{1}{2}})\!)$. Suppose that $\gamma=a+ b\ep^{\frac{1}{2}}$ under the isomorphism $Z_{G(F)}(\gamma)\cong \res_{E/F}E^{\times}$, with $a,b\in \co$. Under the basis $\{\ep^{\frac{1}{2}},1\}$ of $E$ over $\fq(\!(\ep)\!)$, the element $\gamma$ is of the form
$$
\gamma=\begin{bmatrix}
a&b\\
b\ep&a
\end{bmatrix}.
$$
It is clear that the affine Springer fibers $\xx_{\gamma}$ and $\xx_{-a+\gamma}$ are isomorphic, so we can assume that $a=0$. Let $b=b_{0}\ep^{n},\,b_{0}\in \co^{\times}$, we can write
\begin{equation}\label{gamma sl2 ramified}
\gamma=\begin{bmatrix}
&b_{0}\ep^{n}\\
b_{0}\ep^{n+1}&
\end{bmatrix}.
\end{equation}
Put in this form, it has been shown by Goresky, Kottwitz and MacPherson \cite{gkm2} that $\xx_{\gamma}$ admits an affine paving which is induced by the standard Bruhat-Tits decomposition of the affine Grassmannian. More precisely, let $I$ be the standard Iwahori subgroup, i.e. it is the pre-image of $B_{0}$ under the reduction $G(\co)\to G$, then 
$$
\xx_{\gamma}=\bigsqcup_{(a_{1},a_{2})\in \bz^{2}} \xx_{\gamma}\cap I\begin{pmatrix}
\ep^{a_{1}}&\\
&\ep^{a_{2}}
\end{pmatrix}K/K,
$$
and each intersection, denoted $S_{\baa}$, is isomorphic to a standard affine space. We calculate that $S_{\baa}$ is not empty if and only if
$$
-(n+1)\leq a_{1}-a_{2} \leq n,
$$
and that
$$
\dim(S_{\baa})=\begin{cases}
a_{1}-a_{2}, &\text{ if } a_{1}\geq a_{2};\\
a_{2}-a_{1}-1, &\text{ if } a_{1}< a_{2}.
\end{cases}
$$
Notice that this is also the dimension of $I\begin{pmatrix}
\ep^{a_{1}}&\\
&\ep^{a_{2}}
\end{pmatrix}K/K$, so they must be the same. 
Summarize the above calculations, and notice that $T(F)^{1}=T(\co)=T(F)_{1}$ have volume 1, we get

\begin{thm}\label{gl2 cal ramified}

Let $\gamma$ be the matrix $(\ref{gamma sl2 ramified})$. For $(a_{1},a_{2})\in \bz^{2}$, we have
$$
\xx_{\gamma}\cap I\begin{pmatrix}
\ep^{a_{1}}&\\
&\ep^{a_{2}}
\end{pmatrix}K/K=\begin{cases}
I\begin{pmatrix}
\ep^{a_{1}}&\\
&\ep^{a_{2}}
\end{pmatrix}K/K, &\text{ if } -(n+1)\leq a_{1}- a_{2}\leq n;\\
\emptyset, &\text{ if not}.
\end{cases}$$
As a corollary, we have
$$
J_{G}(\gamma)=I^{G}_{\gamma}=|F_{\gamma}(\fq)|=\sum_{i=0}^{n} q^{i}.
$$

\end{thm}

\section{Calculations for $\gl_{3}$--Split case}

Let $G=\gl_{3}$, let $\gamma\in \ggl_{3}(F)$ be a regular semisimple integral element. Assume that $\mathrm{char}(k)>3$ and the splitting field of $\gamma$ is totally ramified over $F$. 
The torus $T$ is isomorphic to either $F^{\times}\times F^{\times}\times F^{\times}$, or $F^{\times}\times \res_{E_{2}/F}E_{2}^{\times}$, or $\res_{E_{3}/F}E_{3}^{\times}$, where $E_{2},\,E_{3}$ are separable totally ramified field extensions over $F$ of degree $2$ and $3$ respectively. We call elements $\gamma$ in these cases \emph{split}, \emph{mixed} and \emph{anisotropic} respectively. Notice that in all these cases $T(F)^{1}=T(\co)=T(F)_{1}$ have volume 1, hence by proposition \ref{cl geom} we have
$$
J_{M}^{\xi}(\gamma)=|\xx_{\gamma}^{0,\,\xi}(\fq)|
$$
and so
$$
J_{M}(\gamma)=\vol\big(\ka_{M_{\SL_{3}}}^{\SL_{3}}/X_{*}(M_{\SL_{3}})\big)
\cdot|\xx_{\gamma}^{0,\,\xi}(\fq)|
$$
by theorem \ref{cl comparision} and the remark following it.

In this section, we restrict ourselves to the split case.
After conjugation, we take $T$ to be the maximal torus of $G$ of the diagonal matrices. Then $M_{0}=T$ and the other proper Levi subgroups in $\cl(T)$ can be parametrized as follows:
For a nonempty subset $I\subsetneq \{1,2,3\}$, let $P_{I}$ be the parabolic subgroup of $G$ which stabilizes the flag
$$
\bigoplus_{i=1}^{3}\fq e_{i}\supsetneq \bigoplus_{i\notin I}\fq e_{i}\supsetneq \emptyset.
$$ 
Let $P_{I}=M_{I}N_{I}$ be the standard Levi factorization, we have $M_{I}\cong \gl_{2}\times \gl_{1}$. As $M_{I}=M_{I^{c}}$ with $I^{c}$ being the complement of $I$, it is enough to calculate $J_{T}(\gamma)$ and $J_{M_{\{i\}}}(\gamma),\,i=1,2,3$.

Let $\gamma\in \kt(\co)$ be a regular element. As we show in the appendix of \cite{chen}, up to conjugation by the Weyl group, we can suppose that
$$
\val(\alpha_{12}(\gamma))\leq\val(\alpha_{23}(\gamma)),\;
\val(\alpha_{13}(\gamma))=\val(\alpha_{12}(\gamma)).
$$
In this case, $\gamma$ is said to be \emph{in minimal form}, and we call 
$$
(n_{1},n_{2})=\big(\val(\alpha_{12}(\gamma)),\val(\alpha_{23}(\gamma))\big)$$ 
the \emph{root valuation} of $\gamma$. The dimension of the affine Springer fiber $\xx_{\gamma}$ is known to be
$$
\dim(\xx_{\gamma})=2n_{1}+n_{2}.
$$

In the remaining of the section, we assume that $n_{1}\geq 1$, as the case $n_{1}=0$ reduces to the group $\gl_{2}$. 
Recall that we have calculated the Poincar\'e polynomial of $F_{\gamma}$ in \cite{chen2}.

\begin{prop}
The fundamental domain $F_{\gamma}$ admits an affine paving, its Poincaré polynomial depends only on the root valuation $(n_{1},n_{2})$, and it is 
\begin{eqnarray*}
P_{(n_{1},n_{2})}(t)&=&\sum_{i=1}^{n_{1}}i(t^{4i-2}+t^{4i-4})+\sum_{i=2n_{1}}^{n_{1}+n_{2}-1}(2n_{1}+1)t^{2i}\\
&&+\sum_{i=n_{1}+n_{2}}^{2n_{1}+n_{2}-1}4(2n_{1}+n_{2}-i)t^{2i}+t^{4n_{1}+2n_{2}}.
\end{eqnarray*}
In particular, $
|F_{\gamma}(\fq)|=P_{(n_{1},n_{2})}(q^{1/2}).
$\end{prop}

\subsection{Calculation of $J_{T}(\gamma)$}

Let $(a_{1},a_{2})\in \bn^{2}$, let $\Pi$ be the positive $(G,T)$-orthogonal family defined by 
$$
\lambda_{w}(\Pi)=\lambda_{w}(\Sigma_{\gamma})+ w\sum_{i=1}^{2}a_{i}\alpha_{i}^{\vee},\quad \forall\,w\in W.
$$
Assume that $\Pi$ is sufficiently regular in the sense of \S \ref{section HN truncated}, which means that $a_{1},\,a_{2}\gg 0$ and 
$$
2a_{1}-a_{2}>0,\quad 2a_{2}-a_{1}>0.
$$
We will calculate
$$
Q_{\gamma}^{0}(a_{1},a_{2}):=|\xx_{\gamma}^{0}(\Pi)(\fq)|
$$ 
following the two approaches that we have explained, and draw conclusions on Arthur's weighted orbital integral.

\subsubsection{Counting points by Arthur-Kottwitz reduction}

We will work out each term in corollary \ref{AK count}. Look at the summands indexed by the Borel subgroups. Each stratum contributes $q^{2n_{1}+n_{2}}$, so it remains to count the number of lattice points 
$$
\sum_{B\in \cp(T)}|\Lambda_{T}^{0}\cap R_{B}\cap \Pi|=6|\Lambda_{T}^{0}\cap R_{B_{0}}\cap \Pi|,
$$ 
where the equality is due to the symmetry of $\Pi$ with respect to $\Sigma_{\gamma}$.
We identify 
$$
\Lambda_{T}^{0}\cong \{(m_{1},m_{2},m_{3})\in \bz^{3}\mid m_{1}+m_{2}+m_{3}=0\}
$$ 
in the usual way. Let
$\overline{\ka}_{B_{0}}^{G}=\{a\in \ka_{T}^{G}\mid \alpha_{1}(a)\geq 0,\,\alpha_{2}(a)\geq 0\}$, let
$$
R_{0}=\{a\in \overline{\ka}_{B_{0}}^{G}\mid \varpi_{1}(a)\leq a_{1}-1,\,\varpi_{2}(a)\leq a_{2}-1\}.
$$  
Up to a suitable translation, we have 
$$
|\Lambda_{T}^{0}\cap R_{B_{0}}\cap \Pi|=|\Lambda_{T}^{0}\cap R_{0}|.
$$
We can express it as the difference of two lattice point counting problems. Let 
\begin{eqnarray*}
R_{1}
&=&
\{a\in \overline{\ka}_{B_{0}}^{G}|\,\varpi_{1}(a)\leq a_{1}-1,\,\varpi_{2}(a)\leq 2(a_{1}-1)\},\\
R_{2}
&=&
\{a\in \overline{\ka}_{B_{0}}^{G}|\,\varpi_{1}(a)\leq a_{1}-1,\,\varpi_{2}(a)\geq a_{2}\},
\end{eqnarray*}
then we have
$$
|\Lambda_{T}^{0}\cap R_{0}|=|\Lambda_{T}^{0}\cap R_{1}|-|\Lambda_{T}^{0}\cap R_{2}|.
$$
We count $|\Lambda_{T}^{0}\cap R_{1}|$ as follows:
\begin{eqnarray*}
|\Lambda_{T}^{0}\cap R_{1}|
&=&
\sum_{n=0}^{+\infty} |R_{1}\cap \{\mu\in \Lambda_{T}^{0}|\,\varpi_{1}(\mu)=n\}|\\
&=&
\sum_{i=1}^{\lfloor \frac {a_{1}}{2}\rfloor}[(3i-2)+(3i-1)]+\frac{1-(-1)^{a_{1}}}{2}\bigg(1+\frac{3(a_{1}-1)}{2}\bigg)\\
&=&
3\big\lfloor \frac {a_{1}}{2}\big\rfloor^{2}+\frac{1}{4}(1-(-1)^{a_{1}})(3a_{1}-1),
\end{eqnarray*}
where $\lfloor x \rfloor$ means the largest integer that is less than or equal to $x$. Similarly, we have
\begin{eqnarray*}
|\Lambda_{T}^{0}\cap R_{2}|
&=&
\sum_{n=2(a_{1}-1)}^{-\infty} |R_{2}\cap \{\mu\in \Lambda_{T}^{0}|\,\varpi_{2}(\mu)=n\}|\\
&=&
\sum_{i=1}^{\lfloor \frac {2a_{1}-a_{2}-1}{2}\rfloor}(i+i)+\frac{1+(-1)^{2a_{1}-a_{2}}}{2}\cdot \frac{2a_{1}-a_{2}}{2}\\
&=& \big\lfloor \frac {2a_{1}-a_{2}-1}{2}\big\rfloor\left(\big\lfloor \frac {2a_{1}-a_{2}-1}{2}\big\rfloor+1\right)+\frac{1}{4}(1+(-1)^{2a_{1}-a_{2}})(2a_{1}-a_{2}).
\end{eqnarray*}
In summary, the summands in corollary \ref{AK count} indexed by the Borel subgroups contributes 
\begin{multline}\label{Borel contribution}
6q^{2n_{1}+n_{2}}\Bigg[
3\big\lfloor \frac {a_{1}}{2}\big\rfloor^{2}
-\big\lfloor \frac {2a_{1}-a_{2}-1}{2}\big\rfloor\left(\big\lfloor \frac {2a_{1}-a_{2}-1}{2}\big\rfloor+1\right)
\\
+\frac{1}{4}(1-(-1)^{a_{1}})(3a_{1}-1)-\frac{1}{4}(1+(-1)^{2a_{1}-a_{2}})(2a_{1}-a_{2})
\Bigg].
\end{multline}

Now we calculate the contributions of the summands indexed by the maximal parabolic subgroups. They are parametrized at the beginning of the section by nonempty subsets $I\subsetneq \{1,2,3\}$. For $\mu\in \Lambda_{M_{I}^{\ad}}\cong \bz/2$, let $q_{I}^{\mu}=|F_{\gamma}^{M_{I},\nu}(\fq)|$, for any $\nu\in \Lambda_{M_{I}}$ which projects to $\mu\in \Lambda_{M_{I}^{\ad}}$. Let $\alpha_{I}$ be the unique element in $\Phi_{B_{0}\cap M_{I}}(M_{I},T)$, a simple calculation with the affine Springer fibers for the group $\gl_{2}$ shows that 
$$
q_{I}^{(0)}=\sum_{i=0}^{\val(\alpha_{I}(\gamma))}q^{i},\quad q_{I}^{(1)}=\sum_{i=0}^{\val(\alpha_{I}(\gamma))-1}q^{i}.
$$

For $I=\{i\},\,i=1,2,3$, it is easy to see that
\begin{eqnarray*}
|\Lambda_{M_{I}}^{0}\cap \pi_{M_{I}}(R_{P_{I}})\cap \pi_{M_{I}}(\Pi)\cap c_{M}^{-1}(0)|
&=&
\big\lfloor \frac{a_{1}}{2}\big\rfloor,\\
|\Lambda_{M_{I}}^{0}\cap \pi_{M_{I}}(R_{P_{I}})\cap \pi_{M_{I}}(\Pi)\cap c_{M}^{-1}(1)|&=&\big\lfloor \frac{a_{1}+1}{2}\big\rfloor.
\end{eqnarray*}
The summands indexed by $P_{I}$ in corollary \ref{AK count} with $|I|=1$ contributes in total
\begin{equation}\label{parabolic cont 1}
\big\lfloor \frac{a_{1}}{2}\big\rfloor\left(q^{2n_{1}}\sum_{i=0}^{n_{2}}q^{i}+2q^{n_{1}+n_{2}}\sum_{i=0}^{n_{1}}q^{i}
\right)+
\big\lfloor \frac{a_{1}+1}{2}\big\rfloor\left(q^{2n_{1}}\sum_{i=0}^{n_{2}-1}q^{i}+2q^{n_{1}+n_{2}}\sum_{i=0}^{n_{1}-1}q^{i}
\right).
\end{equation}
Similarly, the summands indexed by $P_{I}$ with $|I|=2$ contributes in total
\begin{equation}\label{parabolic cont 2}
\big\lfloor \frac{a_{2}}{2}\big\rfloor\left(q^{2n_{1}}\sum_{i=0}^{n_{2}}q^{i}+2q^{n_{1}+n_{2}}\sum_{i=0}^{n_{1}}q^{i}
\right)
+
\big\lfloor \frac{a_{2}+1}{2}\big\rfloor\left(q^{2n_{1}}\sum_{i=0}^{n_{2}-1}q^{i}+2q^{n_{1}+n_{2}}\sum_{i=0}^{n_{1}-1}q^{i}
\right).
\end{equation}

Summing up the contributions from equations (\ref{Borel contribution}), (\ref{parabolic cont 1}), (\ref{parabolic cont 2}), we obtain

\begin{prop}\label{truncation count}
We have 
\begin{eqnarray*}
Q_{\gamma}^{0}(a_{1},a_{2})&=&\sum_{i=1}^{n_{1}}i(q^{2i-1}+q^{2i-2})+\sum_{i=2n_{1}}^{n_{1}+n_{2}-1}(2n_{1}+1)q^{i}\\
&+&\sum_{i=n_{1}+n_{2}}^{2n_{1}+n_{2}-1}4(2n_{1}+n_{2}-i)q^{i}+q^{2n_{1}+n_{2}}\\
&+&6q^{2n_{1}+n_{2}}\Bigg[
3\big\lfloor \frac {a_{1}}{2}\big\rfloor^{2}
-\big\lfloor \frac {2a_{1}-a_{2}-1}{2}\big\rfloor\left(\big\lfloor \frac {2a_{1}-a_{2}-1}{2}\big\rfloor+1\right)
\\
&&+\frac{1}{4}(1-(-1)^{a_{1}})(3a_{1}-1)-\frac{1}{4}(1+(-1)^{2a_{1}-a_{2}})(2a_{1}-a_{2})
\Bigg]\\
&+&\left(\big\lfloor \frac{a_{1}}{2}\big\rfloor+\big\lfloor \frac{a_{2}}{2}\big\rfloor\right)\left(q^{2n_{1}}\sum_{i=0}^{n_{2}}q^{i}+2q^{n_{1}+n_{2}}\sum_{i=0}^{n_{1}}q^{i}
\right)\\
&+&
\left(\big\lfloor \frac{a_{1}+1}{2}\big\rfloor+\big\lfloor \frac{a_{2}+1}{2}\big\rfloor\right)\left(q^{2n_{1}}\sum_{i=0}^{n_{2}-1}q^{i}+2q^{n_{1}+n_{2}}\sum_{i=0}^{n_{1}-1}q^{i}
\right)
\end{eqnarray*}

In particular, it depends quasi-polynomially on $(a_{1},a_{2})$.

\end{prop}

\subsubsection{Counting points by Harder-Narasimhan reduction}

We begin by counting points on the main body, we need to work out each term in theorem \ref{HN main body}. For $L=T$, it is easy to see that $|\xx_{\gamma}^{T,\,0,\,\xi^{T}}(\fq)|=1$, and we need to count the number of lattice points in
$
\Lambda_{T}^{0}\cap \Pi_{0}.
$
Notice that for this we can shrink $\Pi_{0}$ to the convex hull of 
$
\Lambda_{T}^{0}\cap \Pi_{0}
$,
we conserve the notation $\Pi_{0}$ for the shrunk polytope.
In our work \cite{chen2}, \S6, we calculate $\ec(x_{0})$ for a particular choice of regular point $x_{0}\in \xx^{\reg}_{\gamma}$, we can adapt the result to our current setting. Let $(\sigma_{1}\,\sigma_{2}\,\sigma_{3})$ be the permutation sending $(1\,2\,3)$ to $(\sigma_{1}\,\sigma_{2}\,\sigma_{3})$, the vertices of $\Sigma_{\gamma}$ are: 
\begin{eqnarray*}
\lambda_{123}(\Sigma_{\gamma})&=&(0,0,0),\quad \lambda_{321}(\Sigma_{\gamma})=(-2n_{1},n_{1}-n_{2},n_{1}+n_{2}),\\
\lambda_{213}(\Sigma_{\gamma})&=&(-n_{1},n_{1},0), \quad \lambda_{312}(\Sigma_{\gamma})=(-n_{1},-n_{2},n_{1}+n_{2}),\\
\lambda_{132}(\Sigma_{\gamma})&=&(0,-n_{2},n_{2}), \quad \lambda_{231}(\Sigma_{\gamma})=(-2n_{1},n_{1},n_{1}),
\end{eqnarray*}
The vertices of $\Pi_{0}$ can be calculated to be
\begin{eqnarray*}
\lambda_{123}(\Pi_{0})&=&(a_{1}-2n_{1}-1, a_{2}-a_{1}+n_{1}-n_{2}, -a_{2}+n_{1}+n_{2}+1),\\
\lambda_{321}(\Pi_{0})&=&(-a_{2}+1, a_{2}-a_{1},a_{1}-1),\\
\lambda_{213}(\Pi_{0})&=&(a_{2}-a_{1}-n_{1},a_{1}-n_{2}-1, -a_{2}+n_{1}+n_{2}+1),\\
\lambda_{312}(\Pi_{0})&=&(a_{2}-a_{1}-n_{1}, -a_{2}+n_{1}+1,a_{1}-1),\\
\lambda_{132}(\Pi_{0})&=&(a_{1}-2n_{1}-1, -a_{2}+n_{1}+1, a_{2}-a_{1}+n_{1}), \\ 
\lambda_{231}(\Pi_{0})&=&(-a_{2}+1, a_{1}-n_{2}-1, a_{2}-a_{1}+n_{2}).
\end{eqnarray*}
We will count the lattice points in $\Pi_{0}$ indirectly. We complete the hexagon $\Pi_{0}$ to a triangle $T_{0}$, whose vertices are
\begin{eqnarray*}
\lambda_{123}(T_{0})&=& \lambda_{132}(T_{0})
=(2a_{2}-2n_{1}-n_{2}-2, -a_{2}+n_{1}+1, -a_{2}+n_{1}+n_{2}+1),\\
\lambda_{321}(T_{0})&=&\lambda_{312}(T_{0})
=(-a_{2}+1, -a_{2}+n_{1}+1, 2a_{2}-n_{1}-2),\\
\lambda_{213}(T_{0})&=&\lambda_{231}(T_{0})
=(-a_{2}+1,2a_{2}-2-n_{1}-n_{2},-a_{2}+n_{1}+n_{2}+1).
\end{eqnarray*}
Let $T_{1}\cup T_{2}\cup T_{3}$ be the complement of $\Pi_{0}$ in $T_{0}$, as shown in Fig. \ref{completion}. Notice that the $T_{i}$'s don't contain their common boundary with $\Pi_{0}$, so 
$$
|\Lambda_{T}^{0}\cap \Pi_{0}|=|\Lambda_{T}^{0}\cap T_{0}|-\sum_{i=1}^{3}|\Lambda_{T}^{0}\cap T_{i}|.
$$
The right hand side is much easier to calculate.

\begin{figure}[h]
\begin{center}
\begin{tikzpicture}[node distance = 2cm, scale=0.6]
\draw (-5.5, 3.17)--(5.5,3.17);
\draw (-5.5, 3.17)--(0,-6.35);
\draw (0,-6.35)--(5.5,3.17);

\draw (-4,0.57)--(-2.5,3.17);
\draw (-1.5,-3.75)--(1.5,-3.75);
\draw (0.5,3.17)--(3,-1.15);

\node [blue] at (-0.5,0) {$\Pi_{0}$};
\node [blue] at (-4, 2.3) {$T_{1}$};
\node [blue] at (3.5, 1.8) {$T_{2}$};
\node [blue] at (0, -4.7){$T_{3}$};

\end{tikzpicture}
\caption{Complete the hexagon to a triangle.}
\label{completion}
\end{center}
\end{figure}

The length of the edges of $T_{0}$ is $3a_{2}-3-2n_{1}-n_{2}$, so
\begin{eqnarray*}
|\Lambda_{T}^{0}\cap T_{0}|&=& \sum_{i=1}^{3a_{2}-3-2n_{1}-n_{2}+1} i\\
&=& \frac{1}{2}(3a_{2}-2-2n_{1}-n_{2})(3a_{2}-1-2n_{1}-n_{2}).
\end{eqnarray*}

The length of the edges of $T_{1}$ is $2a_{2}-a_{1}-n_{1}-1$. As we don't count the lattice points on the common boundary of $T_{1}$ and $\Pi_{0}$, we have
\begin{eqnarray*}
|\Lambda_{T}^{0}\cap T_{1}|&=& \sum_{i=1}^{2a_{2}-a_{1}-n_{1}-1} i\\
&=& \frac{1}{2}(2a_{2}-a_{1}-n_{1}-1)(2a_{2}-a_{1}-n_{1}).
\end{eqnarray*}

Similarly, the length of the edges of $T_{2}$ is $2a_{2}-a_{1}-n_{2}-1$ and we have
\begin{eqnarray*}
|\Lambda_{T}^{0}\cap T_{2}|
= \frac{1}{2}(2a_{2}-a_{1}-n_{2}-1)(2a_{2}-a_{1}-n_{2}).
\end{eqnarray*}

The triangle $T_{3}$ is of the same size as $T_{1}$, so 
\begin{eqnarray*}
|\Lambda_{T}^{0}\cap T_{3}|
= \frac{1}{2}(2a_{2}-a_{1}-n_{1}-1)(2a_{2}-a_{1}-n_{1}).
\end{eqnarray*}

Finally, 
\begin{eqnarray}\label{stable Borel}
|\Lambda_{T}^{0}\cap \Pi_{0}|&=&|\Lambda_{T}^{0}\cap T_{0}|-\sum_{i=1}^{3}|\Lambda_{T}^{0}\cap T_{i}|\\
&=&\frac{1}{2}(3a_{2}-2-2n_{1}-n_{2})(3a_{2}-1-2n_{1}-n_{2})\nonumber\\
&&-(2a_{2}-a_{1}-n_{1}-1)(2a_{2}-a_{1}-n_{1})\nonumber\\
&&-\frac{1}{2}(2a_{2}-a_{1}-n_{2}-1)(2a_{2}-a_{1}-n_{2}).\nonumber
\end{eqnarray}

We go on to calculate $|\Lambda_{L}^{0}\cap \pi_{L}(\Pi_{0})|$ for the other Levi subgroups $L\in \cl(T)$. Let $d_{L}$ be the distance between the facets $\Pi_{0}^{Q}$ and $\Pi_{0}^{Q^{-}}$, where $\cp(L)=\{Q,Q^{-}\}$. It is easy to see that
$$
|\Lambda_{L}^{0}\cap \pi_{L}(\Pi_{0})|=d_{L}+1.
$$
The set $\cl(T)\backslash \{T,G\}$ consists of $3$ elements, they are Levi factors $M_{\{i\}}$ of the parabolic subgroups $P_{\{i\}},\,i=1,2,3$. Use the explicit expression of the vertices of $\Pi_{0}$, we can calculate
\begin{eqnarray}
\label{para1}
|\Lambda_{M_{\{1\}}}^{0}\cap \pi_{M_{\{1\}}}(\Pi_{0})|
&=& 
\label{para2} 
d_{M_{\{1\}}}+1=a_{1}+a_{2}-2n_{1}-1,\\
|\Lambda_{M_{\{2\}}}^{0}\cap \pi_{M_{\{2\}}}(\Pi_{0})|
&=& 
d_{M_{\{2\}}}+1=a_{1}+a_{2}-n_{1}-n_{2}-1,\\
\label{para3}
|\Lambda_{M_{\{3\}}}^{0}\cap \pi_{M_{\{3\}}}(\Pi_{0})|
&=& d_{M_{\{3\}}}+1=a_{1}+a_{2}-n_{1}-n_{2}-1.
\end{eqnarray}
  
Now that $|\xx_{\gamma}^{L,\,0,\,\xi^{L}}(\fq)|$ has been calculated in equation (\ref{gl2 xi}), we can insert the equations (\ref{stable Borel}) to (\ref{para3}) to the equation in theorem \ref{HN main body}, and get

\begin{prop}\label{count main body gl3}

The number of rational points on the main body is 
\begin{eqnarray*}
|{}^{m}\xx_{\gamma}^{0}(\Pi)(\fq)|
&=&
|\xx_{\gamma}^{0,\,\xi}(\fq)|+ 
q^{2n_{1}+n_{2}}
\Bigg[\frac{1}{2}(3a_{2}-2-2n_{1}-n_{2})(3a_{2}-1-2n_{1}-n_{2})\\
&&-(2a_{2}-a_{1}-n_{1}-1)(2a_{2}-a_{1}-n_{1})\\
&&-\frac{1}{2}(2a_{2}-a_{1}-n_{2}-1)(2a_{2}-a_{1}-n_{2})\Bigg ]\\
&+&q^{2n_{1}} (a_{1}+a_{2}-2n_{1}-1)
\left(n_{2}q^{n_{2}}-\sum_{i=0}^{n_{2}-1}q^{i}\right)\\
&+&2q^{n_{1}+n_{2}}(a_{1}+a_{2}-n_{1}-n_{2}-1)
\left(n_{1}q^{n_{1}}-\sum_{i=0}^{n_{1}-1}q^{i}\right).
\end{eqnarray*}

\end{prop}

Now we proceed to counting points on the tail. To begin with, we write down the vertices of $\Pi$.
\begin{eqnarray*}
\lambda_{123}(\Pi)&=&(a_{1}, a_{2}-a_{1}, -a_{2}),\\
\lambda_{321}(\Pi)&=&(-a_{2}-2n_{1}, a_{2}-a_{1}+n_{1}-n_{2},a_{1}+n_{1}+n_{2}),\\
\lambda_{213}(\Pi)&=&(a_{2}-a_{1}-n_{1},a_{1}+n_{1}, -a_{2}),\\
\lambda_{312}(\Pi)&=&(a_{2}-a_{1}-n_{1}, -a_{2}-n_{2},a_{1}+n_{1}+n_{2}),\\
\lambda_{132}(\Pi)&=&(a_{1}, -a_{2}-n_{2}, a_{2}-a_{1}+n_{2}), \\ 
\lambda_{231}(\Pi)&=&(-a_{2}-2n_{1}, a_{1}+n_{1}, a_{2}-a_{1}+n_{1}).
\end{eqnarray*}
For nonempty subset $I\subsetneq \{1,2,3\}$, we simplify the notation $E_{P_{I}}(\Pi)$ to $E_{I}(\Pi)$. Using the coordinates of vertices of $\Pi$, we can calculate the lengths of the edges of $\Pi$, and find the following expression for $E_{I}(\Pi)$:
When $|I|=1$, we have
\begin{equation}\label{boundary induction 1}
E_{I}(\Pi)=\left(A_{\alpha_{I}}^{G,M_{I}}\right)^{2a_{2}-a_{1}}(\Sigma_{\gamma});
\end{equation}
When $|I|=2$, we have
\begin{equation}\label{boundary induction 2}
E_{I}(\Pi)=\left(A_{\alpha_{I}}^{G,M_{I}}\right)^{2a_{1}-a_{2}}(\Sigma_{\gamma}).
\end{equation}
As explained before, we can use the Arthur-Kottwitz reduction inductively to count the number of rational points on $\xx_{\gamma}^{0}(E_{I}(\Pi))$. We give the details for $I=\{3\}$, the others can be calculated in the same way.

Applying Arthur-Kottwitz reduction to pass from $\left(A_{\alpha_{\{3\}}}^{G,M_{\{3\}}}\right)^{a}(\Sigma_{\gamma})$ to $\left(A_{\alpha_{\{3\}}}^{G,M_{\{3\}}}\right)^{a+1}(\Sigma_{\gamma})$, the picture is similar to Fig.\ref{AK induction picture-semi-infinite}, we obtain
\begin{eqnarray*}
\bigg|\xx_{\gamma}^{0}\left(\left(A_{\alpha_{\{3\}}}^{G,M_{\{3\}}}\right)^{a+1}(\Sigma_{\gamma})\right)(\fq)\bigg|
&=&
\bigg|\xx_{\gamma}^{0}\left(\left(A_{\alpha_{\{3\}}}^{G,M_{\{3\}}}\right)^{a}(\Sigma_{\gamma})\right)(\fq)\bigg|
+q^{2n_{1}}|F_{\gamma}^{M_{\{1\}},1}(\fq)|\\
&& +q^{n_{1}+n_{2}}|F_{\gamma}^{M_{\{13\}},1}(\fq)|+q^{2n_{1}+n_{2}}\\
&=&
\bigg|\xx_{\gamma}^{0}\left(\left(A_{\alpha_{\{3\}}}^{G,M_{\{3\}}}\right)^{a}(\Sigma_{\gamma})\right)(\fq)\bigg|
+q^{2n_{1}}\sum_{i=0}^{n_{2}-1}q^{i}\\
&& +q^{n_{1}+n_{2}}\sum_{i=0}^{n_{1}-1}q^{i}+q^{2n_{1}+n_{2}}.
\end{eqnarray*}

From this relation and the equation (\ref{boundary induction 1}), we deduce that
\begin{eqnarray}\label{boundary count 1}
|\xx_{\gamma}^{0}(E_{\{3\}}(\Pi))(\fq)|&=&|F_{\gamma}(\fq)|+(2a_{2}-a_{1})\Big(q^{2n_{1}}\sum_{i=0}^{n_{2}-1}q^{i}\\
&& +q^{n_{1}+n_{2}}\sum_{i=0}^{n_{1}-1}q^{i}+q^{2n_{1}+n_{2}}\Big).\nonumber
\end{eqnarray}

Similarly, we have 
\begin{eqnarray}\label{boundary count 2}
|\xx_{\gamma}^{0}(E_{\{1\}}(\Pi))(\fq)|&=&|F_{\gamma}(\fq)|+(2a_{2}-a_{1})\Big(
2q^{n_{1}+n_{2}}\sum_{i=0}^{n_{1}-1}q^{i}\\
&&+q^{2n_{1}+n_{2}}\Big),\nonumber
\end{eqnarray}

\begin{eqnarray}\label{boundary count 3}
|\xx_{\gamma}^{0}(E_{\{2\}}(\Pi))(\fq)|&=&|F_{\gamma}(\fq)|+(2a_{2}-a_{1})\Big(q^{2n_{1}}\sum_{i=0}^{n_{2}-1}q^{i}\\
&& +q^{n_{1}+n_{2}}\sum_{i=0}^{n_{1}-1}q^{i}+q^{2n_{1}+n_{2}}\Big),\nonumber
\end{eqnarray}

\begin{eqnarray}\label{boundary count 4}
|\xx_{\gamma}^{0}(E_{\{12\}}(\Pi))(\fq)|&=&|F_{\gamma}(\fq)|+(2a_{1}-a_{2})\Big(q^{2n_{1}}\sum_{i=0}^{n_{2}-1}q^{i}\\
&& +q^{n_{1}+n_{2}}\sum_{i=0}^{n_{1}-1}q^{i}+q^{2n_{1}+n_{2}}\Big),\nonumber
\end{eqnarray}

\begin{eqnarray}\label{boundary count 5}
|\xx_{\gamma}^{0}(E_{\{23\}}(\Pi))(\fq)|&=&|F_{\gamma}(\fq)|+(2a_{1}-a_{2})\Big(
2q^{n_{1}+n_{2}}\sum_{i=0}^{n_{1}-1}q^{i}\\
&&+q^{2n_{1}+n_{2}}\Big),\nonumber
\end{eqnarray}

\begin{eqnarray}\label{boundary count 6}
|\xx_{\gamma}^{0}(E_{\{13\}}(\Pi))(\fq)|&=&|F_{\gamma}(\fq)|+(2a_{1}-a_{2})\Big(q^{2n_{1}}\sum_{i=0}^{n_{2}-1}q^{i}\\
&& +q^{n_{1}+n_{2}}\sum_{i=0}^{n_{1}-1}q^{i}+q^{2n_{1}+n_{2}}\Big).\nonumber
\end{eqnarray}

Inserting the equations (\ref{boundary count 1}) to (\ref{boundary count 6}) to the equation (\ref{tail alternate}), we get

\begin{prop}\label{count tail gl3}

The number of rational points on the tail equals
$$
|{}^{t}\xx_{\gamma}^{0}(\Pi)(\fq)|=(a_{1}+a_{2})\Big(2q^{2n_{1}}\sum_{i=0}^{n_{2}-1}q^{i}
+4q^{n_{1}+n_{2}}\sum_{i=0}^{n_{1}-1}q^{i}+3q^{2n_{1}+n_{2}}
\Big)
$$

\end{prop}

The sum of results in proposition \ref{count main body gl3}, \ref{count tail gl3} gives us another expression for $Q_{\gamma}^{0}(a_{1},a_{2})$.

\begin{cor}\label{truncation count 2}
\begin{eqnarray*}
Q_{\gamma}^{0}(a_{1},a_{2})
&=&
|\xx_{\gamma}^{0,\,\xi}(\fq)|+ 
q^{2n_{1}+n_{2}}
\Bigg[\frac{1}{2}(3a_{2}-2-2n_{1}-n_{2})(3a_{2}-1-2n_{1}-n_{2})\\
&&-(2a_{2}-a_{1}-n_{1}-1)(2a_{2}-a_{1}-n_{1})\\
&&-\frac{1}{2}(2a_{2}-a_{1}-n_{2}-1)(2a_{2}-a_{1}-n_{2})\Bigg ]\\
&+&q^{2n_{1}} (a_{1}+a_{2}-2n_{1}-1)
\left(n_{2}q^{n_{2}}-\sum_{i=0}^{n_{2}-1}q^{i}\right)\\
&+&2q^{n_{1}+n_{2}}(a_{1}+a_{2}-n_{1}-n_{2}-1)
\left(n_{1}q^{n_{1}}-\sum_{i=0}^{n_{1}-1}q^{i}\right)\\
&+&(a_{1}+a_{2})\Big(2q^{2n_{1}}\sum_{i=0}^{n_{2}-1}q^{i}
+4q^{n_{1}+n_{2}}\sum_{i=0}^{n_{1}-1}q^{i}+3q^{2n_{1}+n_{2}}\Big)
\end{eqnarray*}
\end{cor}

In particular, this shows that $Q_{\gamma}^{0}(a_{1},a_{2})$ depends polynomially on $(a_{1},a_{2})\in \bn^{2}$. As a corollary, the expression for $Q_{\gamma}^{0}(a_{1},a_{2})$ in proposition \ref{truncation count} is also a polynomial in $(a_{1},a_{2})$, although it doesn't seem to be so.

\subsubsection{Arthur's weighted orbital integral}

Now we can compare the two expressions in proposition \ref{truncation count} and corollary \ref{truncation count 2} for $Q_{\gamma}^{0}(a_{1},a_{2})$. Look at their constant terms $Q_{\gamma}^{0}(0,0)$, as $J_{T}^{\xi}(\gamma)=|\xx_{\gamma}^{0,\,\xi}(\fq)|$ in this case, we obtain 

\begin{thm}\label{integral gl3}
Chaudouard-Laumon's weighted orbital integral for $\gamma$ equals
\begin{eqnarray*}
J_{T}^{\xi}(\gamma)&=&\sum_{i=1}^{n_{1}}i(q^{2i-1}+q^{2i-2})
+\sum_{i=n_{1}+n_{2}}^{2n_{1}+n_{2}-1}
(4n_{1}+2n_{2}-4i-3)q^{i}
+(n_{1}^{2}+2n_{1}n_{2})q^{2n_{1}+n_{2}}.
\end{eqnarray*}
 
\end{thm}

By theorem \ref{cl comparision} and the remark following it, we get Arthur's weighted orbital integral as well. For the orbital integral $I^{G}_{\gamma}$, as $T$ is split, we can  calculate it easily by equation (\ref{geom orbital reduced}):
$$
I^{G}_{\gamma}=q^{2n_{1}+n_{2}}.
$$

\subsection{Calculation of $J_{M_{1}}(\gamma)$}

We parametrize the Levi groups as before, with the further simplification $M_{i}:=M_{\{i\}}$.
Let $\gamma=\diag(\gamma_{1},\gamma_{2},\gamma_{3})$ and $\gamma'=\diag(\gamma_{1},\gamma_{3},\gamma_{2})$, notice that 
$$
J_{M_{2}}(\gamma)=J_{M_{3}}(\gamma').
$$
Moreover, $\xx_{\gamma}$ and $\xx_{\gamma'}$ have the same geometry as they have the same root valuation (indeed, they have the same affine paving), which implies that
$$
J_{M_{2}}(\gamma)=J_{M_{3}}(\gamma')=J_{M_{3}}(\gamma).
$$
Hence it is enough to calculate $J_{M_{1}}(\gamma)$ and $J_{M_{3}}(\gamma)$. Notice that $M_{1}$ corresponds to the root $\alpha_{2}$ and $M_{3}$ to the root $\alpha_{1}$.

As usual, we identify $X_{*}(T)\cong \bz^{3}$ and $\ka_{T}^{G}$ with the hyperplane $x_{1}+x_{2}+x_{3}=0$ of $\ka_{T}=X_{*}(T)\otimes \br\cong\br^{3}$. The subspace $\ka_{T}^{M_{1}}\subset \ka_{T}^{G}$ becomes the line $\{(0,x,-x)\}$ and the subspace $\ka_{M_{1}}^{G}\subset \ka_{T}^{G}$ becomes $\{(-x,x/2,x/2)\}$. The lattice $\Lambda_{M_{1}}^{0}$ is identified with $\bz$ by the mapping
$$
\Lambda_{M_{1}}^{0}\to \bz:\quad \big(-(a+b),a,b\big)\mapsto a+b.
$$
Its inclusion in $\ka_{M_{1}}^{G}$ is described by the mapping
$$
\Lambda_{M_{1}}^{0}\to \ka_{M_{1}}^{G}:\quad \big(-(a+b),a,b\big)\mapsto \big(-(a+b),(a+b)/2, (a+b)/2\big).
$$
We identify $\ka_{M_{1}}^{G}$ with $\br$ by identifying $(-x,x/2,x/2)$ with $x$, the inclusion $\Lambda_{M_{1}}^{0}\subset \ka_{M_{1}}^{G}$ becomes the natural embedding $\bz\subset \br$. On the other hand, the discrete free abelian group $\Lambda\cong X_{*}(T)$ is naturally identified with $\bz^{3}$, the morphism $H_{M_{1}}:\Lambda\to \ka_{M_{1}}$ can be calculated to be
$$
H_{M_{1}}(a_{1},a_{2},a_{3})=\left(a_{1}, \frac{a_{2}+a_{3}}{2}, \frac{a_{2}+a_{3}}{2} \right).
$$
Hence $\Lambda^{H_{M_{1}}}$ is freely generated by the element $\diag(1,\ep,\ep^{-1})$.

According to proposition \ref{gkm bound}, we can take $\Sigma_{\gamma}^{G,M_{1}}$ to be the interval $[0,2n_{1}]$ in $\ka_{M_{1}}^{G}\cong\br$. For $N\in \bn,\,N\gg 0$, let $\Pi_{N}$ be the interval $[-N,2n_{1}+N]$, regarded as a $(G,M_{1})$-orthogonal family in $\ka_{M_{1}}^{G}$,
we are going to calculate $\big(\Lambda^{H_{M_{1}}}\backslash \xx_{\gamma}(\Pi_{N})\big)(\fq)$ by the two approaches described above.

In the Arthur-Kottwitz approach, we need to calculate $$|\big(\Lambda^{H_{M_{1}}}\backslash F_{\gamma}^{G,M_{1}}\big)(\fq)|\; \text{ and }\;|\big(\Lambda^{H_{M_{1}}}\backslash F_{\gamma}^{M_{1},M_{1}}\big)(\fq)|.$$
Combining proposition \ref{reduce 1} and corollary \ref{reduce 3}, we get
\begin{equation*}
|\big(\Lambda^{H_{M_{1}}}\backslash F_{\gamma}^{G,M_{1}}\big)(\fq)|=|F_{\gamma,\,\mu}^{G,M_{1}}(\fq)|=q^{2n_{1}+n_{2}}+2q^{n_{1}+n_{2}}(1+q+\cdots+q^{n_{1}-1}).
\end{equation*}
For the second one, since $F_{\gamma}^{M_{1},M_{1}}=\xx_{\gamma}^{M_{1},(0)}$ and $\Lambda^{H_{M_{1}}}=\Lambda^{M_{1}}$, we have
\begin{equation*}
|\big(\Lambda^{H_{M_{1}}}\backslash F_{\gamma}^{M_{1},M_{1}}\big)(\fq)|=|\big(\Lambda^{M_{1}}\backslash \xx_{\gamma}^{M_{1},0}\big)(\fq)|=q^{n_{2}}.
\end{equation*}
The reduction process is illustrated by a figure similar to Fig.\ref{AK induction picture-semi-infinite}. By corollary \ref{AK count}, we have
\begin{prop}\label{beta AK}
\begin{eqnarray*}
|\big(\Lambda^{H_{M_{1}}}\backslash \xx_{\gamma}^{0}(\Pi_{N})\big)(\fq)|&=&|\big(\Lambda^{H_{M_{1}}}\backslash F_{\gamma}^{G,M_{1}}\big)(\fq)|+2Nq^{2n_{1}}|\big(\Lambda^{H_{M_{1}}}\backslash F_{\gamma}^{M_{1},M_{1}}\big)(\fq)|\nonumber\\
&=& q^{2n_{1}+n_{2}}+2q^{n_{1}+n_{2}}(1+q+\cdots+q^{n_{1}-1})+2Nq^{2n_{1}+n_{2}}\nonumber \\
&=&(2N+1)q^{2n_{1}+n_{2}}+2q^{n_{1}+n_{2}}(1+q+\cdots+q^{n_{1}-1}).\end{eqnarray*}

\end{prop}

In the Harder-Narasimhan approach, we begin with counting points on the tail. By construction,
\begin{eqnarray}\label{beta tail}
|\big(\Lambda^{H_{M_{1}}}\backslash ^{t}\xx_{\gamma}^{0}(\Pi_{N})\big)(\fq)|&=&2|\big(\Lambda^{H_{M_{1}}}\backslash F_{\gamma}^{G,M_{1}}\big)(\fq)|\nonumber\\
&=& 2[q^{2n_{1}+n_{2}}+2q^{n_{1}+n_{2}}(1+q+\cdots+q^{n_{1}-1})].
\end{eqnarray}
Then, we calculate
\begin{equation*}
|\big(\Lambda^{H_{M_{1}}}\backslash \xx_{\gamma}^{M_{1},0,\,\xi^{M_{1}}}\big)(\fq)|=|\big(\Lambda^{H_{M_{1}}}\backslash \xx_{\gamma}^{M_{1},0}\big)(\fq)|=q^{n_{2}}.
\end{equation*}
By theorem \ref{HN main body}, this implies
\begin{eqnarray}\label{beta main}
|\big(\Lambda^{H_{M_{1}}}\backslash ^{m}\xx_{\gamma}^{0}(\Pi_{N})\big)(\fq)|&=&|\big(\Lambda^{H_{M_{1}}}\backslash \xx_{\gamma}^{\xi,0}\big)(\fq)|+[2N-(2n_{1}+1)]\cdot q^{2n_{1}}\cdot\nonumber\\
&& |\big(\Lambda^{H_{M_{1}}}\backslash \xx_{\gamma}^{M_{1},0,\,\xi^{M_{1}}}\big)(\fq)|\nonumber \\
&=&
|\big(\Lambda^{H_{M_{1}}}\backslash \xx_{\gamma}^{\xi,0}\big)(\fq)|+[2N-(2n_{1}+1)]\cdot q^{2n_{1}+n_{2}}.
\end{eqnarray}

Combining the equation (\ref{beta tail}) and (\ref{beta main}), we obtain
\begin{prop}\label{beta HN}
\begin{eqnarray*}
|\big(\Lambda^{H_{M_{1}}}\backslash \xx_{\gamma}^{0}(\Pi_{N})\big)(\fq)|&=&|\big(\Lambda^{H_{M_{1}}}\backslash \xx_{\gamma}^{\xi,0}\big)(\fq)|+[2N-(2n_{1}+1)]\cdot q^{2n_{1}+n_{2}}
\nonumber \\
&&+2[q^{2n_{1}+n_{2}}+2q^{n_{1}+n_{2}}(1+q+\cdots+q^{n_{1}-1})].
\end{eqnarray*}
\end{prop}

Comparing proposition \ref{beta AK} and \ref{beta HN}, we get

\begin{prop}\label{beta stable}
\begin{equation*}
|\big(\Lambda^{H_{M_{1}}}\backslash \xx_{\gamma}^{\xi,0}\big)(\fq)|=
2n_{1}q^{2n_{1}+n_{2}}-2q^{n_{1}+n_{2}}(1+q+\cdots+q^{n_{1}-1}).
\end{equation*}
\end{prop}

It remains to calculate the volume factor $\vol_{dt}\big(\Lambda^{H_{M_{1}}}\backslash T(F)_{M_{1}}^{1}\big)$. By equation (\ref{volume factor}), it equals $1$ because $S=T$ and the morphism $H_{M_{1}}: X_{*}(T)\to X_{*}(M)$ is surjective. The above calculations can be summarized as:

\begin{thm}\label{beta}

We have 
$$
J_{M_{1}}^{\xi}(\gamma)=|\big(\Lambda^{H_{M_{1}}}\backslash \xx_{\gamma}^{\xi,0}\big)(\fq)|=
2n_{1}q^{2n_{1}+n_{2}}-2q^{n_{1}+n_{2}}(1+q+\cdots+q^{n_{1}-1}).
$$

\end{thm}

\subsection{Calculation of $J_{M_{3}}(\gamma)$}

We make identifications as before. The subspace $\ka_{T}^{M_{3}}\subset \ka_{T}^{G}$ becomes the line $\{(x,-x,0)\}$ and the subspace $\ka_{M_{3}}^{G}\subset \ka_{T}^{G}$ becomes $\{(x/2,x/2,-x)\}$, they are identified with $\br$ as before. The lattice $\Lambda_{M_{3}}^{0}$ is identified with $\bz$ by the mapping
$$
\Lambda_{M_{3}}^{0}\to \bz:\quad \big(a,b,-(a+b)\big)\mapsto a+b.
$$
The inclusion $\Lambda_{M_{3}}^{0}\subset \ka_{M_{3}}^{G}$ becomes again the natural embedding $\bz\subset \br$. Similarly, the group $\Lambda^{H_{M_{3}}}$ is freely generated by the element $\diag(\ep,\ep^{-1},1)$.

By proposition \ref{gkm bound}, we can take $\Sigma_{\gamma}^{G,M_{3}}$ to be the interval $[0,n_{1}+n_{2}]$ in $\ka_{M_{3}}^{G}\cong\br$. For $N\in \bn,\,N\gg 0$, let $\Pi_{N}$ be the interval $[-N,n_{1}+n_{2}+N]$, regarded as a $(G,M_{3})$-orthogonal family in $\ka_{M_{3}}^{G}$.
We calculate $\big(\Lambda^{H_{M_{3}}}\backslash \xx_{\gamma}(\Pi_{N})\big)(\fq)$ in two ways as before.

Similar calculation as above, we get
\begin{equation*}
|\big(\Lambda^{H_{M_{3}}}\backslash F_{\gamma}^{M_{3},M_{3}}\big)(\fq)|=|\big(\Lambda^{M_{3}}\backslash \xx_{\gamma}^{M_{3},0}\big)(\fq)|=q^{n_{1}},
\end{equation*}
and
\begin{eqnarray*}
|\big(\Lambda^{H_{M_{3}}}\backslash F_{\gamma}^{G,M_{3}}\big)(\fq)|&=&|F_{\gamma,\,\mu}^{G,M_{3}}(\fq)|=q^{2n_{1}+n_{2}}+q^{n_{1}+n_{2}}(1+q+\cdots+q^{n_{1}-1})\nonumber\\
&&+q^{2n_{1}}(1+q+\cdots +q^{n_{2}-1}).
\end{eqnarray*}

With Arthur-Kottwitz reduction, we obtain

\begin{prop}\label{alpha AK}
\begin{eqnarray*}
|\big(\Lambda^{H_{M_{3}}}\backslash \xx_{\gamma}^{0}(\Pi_{N})\big)(\fq)|&=&|\big(\Lambda^{H_{M_{3}}}\backslash F_{\gamma}^{G,M_{3}}\big)(\fq)|+2Nq^{n_{1}+n_{2}}|\big(\Lambda^{H_{M_{3}}}\backslash F_{\gamma}^{M_{3},M_{3}}\big)(\fq)|\nonumber\\
&=& q^{2n_{1}+n_{2}}+q^{n_{1}+n_{2}}(1+q+\cdots+q^{n_{1}-1})\nonumber\\
&&+q^{2n_{1}}(1+q+\cdots +q^{n_{2}-1})+2Nq^{2n_{1}+n_{2}}\nonumber \\
&=&(2N+1)q^{2n_{1}+n_{2}}+q^{n_{1}+n_{2}}(1+q+\cdots+q^{n_{1}-1})\nonumber\\
&&+q^{2n_{1}}(1+q+\cdots +q^{n_{2}-1}).\end{eqnarray*}

\end{prop}

For the Harder-Narasimhan reduction, we count the points on the tail
\begin{eqnarray*}
|\big(\Lambda^{H_{M_{3}}}\backslash ^{t}\xx_{\gamma}^{0}(\Pi_{N})\big)(\fq)|&=&2|\big(\Lambda^{H_{M_{3}}}\backslash F_{\gamma}^{G,M_{3}}\big)(\fq)|\nonumber\\
&=& 2[q^{2n_{1}+n_{2}}+q^{n_{1}+n_{2}}(1+q+\cdots+q^{n_{1}-1})\nonumber\\
&&+q^{2n_{1}}(1+q+\cdots+q^{n_{2}-1})],
\end{eqnarray*}
and the $\xi$-stable points
\begin{equation*}
|\big(\Lambda^{H_{M_{3}}}\backslash \xx_{\gamma}^{M_{3},0,\,\xi^{M_{3}}}\big)(\fq)|=|\big(\Lambda^{H_{M_{3}}}\backslash \xx_{\gamma}^{M_{3},0}\big)(\fq)|=q^{n_{1}},
\end{equation*}
hence the points in the main body
\begin{eqnarray*}
|\big(\Lambda^{H_{M_{3}}}\backslash ^{m}\xx_{\gamma}^{0}(\Pi_{N})\big)(\fq)|&=&|\big(\Lambda^{H_{M_{3}}}\backslash \xx_{\gamma}^{\xi,0}\big)(\fq)|+[2N-(n_{1}+n_{2}+1)]\cdot q^{n_{1}+n_{2}}\cdot \nonumber\\
&&|\big(\Lambda^{H_{M_{3}}}\backslash \xx_{\gamma}^{M_{3},0,\,\xi^{M_{3}}}\big)(\fq)|\nonumber \\
&=&
|\big(\Lambda^{H_{M_{3}}}\backslash \xx_{\gamma}^{\xi,0}\big)(\fq)|+[2N-(n_{1}+n_{2}+1)]\cdot q^{2n_{1}+n_{2}}.
\end{eqnarray*}

Combining them, we get
\begin{prop}\label{alpha HN}
\begin{eqnarray*}
|\big(\Lambda^{H_{M_{3}}}\backslash \xx_{\gamma}^{0}(\Pi_{N})\big)(\fq)|&=&|\big(\Lambda^{H_{M_{3}}}\backslash ^{m}\xx_{\gamma}^{0}(\Pi_{N})\big)(\fq)|+|\big(\Lambda^{H_{M_{3}}}\backslash ^{t}\xx_{\gamma}^{0}(\Pi_{N})\big)(\fq)|\nonumber\\
&=&
|\big(\Lambda^{H_{M_{3}}}\backslash \xx_{\gamma}^{\xi,0}\big)(\fq)|+[2N-(n_{1}+n_{2}+1)]\cdot q^{2n_{1}+n_{2}}\nonumber\\
&&+2[q^{2n_{1}+n_{2}}+q^{n_{1}+n_{2}}(1+q+\cdots+q^{n_{1}-1})\nonumber\\
&&+q^{2n_{1}}(1+q+\cdots+q^{n_{2}-1})].
\end{eqnarray*}
\end{prop}

Comparing proposition \ref{alpha AK} and \ref{alpha HN}, we get
\begin{prop}\label{alpha stable}
\begin{eqnarray*}
|\big(\Lambda^{H_{M_{3}}}\backslash \xx_{\gamma}^{\xi,0}\big)(\fq)|&=&
(n_{1}+n_{2})q^{2n_{1}+n_{2}}-q^{n_{1}+n_{2}}(1+q+\cdots+q^{n_{1}-1})\nonumber\\
&&-q^{2n_{1}}(1+q+\cdots+q^{n_{2}-1}).
\end{eqnarray*}
\end{prop}

As before, the volume factor $\vol_{dt}\big(\Lambda^{H_{M_{3}}}\backslash T(F)_{M_{3}}^{1}\big)$ equals 1, and so

\begin{thm}\label{alpha}

We have 
\begin{eqnarray*}
J_{M_{3}}^{\xi}(\gamma)=|\big(\Lambda^{H_{M_{3}}}\backslash \xx_{\gamma}^{\xi,0}\big)(\fq)|&=&(n_{1}+n_{2})q^{2n_{1}+n_{2}}-q^{n_{1}+n_{2}}(1+q+\cdots+q^{n_{1}-1})\\
&&-q^{2n_{1}}(1+q+\cdots+q^{n_{2}-1}).
\end{eqnarray*}

\end{thm}

\section{Calculations for $\gl_{3}$--Mixed case}

Let $G=\gl_{3}$, $\gamma\in \ggl_{3}(F)$ a regular semi-simple integral element.
Assume that $T\cong F^{\times}\times \res_{E_{2}/F}E_{2}^{\times}$, with $E_{2}$ a separable totally ramified field extension over $F$ of degree $2$. As before, we can reduce to the case that $\gamma$ is a matrix of the form
\begin{equation}\label{gl3 mixed ramified}
\gamma=\begin{bmatrix}
a&&\\
&&b_{0}\ep^{n}\\
&b_{0}\ep^{n+1}&
\end{bmatrix},
\end{equation}
with $a\in \co$, $b_{0}\in \co^{\times}$. Let $m=\val(a)$.

Let $P$ be the parabolic subgroup $P=B_{0}\cup B_{0}s_{2}B_{0}$. Let $P=MN$ be the standard Levi decomposition. We identify $X_{*}(A)\cong \bz^{3}$ in the usual way. This gives us an identification $\Lambda_{M}\cong \bz^{2}$ and hence $\Lambda_{M}\otimes \br\cong \br^{2}$, we also identify $\ka_{M}^{G}$ with the line $x+y=0$ in $\br^{2}$, which can be further identified with $\br$ by taking the coordinate $x$. Under these identifications, the moment polytope $\Sigma_{\gamma}$ of the fundamental domain $F_{\gamma}$ can be taken to be the closed interval
$$
\Sigma_{\gamma}=[\,-n(\gamma, P,P^{-}),\,0\,]\subset \br\cong \ka_{M}^{G}.
$$ 
To simplify the notations, we abbreviate $n(\gamma, P,P^{-})$ to $n_{\gamma}$. We have 
$$
n_{\gamma}=\min\{2m, 2n+1\}.
$$ 
By definition, we can take
$$
F_{\gamma}=\xx_{\gamma}\cap \xx^{n+1}(\Sigma_{\gamma}).
$$
This can be refined a little bit. For $(n_{1},n_{2})\in \bn^{2}$, let $\Pi_{n_{1},n_{2}}$ be the positive $(G,A)$-orthogonal family as indicated in Fig. \ref{pi n} (The dashed part not included).

\begin{figure}[h]
\begin{center}
\begin{tikzpicture}[scale=0.6, every node/.style={scale=0.7}]

\draw [dashed] (-0.58,3)--(5.18,3);
\draw [dashed] (2.3,-2)--(5.18,3);


\draw (-0.58,3)--(-5.18,3);
\draw (-5.18,3)--(0,-6);
\draw (0,-6)--(2.3,-2);
\draw (2.3,-2)--(-0.58,3);

\node [blue] at (-0.58,3.5) {$(0,n_{2},0)$};
\node [black] at (-0.58,3) {$\bullet$};

\node [blue] at (-5.18,3.5) {$(-n_{1}, n_{1}+n_{2},0)$};
\node [black] at (-5.18,3) {$\bullet$};

\node [blue] at (0,-6.5) {$(-n_{1},0, n_{1}+n_{2})$};
\node [black] at (0,-6) {$\bullet$};

\node [blue] at (3.7,-2) {$(0,0,n_{2})$};
\node [black] at (2.3,-2) {$\bullet$};

\node [blue] at (5.18,3.5) {$(n_{2},0,0)$};
\node [black] at (5.18,3) {$\bullet$};







\end{tikzpicture}
\caption{The $(G,A)$-orthogonal family $\Pi_{n_{1},n_{2}}$ and its extension to a triangle.}
\label{pi n}
\end{center}
\end{figure}

Consider the positive $(G,A)$-orthogonal families $\Pi_{n_{\gamma},n+1}$. For $i\in \bz,\,-n_{\gamma}\leq i\leq 0$, let 
$$
\Pi_{n_{\gamma},n+1}^{i}=\Pi_{n_{\gamma},n+1}\cap \pi_{M}^{-1}(i),
$$ 
where $i\in \bz$ is considered as an element of $\ka_{M}^{G}$ by the identification $\br\cong \ka_{M}^{G}$. By theorem \ref{gl2 cal ramified}, we have 
$$
\xx_{\gamma}^{M,(i,n+1-i)}\subset \xx^{M,(i,n+1-i)}(\Pi_{n_{\gamma},n+1}^{i}),\quad \text{ for } i=-n_{\gamma},\cdots, 0.
$$
This implies that 
$$
F_{\gamma}=\xx_{\gamma}\cap \xx^{n+1}(\Sigma_{\gamma})=\xx_{\gamma}\cap \xx^{n+1}(\Pi_{n_{\gamma},n+1}).
$$

It is possible but quite hard to construct an affine paving of $F_{\gamma}$ and to count the number of rational points with it. Instead, we take an indirect way. Let $\Delta_{n_{\gamma},n+1}$ be the completion of $\Pi_{n_{\gamma},n+1}$ into a triangle as indicated in Fig. \ref{pi n}. We can count the number of rational points on $\xx_{\gamma}^{n+1}(\Delta_{n_{\gamma},n+1})$ quite easily, using the affine pavings in \cite{chen}, proposition 3.6. The complementary $\xx_{\gamma}^{n+1}(\Delta_{n_{\gamma},n+1})\backslash F_{\gamma}$ can be treated by the Arthur-Kottwitz reduction. Taking their difference, we find $|F_{\gamma}(\fq)|$.

We calculate the number of rational points on $\xx_{\gamma}^{n+1}(\Delta_{n_{\gamma},n+1})$. 
For $N\in \bn$, let 
$$
I_{N}=\Ad(\diag(\ep^{N},1,1)) I.
$$
According to \cite{chen}, proposition 3.6, when $N\gg 0$, we have an affine paving
$$
\xx_{\gamma}^{n+1}(\Delta_{n_{\gamma},n+1})=\bigcup_{\ep^{\baa}\in \xx_{\gamma}^{n+1}(\Delta_{n_{\gamma},n+1})^{A}} \xx_{\gamma}^{n+1}(\Delta_{n_{\gamma},n+1})\cap I_{N}\ep^{\baa}K/K.
$$
The dimension of the affine paving can be calculated using \cite{chen}, lemma 3.1, together with theorem \ref{gl2 cal ramified}. When $m\leq n$, i.e. $\gamma$ is not equivalued, the dimension of the paving is
$$
\min\{a_{2},m\}+\min\{a_{3},m\}+\begin{cases}
a_{2}-a_{3},& \text{ if } 0\leq a_{2}-a_{3}\leq n;\\
a_{3}-a_{2}-1, &\text{ if } 1\leq a_{3}-a_{2}\leq n+1
\end{cases}
$$
Otherwise, the intersection is empty.
When $m\geq n+1$, i.e. $\gamma$ is equivalued, the dimension of the paving is
$$
\min\{a_{2},n\}+\min\{a_{3},n+1\}+\begin{cases}
a_{2}-a_{3},& \text{ if } 0\leq a_{2}-a_{3}\leq n;\\
a_{3}-a_{2}-1, &\text{ if } 1\leq a_{3}-a_{2}\leq n+1.
\end{cases}
$$
Otherwise, the intersection is empty. We summarize the situation in the following two pictures, Fig. \ref{ram count 1} and Fig. \ref{ram count 2}. The triangle is cut into 4 parts by the two long red lines, the dimension of the fibration $f_{P}$ restricted to the affine pavings in different parts are given by different formulas. The two dashed lines bound the region where $\xx_{\gamma}^{M}\cap I_{N}\ep^{\baa}K/K$ is non-empty.

\begin{figure}[h]
\begin{center}
\begin{tikzpicture}[scale=0.6, every node/.style={scale=0.7}]



\draw (5.18,3)--(-5.18,3);
\draw (-5.18,3)--(0,-6);
\draw (0,-6)--(5.18,3);

\draw [dashed] (-1.15,-4)--(2.3,-2);
\draw [dashed] (-3.74,0.5)--(0.58,3);

\draw [red] (-4.03,1)--(4.03,1);
\draw [red] (-1.15,-4)--(2.88,3);

\draw [blue, densely dotted] (-1.73,1)--(0,-2);
\draw [blue, densely dotted] (4.6,2)--(4.03,3);
\draw [blue, densely dotted] (-1.15,2)--(2.3,2);
\draw [blue, densely dotted] (0.58,-3)--(2.88,1);

\node [blue] at (0.58,3.5) {$(1,n,0)$};
\node [black] at (0.58,3) {$\bullet$};

\node [blue] at (-5.18,3.5) {$(-2m, 2m+n+1,0)$};
\node [black] at (-5.18,3) {$\bullet$};

\node [blue] at (0,-6.5) {$(-2m,0, 2m+n+1)$};
\node [black] at (0,-6) {$\bullet$};

\node [blue] at (3.7,-2) {$(0,0,n+1)$};
\node [black] at (2.3,-2) {$\bullet$};

\node [blue] at (6.5,3.5) {$(n+1,0,0)$};
\node [black] at (5.18,3) {$\bullet$};

\node [red] at (-6.5,1) {$(-2m,m+n+1,m)$};
\node [red] at (-4.03,1) {$\bullet$};

\node [red] at (6,1) {$(n+1-m,0,m)$};
\node [red] at (4.03,1) {$\bullet$};

\node [red] at (3,3.5) {$(n+1-m,m,0)$};
\node [red] at (2.88,3) {$\bullet$};

\node [red] at (-3.5,-4) {$(-2m,m,m+n+1)$};
\node [red] at (-1.15,-4) {$\bullet$};







\end{tikzpicture}
\caption{Counting points in the non-equivalued case.}
\label{ram count 1}
\end{center}
\end{figure}

\begin{figure}[h]
\begin{center}
\begin{tikzpicture}[scale=0.5, every node/.style={scale=0.7}]

\draw [dashed] (-4.03,-0.33)--(2.88,3.67);
\draw [dashed] (-2.01,-3.83)--(4.03,-0.33);


\draw (6.33,3.67)--(-6.33,3.67);
\draw (-6.33,3.67)--(0,-7.33);
\draw (0,-7.33)--(6.33,3.67);


\draw [red] (-4.03,-0.33)--(4.03,-0.33);
\draw [red] (-1.73,-4.33)--(2.88,3.67);

\draw [blue, densely dotted] (-1.73,-0.33)--(-0.58,-2.33);
\draw [blue, densely dotted] (1.73,-0.33)--(0.58,-2.33);
\draw [blue, densely dotted] (4.03,3.67)--(5.18,1.67);
\draw [blue, densely dotted] (-0.58,1.67)--(1.73,1.67);

\node [blue] at (2.88,4.2) {$(1,n,0)$};
\node [black] at (2.88,3.67) {$\bullet$};

\node [blue] at (-6.33,4.2) {$(-(2n+1), 3n+2,0)$};
\node [black] at (-6.33,3.67) {$\bullet$};

\node [blue] at (0,-8) {$(-(2n+1),0, 3n+2)$};
\node [black] at (0,-7.33) {$\bullet$};

\node [blue] at (6,-0.33) {$(0,0,n+1)$};
\node [black] at (4.03,-0.33) {$\bullet$};

\node [blue] at (6.33,4.2) {$(n+1,0,0)$};
\node [black] at (6.33,3.67) {$\bullet$};

\node [red] at (-7.5,-0.33) {$(-(2n+1),2n+1,n+1)$};
\node [red] at (-4.03,-0.33) {$\bullet$};

\node [red] at (-5,-4.33) {$(-(2n+1),n,2(n+1))$};
\node [red] at (-1.73,-4.33) {$\bullet$};







\end{tikzpicture}
\caption{Counting points in the equivalued case.}
\label{ram count 2}
\end{center}
\end{figure}

\begin{prop}

Let $\gamma$ be a matrix in the form $(\ref{gl3 mixed ramified})$. 
When $\val(a)=m\leq n$, we have
\begin{eqnarray*}
|\xx_{\gamma}^{n+1}(\Delta_{2m,n+1})(\fq)|&=&
\sum_{j=0}^{2m-1}\left(\bigg\lfloor \frac{j}{2} \bigg\rfloor+1\right) q^{j}
+(2m+n+1)q^{2m}\\
&&+\sum_{j=2m+1}^{m+n}(4m+n+1-j)q^{j}
+\sum_{j=m+n+1}^{2m+n}\big(3(2m+n-j)+1\big)q^{j}\\
&&+q^{2m}\sum_{j=0}^{n}q^{j}.
\end{eqnarray*}

\end{prop}

In the summation, we use the convention that a summand is empty if its subscript is greater than its superscript.

\begin{proof}

Summing along the dotted blue lines in the 4 regions of Fig. \ref{ram count 1}, we get:
\begin{eqnarray*}
|\xx_{\gamma}^{n+1}(\Delta_{2m,n+1})(\fq)|&=&
\sum_{i=0}^{m}q^{i}(1+q+\cdots+q^{i})+
\sum_{i=m+1}^{2m}q^{i}(1+q+\cdots+q^{2m-i})\\
&&+\sum_{i=0}^{m-1}q^{i+m}(q^{m+1-i}+\cdots+q^{n})
+\sum_{i=0}^{m-1}q^{i+m}(q^{m-i}+\cdots+q^{n})\\
&&+ q^{2m}\sum_{i=2m+1}^{2m+n}(1+q+\cdots+q^{i-2m})+q^{2m}\sum_{i=0}^{n}q^{i}.
\end{eqnarray*}
After rearranging the summand, we get the proposition.

\end{proof}

\begin{prop}

Let $\gamma$ be a matrix in the form $(\ref{gl3 mixed ramified})$. 
When $\val(a)=m> n$, we have
\begin{eqnarray*}
|\xx_{\gamma}^{n+1}(\Delta_{2n+1,n+1})(\fq)|&=&
\sum_{j=0}^{2n}\left(\bigg\lfloor \frac{j}{2} \bigg\rfloor+1\right) q^{j}
+\sum_{j=2n+1}^{3n+1}\big(3(3n-j+1)+1\big)q^{j}\\
&&+q^{2n+1}\sum_{j=0}^{n}q^{j}.
\end{eqnarray*}

\end{prop}

\begin{proof}

Summing along the dotted blue lines in the 4 regions of Fig. \ref{ram count 2}, we get:
\begin{eqnarray*}
|\xx_{\gamma}^{n+1}(\Delta_{2n+1,n+1})(\fq)|&=&
\sum_{i=0}^{n}q^{i}(1+q+\cdots+q^{i})+
\sum_{i=n+1}^{2n+1}q^{i}(1+q+\cdots+q^{2n+1-i})\\
&&+\sum_{i=1}^{n}q^{i+n}(q^{n+1-i}+\cdots+q^{n})
+\sum_{i=1}^{n-1}q^{i+n+1}(q^{n+1-i}+\cdots+q^{n})\\
&&+ q^{2n+1}\sum_{i=2n+2}^{3n+1}(1+q+\cdots+q^{i-(2n+1)})+q^{2n+1}\sum_{i=0}^{n}q^{i}.
\end{eqnarray*}
After rearranging the summand, we get the proposition.

\end{proof}

Now we calculate the number of rational points on the complement $\xx_{\gamma}^{n+1}(\Delta_{n_{\gamma},n+1})\backslash F_{\gamma}$.
For $i\in \bz,\,1\leq i\leq n+1$, let 
$$
\Delta_{n_{\gamma},n+1}^{i}=\Delta_{n_{\gamma},n+1}\cap \pi_{M}^{-1}(i),
$$ 
where $i\in \bz$ is considered as an element of $\ka_{M}^{G}$ by the identification $\br\cong \ka_{M}^{G}$.

\begin{prop}

Let $\gamma$ be a matrix in the form $(\ref{gl3 mixed ramified})$. We have
$$
\xx_{\gamma}^{n+1}(\Delta_{n_{\gamma},n+1})\backslash F_{\gamma}
=
\bigcup_{i=1}^{n+1}f_{P}^{-1}\big(\xx_{\gamma}^{M,(i,n+1-i)}(\Delta^{i}_{n_{\gamma},n+1})\big)\cap \xx_{\gamma},
$$
where $(i,n+1-i)\in \bz^{2}$ is regarded as an element in $\Lambda_{M}$ by the identification $\bz^{2}\cong \Lambda_{M}$.
Its number of rational points over $\fq$ equals
$$
q^{n_{\gamma}}\sum_{j=0}^{n}(1+q+\cdots+q^{j}).
$$

\end{prop}

\begin{proof}

Observe that the second assertion is a direct consequence of the first one by proposition \ref{KL retraction}, it is hence enough to show the first assertion.

Let $x\in \xx_{\gamma}^{n+1}(\Delta_{n_{\gamma},n+1})$, notice that it doesn't belong to $F_{\gamma}$ if and only if 
\begin{equation}\label{outside}
H_{P}(x)\in [1,n+1]\subset \br\cong \ka_{M}^{G},
\end{equation}
because $H_{P^{-}}(x)\leq H_{P}(x)$. This implies that 
$$
\xx_{\gamma}^{n+1}(\Delta_{n_{\gamma},n+1})\backslash F_{\gamma}
=
\bigcup_{i=1}^{n+1}f_{P}^{-1}(\xx_{\gamma}^{M,(i,n+1-i)}(\Delta^{i}_{n_{\gamma},n+1}))\cap \xx_{\gamma}(\Delta_{n_{\gamma},n+1}).
$$
To finish the proof, we only need to show that 
\begin{equation*}
f_{P}^{-1}(\xx_{\gamma}^{M,(i,n+1-i)}(\Delta^{i}_{n_{\gamma},n+1}))\cap \xx_{\gamma}(\Delta_{n_{\gamma},n+1})=f_{P}^{-1}(\xx_{\gamma}^{M,(i,n+1-i)}(\Delta^{i}_{n_{\gamma},n+1}))\cap \xx_{\gamma},
\end{equation*} 
for $i=1,\cdots, n+1$. The inclusion ``$\subset$'' is obvious, we only need to show its inverse. For any point $x\in f_{P}^{-1}(\xx_{\gamma}^{M,(i,n+1-i)}(\Delta^{i}_{n_{\gamma},n+1}))\cap \xx_{\gamma}$, the inclusion (\ref{outside}) holds. By proposition \ref{gkm bound}, together with the fact that $\ec(x)$ is a positive $(G,A)$-orthogonal family, we have
$$
\ec(x)\subset \Delta_{n_{\gamma},n+1},
$$
whence the equality we want.

\end{proof}

Summarize all the above discussions, we get:

\begin{thm}\label{count mixed ram gl3}

Let $\gamma$ be a matrix in the form $(\ref{gl3 mixed ramified})$. When $\val(a)=m\leq n$, we have
\begin{eqnarray*}
|F_{\gamma}(\fq)|&=&
\sum_{j=0}^{2m-1}\left(\bigg\lfloor \frac{j}{2} \bigg\rfloor+1\right) q^{j}
+(2m+1)\sum_{j=2m}^{m+n} q^{j}\\
&&+\sum_{j=m+n+1}^{2m+n-1}\big(2(2m+n-j)+1\big)q^{j}
+q^{2m+n}.
\end{eqnarray*}

When $\val(a)=m> n$, we have
\begin{eqnarray*}
|F_{\gamma}(\fq)|&=&
\sum_{j=0}^{2n}\left(\bigg\lfloor \frac{j}{2} \bigg\rfloor+1\right) q^{j}
+\sum_{j=2n+1}^{3n}\big(2(3n+1-j)+1\big)q^{j}
+q^{3n+1}.
\end{eqnarray*}

\end{thm}

Now it is easy to deduce the weighted orbital integral $J_{M}^{\xi}(\gamma)$. For $N\in \bn,\,N\gg 0$, let 
$$
\Pi_{N}=[-n_{\gamma}-N, N]\subset \br\cong \ka_{M}^{G}.
$$
We can count the number of rational points $|\xx_{\gamma}^{n+1}(\Pi_{N})(\fq)|$ in two ways. By the Arthur-Kottwitz reduction, we have
$$
|\xx_{\gamma}^{n+1}(\Pi_{N})(\fq)|=|F_{\gamma}(\fq)|+2N q^{n_{\gamma}}\cdot |F_{\gamma}^{M}(\fq)|.
$$
By the Harder-Narasimhan reduction, we have
\begin{eqnarray*}
|\xx_{\gamma}^{n+1}(\Pi_{N})(\fq)|&=&2|F_{\gamma}(\fq)|+|\xx_{\gamma}^{n+1,\,\xi}(\fq)|\\
&&+(2N-n_{\gamma}-1) q^{n_{\gamma}}\cdot |F_{\gamma}^{M}(\fq)|,
\end{eqnarray*}
here we use the fact that $\xx_{\gamma}^{M,\,\nu,\,\xi}=F_{\gamma}^{M}$ for any $\nu\in \Lambda_{M}$ because $\gamma$ is anisotropic in $\km(F)$. 
The comparison of the two expressions implies
$$
|\xx_{\gamma}^{n+1,\xi}(\fq)|=(n_{\gamma}+1)q^{n_{\gamma}} \cdot |F_{\gamma}^{M}(\fq)|-|F_{\gamma}(\fq)|.
$$
By theorem \ref{count mixed ram gl3} and \ref{gl2 cal ramified}, we have

\begin{thm}

Let $\gamma$ be a matrix in the form $(\ref{gl3 mixed ramified})$. When $\val(a)=m\leq n$, we have
\begin{eqnarray*}
J_{M}^{\xi}(\gamma)=|\xx_{\gamma}^{n+1,\,\xi}(\fq)|
&=&
2mq^{2m+n}+\sum_{j=m+n+1}^{2m+n-1}2(j-m-n)q^{j}\\
&&
-\sum_{j=0}^{2m-1}\left(\bigg\lfloor \frac{j}{2} \bigg\rfloor+1\right) q^{j}.
\end{eqnarray*}

When $\val(a)=m> n$, we have
\begin{eqnarray*}
J_{M}^{\xi}(\gamma)=|\xx_{\gamma}^{n+1,\,\xi}(\fq)|
&=&
(2n+1)q^{3n+1}+\sum_{j=2n+1}^{3n}(2j-4n-1)q^{j}\\
&&
-\sum_{j=0}^{2n}\left(\bigg\lfloor \frac{j}{2} \bigg\rfloor+1\right) q^{j}.
\end{eqnarray*}

\end{thm}

By theorem \ref{cl comparision} and the remark following it, we get Arthur's weighted orbital integral.
As before, the orbital integral $I^{G}_{\gamma}$ can be calculated by equation (\ref{geom orbital reduced}):
$$
I_{\gamma}^{G}=q^{n_{\gamma}}\sum_{i=0}^{n}q^{i}=\begin{cases}
q^{2m}\sum_{i=0}^{n}q^{i},&\text{ if }m\leq n;\\
q^{2n+1}\sum_{i=0}^{n}q^{i},&\text{ if }m>n.
\end{cases}
$$

\section{Calculations for $\gl_{3}$--Anisotropic case}

Let $G=\gl_{3}$, $\gamma\in \ggl_{3}(F)$ a regular semisimple integral element. Assume that $\mathrm{char}(k)>3$ and $T\cong \res_{E_{3}/F}E_{3}^{\times}$, with  $E_{3}=\fq(\!(\ep^{\frac{1}{3}})\!)$. As before, take the basis $\{\ep^{\frac{2}{3}}, \ep^{\frac{1}{3}},1\}$ of $E_{3}$ over $F$, we can assume that $\gamma$ is of the form
\begin{equation}\label{gamma sl3 ramified}
\gamma=\begin{bmatrix}
&b_{0}\ep^{n_{1}}&c_{0}\ep^{n_{2}}\\
c_{0}\ep^{n_{2}+1} &&b_{0}\ep^{n_{1}}\\
b_{0}\ep^{n_{1}+1}&c_{0}\ep^{n_{2}+1}&
\end{bmatrix},
\end{equation}
with $b_{0},c_{0}\in \co^{\times}$ and $n_{1},n_{2}\in \bn$.
In this case, Arthur's weighted orbital integral is the same as the orbital integral, and both are equal to $|\xx_{\gamma}^{0}(\fq)|$. The matrix $\gamma$ is equivalued of valuation $n_{1}+\frac{1}{3}$ if $n_{1}\leq n_{2}$, and equivalued of valuation $n_{2}+\frac{2}{3}$ if $n_{2}< n_{1}$. 
According to Goresky, Kottwitz and MacPherson \cite{gkm2}, the affine Springer fiber $\xx_{\gamma}$ admits affine paving
$$
\xx_{\gamma}=\bigcup_{\baa=(a_{1},a_{2},a_{3})\in \bz^{3}} \xx_{\gamma} \cap I\ep^{\baa}K/K.
$$
Let $S_{\baa}$ be the cell $\xx_{\gamma} \cap I\ep^{\baa}K/K$. Restricted to the connected component $\xx_{\gamma}^{0}$, we can calculate that $S_{\baa}$ is non-empty if and only if 
\begin{equation}\label{non-empty domain}
a_{1}-a_{2}\leq n_{1},\quad a_{2}-a_{3}\leq n_{1},\quad a_{3}-a_{1}\leq n_{1}+1,
\end{equation}
and that it is of dimension
$$
|\{(m,\alpha)\in \bz\times \Phi(G,A)\mid 0\leq m+\alpha(x)<n_{1}+\frac{1}{3},\, m+\alpha(y_{\baa})<0\}|,
$$
with $x=(1,2/3,1/3), y_{\baa}=(-a_{1},-a_{2},-a_{3})\in X_{*}(A)\otimes \br$. The results are summarized in Fig. \ref{gl3 ram elliptic}. Summing up, we get

\begin{figure}[t]
\begin{center}
\begin{tikzpicture}[scale=0.6, every node/.style={scale=0.7}]

\draw [dashed] (-6.93,-4)--(6.93,4);
\draw [dashed] (0,-8)--(0,8);
\draw [dashed] (-6.93,4)--(6.93,-4);

\node [blue] at (0,0,0){$(0,0,0)$};

\draw (5.2,-9)--(5.2,9);
\draw (-10.4,0)--(5.2,9);
\draw (-10.4,0)--(5.2,-9);
\draw (-10.4,0)--(-9.8,-1);
\draw (4.05,-9)--(-9.8,-1);
\draw (4.05,-9)--(5.2,-9);

\draw (0,6)--(5.2,3);
\draw (5.2,-3)--(0,-6);
\draw (-5.2,-3)--(-5.2,3);

\node [blue] at (5.2, 10) {$(n_{1},0,-n_{1})$};
\node [black] at (5.2, 9) {$\bullet$};

\node [blue] at (7.2, -9) {$(0,-n_{1},n_{1})$};
\node [black] at (5.2, -9) {$\bullet$};

\node [blue] at (4.05, -10) {$(-1,1-n_{1},n_{1})$};
\node [black] at (4.05, -9) {$\bullet$};

\node [blue] at (-12.5,0) {$(-n_{1},n_{1},0)$};
\node [black] at (-10.4,0) {$\bullet$};

\node [blue] at (-11.5,-1.5) {$(-n_{1},n_{1}-1,1)$};
\node [black] at (-9.8,-1) {$\bullet$};







\node [blue] at (2.5,3.5) {$q^{2(a_{1}-a_{3})}$};
\node [black] at (1.73, 3) {$\bullet$};

\node [blue] at (-2.5,-3.5) {$q^{2(a_{3}-a_{1})-3}$};
\node [black] at (-1.73, -3) {$\bullet$};

\node [blue] at (-2.5,3.5) {$q^{2(a_{2}-a_{3})-1}$};
\node [black] at (-1.73, 3) {$\bullet$};

\node [blue] at (2.5,-3.5) {$q^{2(a_{3}-a_{2})-2}$};
\node [black] at (1.73, -3) {$\bullet$};

\node [blue] at (3,0.5) {$q^{2(a_{1}-a_{2})-1}$};
\node [black] at (3.45, 0) {$\bullet$};

\node [blue] at (-3.5,0.5) {$q^{2(a_{2}-a_{1})-2}$};
\node [black] at (-3.45, 0) {$\bullet$};

\node [blue] at (3.46,6.5) {$q^{a_{1}-a_{3}+n_{1}}$};
\node [black] at (3.46, 6) {$\bullet$};

\node [blue] at (3.46,-6.5) {$q^{a_{3}-a_{2}+n_{1}-1}$};
\node [black] at (3.46, -6) {$\bullet$};

\node [blue] at (-6.9,0.5) {$q^{a_{2}-a_{1}+n_{1}-1}$};
\node [black] at (-6.9,0) {$\bullet$};

\node [blue] at (-9,-2.5) {$q^{a_{2}-a_{1}+n_{1}-1}$};
\node [black] at (-8.05,-2) {$\bullet$};

\node [blue] at (-4,-5.5) {$q^{2(a_{3}-a_{1})-3}$};
\node [black] at (-2.88,-5) {$\bullet$};

\node [blue] at (1.3, -8.5) {$q^{a_{3}-a_{2}+n_{1}-1}$};
\node [black] at (2.3,-8) {$\bullet$};

\end{tikzpicture}
\caption{Counting points for ramified anisotropic $\gamma\in \ggl_{3}(\co)$-first case.}
\label{gl3 ram elliptic}
\end{center}
\end{figure}

\begin{thm}\label{gl3 ram count 1}

Let $\gamma\in \ggl_{3}(\co)$ be the matrix in the form $(\ref{gamma sl3 ramified})$. Suppose that $n_{1}\leq n_{2}$, it is then equivalued of valuation $n_{1}+\frac{1}{3}$. The orbital integral associated to $\gamma$ equals

\begin{eqnarray*}
I_{\gamma}^{G}&=&|\xx_{\gamma}^{0}(\fq)|=1+2\sum_{i=1}^{\lfloor \frac{n_{1}}{3}\rfloor} q^{2(3i-1)}(q^{2}+q+1)\\
&&+\sum_{i=1}^{n_{1}} \left(i-2 \bigg\lfloor \frac{i}{3} \bigg\rfloor -1 \right) q^{2i-3}(q^{3}+2q^{2}+2q+1)\\
&&+\sum_{i=n_{1}+1}^{2n_{1}}\left(2n_{1}-i-2 \bigg\lceil \frac{2n_{1}-i}{3} \bigg\rceil+1\right)q^{i+n_{1}-1}(q+2)\\
&&
+q^{2n_{1}-1}\left(\bigg\lceil \frac{2n_{1}-1}{3} \bigg\rceil-\bigg\lfloor \frac{n_{1}-2}{3} \bigg\rfloor -1\right) 
+2\sum_{i=\lceil \frac{2n_{1}-1}{3} \rceil}^{n_{1}-1}q^{3i+1}.
\end{eqnarray*}
where $\lfloor x\rfloor$ denotes the largest integer less than or equal to $x$, and $\lceil x\rceil$ denotes the smallest integer greater than or equal to $x$.

\end{thm}

Same calculations for $n_{2}< n_{1}$, with the difference that $S_{\baa}$ is non-empty if and only if 
$$
a_{1}-a_{3}\leq n_{2},\quad a_{2}-a_{1}\leq n_{1}+1,\quad a_{3}-a_{2}\leq n_{1}+1,
$$
and that it is isomorphic to an affine space of dimension
$$
|\{(m,\alpha)\in \bz\times \Phi(G,A)\mid 0\leq m+\alpha(x)<n_{2}+\frac{2}{3},\, m+\alpha(y_{\baa})<0\}|.
$$ 
These are summarized schematically in Fig. \ref{gl3 ram elliptic 2}. Summing up, we get

\begin{figure}[b]
\begin{center}
\begin{tikzpicture}[scale=0.6, every node/.style={scale=0.7}]

\draw [dashed] (-6.93,-4)--(6.93,4);
\draw [dashed] (0,-8)--(0,8);
\draw [dashed] (-6.93,4)--(6.93,-4);

\node [blue] at (0,0,0){$(0,0,0)$};

\draw (-5.75,8)--(-5.75,-10);
\draw (-5.75,-10)--(9.78,-1);
\draw (9.78,-1)--(10.35,0);
\draw (-5.18,9)--(10.35,0);
\draw (-5.18,9)--(10.35,0);
\draw (-5.18,9)--(-5.75,8);
\draw (-5.18,9)--(-5.18,-9);
\draw (-5.18,-9)--(10.35,0);

\draw (0,6)--(-5.2,3);
\draw (-5.2,-3)--(0,-6);
\draw (5.2,-3)--(5.2,3);

\draw (-0.58,6.33)--(-5.18,3.67);
\draw (-5.18,-3.67)--(-0.58,-6.33);
\draw (5.75,2.67)--(5.75,-2.67);

\node [blue] at (-5, 9.5) {$(0,n_{2},-n_{2})$};
\node [black] at (-5.18, 9) {$\bullet$};

\node [blue] at (-8, 8) {$(-1,n_{2},-n_{2}+1)$};
\node [black] at (-5.75, 8) {$\bullet$};

\node [blue] at (-5.75,-10.5) {$(-n_{2}-1,0,n_{2}+1)$};
\node [black] at (-5.75,-10) {$\bullet$};

\node [blue] at (-5,-9.5) {$(-n_{2},0,n_{2})$};
\node [black] at (-5.18,-9) {$\bullet$};

\node [blue] at (12,-1) {$(n_{2}-1,-n_{2},1)$};
\node [black] at (9.78,-1) {$\bullet$};

\node [blue] at (12,0) {$(n_{2},-n_{2},0)$};
\node [black] at (10.35,0) {$\bullet$};







\node [blue] at (2.5,3.5) {$q^{2(a_{1}-a_{3})}$};
\node [black] at (1.73, 3) {$\bullet$};

\node [blue] at (-2.5,-3.5) {$q^{2(a_{3}-a_{1})-3}$};
\node [black] at (-1.73, -3) {$\bullet$};

\node [blue] at (-2.5,3.5) {$q^{2(a_{2}-a_{3})-1}$};
\node [black] at (-1.73, 3) {$\bullet$};

\node [blue] at (2.5,-3.5) {$q^{2(a_{3}-a_{2})-2}$};
\node [black] at (1.73, -3) {$\bullet$};

\node [blue] at (3,0.5) {$q^{2(a_{1}-a_{2})-1}$};
\node [black] at (3.45, 0) {$\bullet$};

\node [blue] at (-3.5,0.5) {$q^{2(a_{2}-a_{1})-2}$};
\node [black] at (-3.45, 0) {$\bullet$};

\node [blue] at (8,0.5) {$q^{a_{1}-a_{2}+n_{2}}$};
\node [black] at (8.05,0) {$\bullet$};

\node [blue] at (-3.8,6.5) {$q^{a_{2}-a_{3}+n_{2}}$};
\node [black] at (-4.03,7) {$\bullet$};

\node [blue] at (-3.5,-6.5) {$q^{a_{3}-a_{1}+n_{2}-1}$};
\node [black] at (-4.03,-7) {$\bullet$};

\node [blue] at (-1,-8.5) {$q^{a_{3}-a_{1}+n_{2}-1}$};
\node [black] at (-2.3,-8) {$\bullet$};

\node [blue] at (4,-5) {$q^{2n_{2}}$};
\node [black] at (2.88,-5) {$\bullet$};

\node [blue] at (9.5,-2) {$q^{a_{1}-a_{2}+n_{2}}$};
\node [black] at (8.05,-2) {$\bullet$};

\end{tikzpicture}
\caption{Counting points for ramified anisotropic $\gamma\in \ggl_{3}(\co)$-second case.}
\label{gl3 ram elliptic 2}
\end{center}
\end{figure}

\begin{thm}\label{gl3 ram count 2}

Let $\gamma\in \ggl_{3}(\co)$ be the matrix in the form $(\ref{gamma sl3 ramified})$. Suppose that $n_{2}< n_{1}$, it is then equivalued of valuation $n_{2}+\frac{2}{3}$. The orbital integral associated to $\gamma$ equals
\begin{eqnarray*}
I_{\gamma}^{G}&=&|\xx_{\gamma}^{0}(\fq)|=1+2\sum_{i=1}^{\lfloor \frac{n_{2}}{3}\rfloor} q^{2(3i-1)}(q^{2}+q+1)\\
&&+\sum_{i=1}^{n_{2}} \left(i-2 \bigg\lfloor \frac{i}{3} \bigg\rfloor -1 \right) q^{2i-3}(q^{3}+2q^{2}+2q+1)\\
&&+\left(n_{2}-2 \bigg\lceil \frac{n_{2}-1}{3} \bigg\rceil \right) q^{2n_{2}-1}(1+2q^{2})\\
&&+\sum_{i=n_{2}+2}^{2n_{2}}\left(2n_{2}-i-2 \bigg\lceil \frac{2n_{2}-i}{3} \bigg\rceil+1\right)q^{i+n_{2}-1}(1+2q)\\
&&+2\sum_{i=0}^{\lfloor \frac{n_{2}-2}{3} \rfloor} q^{3(n_{2}-i)-2}(1+q)
+2q^{2n_{2}}\left(\bigg\lceil \frac{2n_{2}-1}{3} \bigg\rceil-\bigg\lfloor \frac{n_{2}-2}{3} \bigg\rfloor -1\right) 
+q^{3n_{2}+1}.
\end{eqnarray*}

\end{thm}


\clearpage


































\begin{thebibliography}{100}
\labelwidth=4em
\addtolength\leftskip{25pt}
\setlength\labelsep{0pt}
\addtolength\parskip{\smallskipamount}



\bibitem[A1]{a}{J. Arthur, \textit{The characters of discrete series as orbital integrals}, Invent. math. \textbf{32}(1976), 205-261.}  

\bibitem[A2]{a4}{J. Arthur, \textit{The trace formula in invariant form}, Ann. of Math., 114 (1981), 1-74.}


\bibitem[A3]{a2}{J. Arthur, \textit{A local trace formula}, Inst. Hautes \'Etudes Sci. Publ. Math. No. 73 (1991), 5–96.}


\bibitem[A4]{a3}{J. Arthur, \textit{An introduction to the trace formula}. Harmonic analysis, the trace formula, and Shimura varieties, 1–263, Clay Math. Proc., 4, Amer. Math. Soc., Providence, RI, 2005.} 


\bibitem[B]{b}{R. Bezrukavnikov, \textit{The dimension of the fixed point set on affine flag manifolds}, Math. Res. Lett. 3 (1996), 185–189.}


\bibitem[Bo]{bo}{M. Borovoi, \textit{Abelian Galois cohomology of reductive groups}. Mem. Amer. Math. Soc. 132 (1998), no. 626, viii+50 pp.}









\bibitem[C1]{chen}{Z. Chen, \textit{Pureté des fibres de Springer affines pour $\mathrm{GL}_{4}$}, Bull. SMF 142, fascicule 2 (2014), 193-222. }


\bibitem[C2]{chen1}{Z. Chen, \textit{The $\xi$-stability on the affine grassmannian}. Math. Z. 280 (2015), no. 3-4, 1163–1184.}



\bibitem[C3]{chen2}{Z. Chen, \textit{On the fundamental domain of affine Springer fibers}. Math. Z. 286 (2017), no. 3-4, 1323–1356.}




\bibitem[C4]{chen3}{Z. Chen, \textit{Truncated affine grassmannians and truncated affine Springer fibers for $\gl_{3}$}, \url{http://arxiv.org/abs/1401.1930.}}



\bibitem[CL1]{cl1}{P-H. Chaudouard, G. Laumon, \textit{Sur l'homologie des fibres de Springer affines tronquées}, Duke Math. J. 145 (2008), no. 3, 443–535.}


\bibitem[CL2]{cl2}{P-H. Chaudouard, G. Laumon, \textit{Le lemme fondamental pondéré. I. Constructions géométriques.} Compos. Math. 146 (2010), no. 6, 1416–1506.} 



\bibitem[GKM1]{gkm1}{M. Goresky, R. Kottwitz, R. MacPherson, \textit{Homology of affine Springer fibers in the unramified case}, Duke Math. J. 121 (2004), no. 3, 509-561.}

\bibitem[GKM2]{gkm2}{M. Goresky, R. Kottwitz, R. MacPherson, \textit{Purity of equivalued affine Springer fibers}, Represent. Theory 10 (2006), 130-146.}

\bibitem[GKM3]{gkm3}{M. Goresky, R. Kottwitz, R. MacPherson, \textit{Regular points in affine Springer fibers}, Michigan Math. J. 53 (2005), no. 1, 97-107.}


\bibitem[GKM4]{gkm local}{M. Goresky, R. Kottwitz, R. MacPherson, \emph{Equivariant cohomology, Koszul duality, and the localization theorem}. Invent. Math. 131 (1998), no. 1, 25–83. }







\bibitem[K1]{k1}{R. E. Kottwitz, \textit{Isocrystals with additional structure}. Compositio Math. 56 (1985), no. 2, 201–220.}


\bibitem[K2]{k2}{R. E. Kottwitz, \textit{Isocrystals with additional structure. II}. Compositio Math. 109 (1997), no. 3, 255–339.}




\bibitem[KL]{kl}{D. Kazhdan, G. Lusztig, \textit{Fixed point varieties on affine flag manifolds}, Israel. J. Math. \textbf{62}(1988), 129-168.}









\bibitem[MP]{mp}{A. Moy, G. Prasad, \textit{Unrefined minimal K-types for p-adic groups}, Invent. Math. 116 (1994), 393–408.}




\bibitem[N]{ngo}{Bau Ch\^{a}u Ng\^o, \textit{Le lemme fondamental pour les algèbres de Lie.} Publ. Math. IHES. No. 111 (2010), 1–169.}










\end{thebibliography}
\end{document}